\documentclass[11pt,final]{article}%
\usepackage[para,ragged]{footmisc}
\usepackage{amsfonts}
\usepackage{amsmath}
\usepackage{amssymb}
\usepackage{amscd}
\usepackage[amsmath,hyperref,thmmarks]{ntheorem}
\usepackage{makeidx}
\usepackage{graphicx}%
\usepackage[OT2,T1]{fontenc}
\newcommand\textcyr[1]{{\fontencoding{OT2}\fontfamily{wncyr}\selectfont #1}}
\newcommand{\TS}{{\mbox{\textcyr{Sh}}}}
\newcommand{\TB}{{\mbox{\textcyr{B}}}}

\let\oldmarginpar\marginpar
\renewcommand\marginpar[1]{\-\oldmarginpar[\raggedleft\footnotesize #1]%
{\raggedright\footnotesize #1}}

\usepackage[textwidth=5.5in,textheight=9.0in,centering]{geometry}

\usepackage{mathptmx}

\usepackage{array}

\usepackage{microtype}

\usepackage{amsxtra,amscd}
\usepackage{tikz}
\usetikzlibrary{matrix,arrows,positioning,decorations.pathmorphing}
 
\usepackage{enumitem}
\setitemize[1]{label=$\diamond$\hspace{0.07in}}
\setlist{nolistsep}

\usepackage{etex}

\usepackage{url}
\usepackage{verbatim}

\usepackage{mathtools}

\usepackage[overload]{textcase}

\usepackage{titlesec}                                                                            
\titleformat*{\section}{\LARGE\bfseries}
\titleformat*{\subsection}{\Large\itshape}
\titleformat*{\subsubsection}{\scshape}
\titleformat*{\paragraph}{\itshape}

\usepackage{natbib}
\let\cite\citealt

\newcommand{\bcomment}{}
\newcommand{\bfootnotesize}{\begin{footnotesize}}\newcommand\efootnotesize{\end{footnotesize}}
\newcommand{\bquote}{\begin{quote}}\newcommand\equote{\end{quote}}
\newcommand{\bsmall}{\begin{small}}\newcommand\esmall{\end{small}}
\newcommand{\btable}{\begin{table}}\newcommand{\etable}{\end{table}}
\newcommand{\dstyle}{\displaystyle}
\newcommand{\edocument}{\end{document}}

\newcommand{\tableoc}{\tableofcontents}

\def\1{{1\mkern-7mu1}}

\DeclareMathOperator{\Aut}{Aut}

\DeclareMathOperator{\Br}{Br}

\DeclareMathOperator{\Coker}{Coker}

\DeclareMathOperator{\End}{End}

\DeclareMathOperator{\Ext}{Ext}

\DeclareMathOperator{\Gal}{Gal}
\DeclareMathOperator{\GL}{GL}

\DeclareMathOperator{\Hom}{Hom}

\DeclareMathOperator{\im}{Im}  
\DeclareMathOperator{\Ind}{Ind}

\DeclareMathOperator{\inv}{inv}

\DeclareMathOperator{\Ker}{Ker}

\DeclareMathOperator{\Nm}{Nm}

\DeclareMathOperator{\ord}{ord}

\DeclareMathOperator{\Pic}{\mathrm{Pic}}

\DeclareMathOperator{\rank}{rank}
\DeclareMathOperator{\rec}{rec}
\DeclareMathOperator{\Res}{Res}
\DeclareMathOperator{\res}{res}

\DeclareMathOperator{\Spec}{Spec}

\DeclareMathOperator{\Tgt}{Tgt}

\DeclareMathOperator{\Tr}{Tr}

\theoremnumbering{arabic}
\theoremheaderfont{\scshape}
\RequirePackage{latexsym}
\theorembodyfont{\slshape}
\theoremseparator{.}
\newtheorem{X}{X}[subsection]

\newtheorem{E}[X]{}

\theorembodyfont{\upshape}

\theorembodyfont{\small}

\newtheorem*{nt}{Notes}
\theorembodyfont{\normalsize}
\theoremstyle{nonumberplain}
\theoremsymbol{\ensuremath{_\Box}}

\qedsymbol{\ensuremath{_\Box}}
\makeindex
\begin{document}

\title{The Work of John Tate}
\date{March 18, 2012/September 23, 2012}
\author{J.S. Milne}
\maketitle

\hfill\begin{minipage}{3.5in}
\textit{Tate helped shape the great reformulation of arithmetic and geometry which has taken
place since the 1950s}\\
\hspace*{2in}Andrew Wiles.\footnotemark
\end{minipage}\footnotetext{Introduction to Tate's talk at the conference on
the Millenium Prizes, 2000.}\bigskip

\begin{quote}
{\small This is my article on Tate's work for the second volume in the book
series on the Abel Prize winners. True to the epigraph, I have attempted to
explain it in the context of the \textquotedblleft great
reformulation\textquotedblright.}
\end{quote}

\tableoc

\subsection*{Notations}

We speak of the primes of a global field where others speak of the places.

$M_{S}=S\otimes_{R}M$ for $M$ an $R$-module and $S$ and $R$-algebra.

$|S|$ is the cardinality of $S$.

$X_{n}=\Ker(x\mapsto nx\colon X\rightarrow X)$ and $X(\ell)=\bigcup_{m\geq
0}X_{\ell^{m}}$ ($\ell$-primary component, $\ell$ a prime).

$\Gal(K/k)$ or $G(K/k)$ denotes the Galois group of $K/k$.

$\mu(R)$ is the group of roots of $1$ in $R$.

$R^{\times}$ denotes the group of invertible elements of a ring $R$ (apologies
to Bourbaki).

$V^{\vee}$ denotes the dual of a vector space or the contragredient of a representation.

$K^{\mathrm{al}}$, $K^{\mathrm{sep}}$, $K^{\mathrm{ab}}$,$K^{\mathrm{un}}%
$\ldots\ denote an algebraic, separable, abelian, unramified$\ldots$closure of
a field $K$.

$\mathcal{O}{}_{K}$ denotes the ring of integers in a local or global field
$K$.

\section{Hecke $L$-series and the cohomology of number
fields\label{a}}

\subsection{Background\label{a1}}

\subsubsection{Kronecker, Weber, Hilbert, and ray class groups}

For every abelian extension $L$ of $\mathbb{Q}{}$, there is an integer $m$
such that $L$ is contained in the cyclotomic field $\mathbb{Q}{}[\zeta_{m}]$;
it follows that the abelian extensions of $\mathbb{Q}{}$ are classified by the
subgroups of the groups $\mathbb{(}\mathbb{Z}{}/m\mathbb{Z}{})^{\times}\simeq
G(\mathbb{Q}{}[\zeta_{m}]/\mathbb{Q}{})$ (Kronecker-Weber). On the other hand,
the unramified abelian extensions of a number field $K$ are classified by the
subgroups of the ideal class group $C$ of $K$ (Hilbert). In order to be able
to state a common generalization of these two results, Weber introduced the
ray class groups. A modulus $\mathfrak{m}{}$ for a number field $K$ is the
formal product of an ideal $\mathfrak{m}{}_{0}$ in $\mathcal{O}{}_{K}$ with a
certain number of real primes of $K$. The corresponding ray class group
$C_{\mathfrak{m}{}}$ is the quotient of the group of ideals relatively prime
to $\mathfrak{m}{}_{0}$ by the principal ideals generated by elements
congruent to $1$ modulo $\mathfrak{m}{}_{0}$ and positive at the real primes
dividing $\mathfrak{m}{}$. For $\mathfrak{m}{}=(m)\infty$ and $K=\mathbb{Q}{}%
$, $C_{\mathfrak{m}{}}\simeq(\mathbb{Z}{}/m\mathbb{Z}{})^{\times}$. For
$\mathfrak{m}{}=1$, $C_{\mathfrak{m}{}}=C$.

\subsubsection{Takagi and the classification of abelian extensions}

Let $K$ be a number field. Takagi showed that the abelian extensions of $K$
are classified by the ray class groups: for each modulus $\mathfrak{m}{}$,
there is a well-defined \textquotedblleft ray class field\textquotedblright%
\ $L_{\mathfrak{m}{}}$ with $G(L_{\mathfrak{m}{}}/K)\approx C_{\mathfrak{m}{}%
}$, and every abelian extension of $K$ is contained in a ray class field for
some modulus $\mathfrak{m}{}$. Takagi also proved precise decomposition rules
for the primes in an extension $L/K$ in terms of the associated ray class
group. These would follow from knowing that the map sending a prime ideal to
its Frobenius element gives an isomorphism $C_{\mathfrak{m}{}}\rightarrow
G(L_{\mathfrak{m}{}}/K)$, but Takagi didn't prove that.

\subsubsection{Dirichlet, Hecke, and $L$-series}

For a character $\chi$ of $(\mathbb{Z}{}/m\mathbb{Z}{})^{\times}$, Dirichlet
introduced the $L$-series%
\[
L(s,\chi)=\prod\nolimits_{(p,m)=1}\frac{1}{1-\chi(p)p^{-s}}=\sum
\nolimits_{(n,m)=1}\chi(n)n^{-s}%
\]
in order to prove that each arithmetic progression, $a$, $a+m$, $a+2m$,
$\ldots$ with $a$ relatively prime to $m$ has infinitely many primes. When
$\chi$ is the trivial character, $L(s,\chi)$ differs from the zeta function
$\zeta(s)$ by a finite number of factors, and so has a pole at $s=1$.
Otherwise $L(s,\chi)$ can be continued to a holomorphic function on the entire
complex plane and satisfies a functional equation relating $L(s,\chi)$ and
$L(1-s,\bar{\chi})$.

Hecke proved that the $L$-series of characters of the ray class groups
$C_{\mathfrak{m}{}}$ had similar properties to Dirichlet $L$-series, and noted
that his methods apply to the $L$-series of even more general characters, now
called Hecke characters (Hecke 1918, 1920).\footnote{Hecke, E., Eine neue Art
von Zetafunktionen und ihre Beziehungen zur Verteilung der Primzahlen. I,
Math. Zeit. 1, 357-376 (1918); II, Math. Zeit. 6, 11-51 (1920).} The
$L$-series of Hecke characters are of fundamental importance. For example,
Deuring (1953)\footnote{Deuring, Max, Die Zetafunktion einer algebraischen
Kurve vom Geschlechte Eins. Nachr. Akad. Wiss. G\"{o}ttingen. 1953, 85--94.}
showed that the $L$-series of an elliptic curve with complex multiplication is
a product of two Hecke $L$-series.

\subsubsection{Artin and the reciprocity law}

Let $K/k$ be an abelian extension of number fields, corresponding to a
subgroup $H$ of a ray class group $C_{\mathfrak{m}{}}$. Then%
\begin{equation}
\zeta_{K}(s)/\zeta_{k}(s)=\prod\nolimits_{\chi}L(s,\chi)\quad\text{(up to a
finite number of factors)} \label{e42}%
\end{equation}
where $\chi$ runs through the nontrivial characters of $C_{\mathfrak{m}{}}/H$.
From this and the results of Dirichlet and Hecke, it follows that $\zeta
_{K}(s)/\zeta_{k}(s)$ is holomorphic on the entire complex plane. In the hope
of extending this statement to nonabelian extensions $K/k$, Artin
(1923)\footnote{Artin, E. \"{U}ber die Zetafunktionen gewisser algebraischer
Zahlk\"{o}rper. Math. Ann. 89 (1923), no. 1-2, 147--156.} introduced what are
now called Artin $L$-series.

Let $K/k$ be a Galois extension of number fields with Galois group $G$, and
let $\rho\colon G\rightarrow\GL(V)$ be a representation of $G$ on a finite
dimensional complex vector space $V$. The Artin $L$-series of $\rho$ is%
\[
L(s,\rho)=\prod\nolimits_{\mathfrak{p}{}}\frac{1}{\det(1-\rho(\sigma
_{\mathfrak{p}{}})\mathbb{N}{}\mathfrak{P}^{-s}\mid V^{I_{\mathfrak{P}{}}})}%
\]
where $\mathfrak{p}{}$ runs through the prime ideals of $K$, $\mathfrak{P}{}$
is a prime ideal of $K$ lying over $\mathfrak{p}{}$, $\sigma_{\mathfrak{p}{}}$
is the Frobenius element of $\mathfrak{P}{}$, $\mathbb{N}{}\mathfrak{P}%
{}=(\mathcal{O}{}_{K}\colon\mathfrak{P}{})$, and $I_{\mathfrak{P}{}}$ is the
inertia group.

Artin observed that his $L$-series for one-dimensional representations would
coincide with the $L$-series of characters on ray class groups if the
following \textquotedblleft theorem\textquotedblright\ were true:

\begin{quote}
for the field $L$ corresponding to a subgroup $H$ of a ray class group
$C_{\mathfrak{m}{}}$, the map $\mathfrak{p}{}\mapsto(\mathfrak{p}{},L/K)$
sending a prime ideal $\mathfrak{p}{}$ not dividing $\mathfrak{m}{}$ to its
Frobenius element induces an isomorphism $C_{\mathfrak{m}{}}/H\rightarrow
G(L/K)$.
\end{quote}

\noindent Initially, Artin was able to prove this statement only for certain
extensions. After Chebotarev had proved his density theorem by a reduction to
the cyclotomic case, Artin (1927)\footnote{Artin, E., Beweis des allgemeinen
Reziprozit\"{a}tsgesetzes. Abhandlungen Hamburg 5, 353-363 (1927).} proved the
statement\ in general. He called it the reciprocity law because, when $K$
contains a primitive $m$th root of $1$, it directly implies the classical
$m$th power reciprocity law.

Artin noted that $L(s,\rho)$ can be analytically continued to a meromorphic
function on the whole complex plane if its character $\chi$ can be expressed
in the form%
\begin{equation}
\chi=\sum\nolimits_{i}n_{i}\Ind\chi_{i},\quad n_{i}\in\mathbb{Z}{},
\label{e26}%
\end{equation}
with the $\chi_{i}$ one-dimensional characters on subgroups of $G$, because
then%
\[
L(s,\rho)=\prod\nolimits_{i}L(s,\chi_{i})^{n_{i}}%
\]
with the $L(s,\chi_{i})$ abelian $L$-series. Brauer (1947)\label{Brauer1947}%
\footnote{Brauer, Richard, On Artin's L-series with general group characters.
Ann. of Math. (2) 48, (1947). 502--514.} proved that the character of a
representation can always be expressed in the form (1), and Brauer and Tate
found what is probably the simplest known proof of this fact (see
p.\pageref{1955b}).

To complete his program, Artin conjectured that, for every nontrivial
irreducible representation $\rho$, $L(s,\rho)$ is holomorphic on the entire
complex plane. \noindent This is called the Artin conjecture. It is known to
be true if the character of $\rho$ can be expressed in the form (\ref{e26})
with $n_{i}\geq0{}$, and in a few other cases.

\subsubsection{Chevalley and id\`{e}les}

Chevalley gave a purely local proof of local class field theory, and a purely
algebraic proof of global class field theory, but probably his most lasting
contribution was to reformulate class field theory in terms of id\`{e}les.

An id\`{e}le of a number field $K$ is an element $(a_{v})_{v}$ of
$\prod\nolimits_{v}K_{v}^{\times}$ such that $a_{v}\in\mathcal{O}{}%
_{v}^{\times}$ for all but finitely many primes $v$. The id\`{e}les form a
group $J_{K}$, which becomes a locally compact topological group when endowed
with the topology for which the subgroup%
\[
\prod\nolimits_{v|\infty}K_{v}^{\times}\times\prod\nolimits_{v\text{ finite}%
}\mathcal{O}{}_{v}^{\times}%
\]
is open and has the product topology.\footnote{The original topology defined
by Chevalley is not Hausdorff. It was Weil who pointed out the need for a
topology in which the Hecke characters become the characters on $J{}$ (Weil,
A., Remarques sur des r\'{e}sultats recents de C. Chevalley. C. R. Acad. Sci.,
Paris 203, 1208-1210 1936). By the time of Tate's thesis, the correct
definition seems to have been common knowledge.}

Let $K$ be number field. In Chevalley's reinterpretation, global class field
theory provides a homomorphism $\phi\colon J_{K}/K^{\times}\rightarrow
G(K^{\mathrm{ab}}/K)$ that induces an isomorphism%
\[
J_{K}/(K^{\times}\cdot\Nm J_{L})\overset{}{\longrightarrow}G(L/K)
\]
for each finite abelian extension $L/K$. For each prime $v$ of $K$, local
class field theory provides a homomorphism $\phi_{v}\colon K_{v}^{\times
}\rightarrow G(K_{v}^{\mathrm{ab}}/K_{v})$ that induces an isomorphism%
\[
K_{v}^{\times}/\Nm L^{\times}\rightarrow G(L/K_{v})
\]
for each finite abelian extension $L/K_{v}$. The maps $\phi_{v}$ and $\phi$
are related by the diagram:%
\[
\begin{CD}
K_{v}^{\times} @>{\phi_{v}}>> G(K_{v}^{\mathrm{ab}}/K_{v})\\
@VV{i_{v}}V@VVV\\
J_K@>{\phi}>>G(K^{\mathrm{ab}}/K).
\end{CD}
\]
Beyond allowing class field theory to be stated for infinite extensions,
Chevalley's id\`{e}lic approach greatly clarified the relation between the
local and global reciprocity maps.

\subsection{Tate's thesis and the local constants\label{a2}}

The modern definition is that a Hecke character is a quasicharacter of
$J/K^{\times}$, i.e., a continuous homomorphism $\chi\colon J\rightarrow
\mathbb{C}{}^{\times}$ such that $\chi(x)=1$ for all $x\in K^{\times}$. We
explain how to interpret $\chi$ as a map on a group of ideals, which is the
classical definition.

For a finite set $S$ of primes, including the infinite primes, let $J^{S}$
denote the subgroup of $J_{K}{}$ consisting of the id\`{e}les $(a_{v})_{v}$
with $a_{v}=1$ for all $v\in S$, and let $I^{S}$ denote the group of
fractional ideals generated by those prime ideals not in $S$. There is a
canonical surjection $J^{S}\rightarrow I^{S}$. For each Hecke character $\chi
$, there exists a finite set $S$ such that $\chi$ factors through
$J^{S}\overset{\text{can.}}{\longrightarrow}I^{S}$, and a homomorphism
$\varphi\colon I^{S}\rightarrow\mathbb{C}{}^{\times}$ arises from a Hecke
character if and only if there exists an integral ideal $\mathfrak{m}{}$ with
support in $S$, complex numbers $(s_{\sigma})_{\sigma\in\Hom(K,\mathbb{C}{})}%
$, and integers $(m_{\sigma})_{\sigma\in\Hom(K,\mathbb{C}{})}$ such that%
\[
\varphi((\alpha))=\prod_{\sigma\in\Hom(K,\mathbb{C}{})}\sigma(\alpha
)^{m_{\sigma}}\left\vert \sigma(\alpha)\right\vert ^{s_{\sigma}}%
\]
for all $\alpha\in K^{\times}$ with $(\alpha)\in I^{S}$ and $\alpha
\equiv1\,\,$(mod $\mathfrak{m}{})$.

\subsubsection{Hecke's proof}

The classical proof uses that $\mathbb{R}^{n}{}$ is self-dual as an additive
topological group, and that the discrete subgroup $\mathbb{Z}{}^{n}$ of
$\mathbb{R}{}^{n}$ is its own orthogonal complement under the duality. The
Poisson summation formula follows easily\footnote{Let $f$ be a Schwartz
function on $\mathbb{R}{}$, and let $\hat{f}$ be its Fourier transform on
$\mathbb{\hat{R}}=\mathbb{R}{}{}$. Let $\phi$ be the function $x+\mathbb{Z}%
{}\mapsto\sum_{n\in\mathbb{Z}{}}f(x+n)$ on $\mathbb{R}{}/\mathbb{Z}$, and let
$\hat{\phi}$ be its Fourier transform on $\widehat{\mathbb{R}{}/\mathbb{Z}{}%
}=\mathbb{Z}{}$. A direct computation shows that $\hat{f}(n)=\hat{\phi}(n)$
for all $n\in\mathbb{Z}{}$. The Fourier inversion formula says that
$\phi(x)=\sum_{n\in\mathbb{Z}{}}\hat{\phi}(n)\chi(x)$; in particular,
$\phi(0)=\sum_{n\in\mathbb{Z}{}}\hat{\phi}(n)=\sum_{n\in\mathbb{Z}{}}\hat
{f}(n)$. But, by definition, $\phi(0)=\sum_{n\in\mathbb{Z}{}}f(n)$.} from
this: for any Schwartz function $f$ on $\mathbb{R}{}^{n}$ and its Fourier
transform $\hat{f}$,%
\[
\sum\nolimits_{m\in\mathbb{Z}{}^{n}}f(m)=\sum\nolimits_{m\in\mathbb{Z}{}^{n}%
}\hat{f}(m)\text{.}%
\]
Write the $L$-series as a sum over integral ideals, and decompose it into a
finite family of sums, each of which is over the integral ideals in a fixed
element of an ideal class group. The individual series are Mellin transforms
of theta series, and the functional equation follows from the transformation
properties of the theta series, which, in turn, follow from the Poisson
summation formula.

\subsubsection{Tate's proof}

An ad\`{e}le of $K$ is an element $(a_{v})_{v}$ of $\prod K_{v}$ such that
$a_{v}\in\mathcal{O}{}_{v}$ for all but finitely many primes $v$. The
ad\`{e}les form a ring $\mathbb{A}{}$, which becomes a locally compact
topological ring when endowed with its natural topology.

Tate proved that the ring of ad\`{e}les $\mathbb{A}{}$ of $K$ is self-dual as
an additive topological group, and that the discrete subgroup $K$ of
$\mathbb{A}{}$ is its own orthogonal complement under the duality. As in the
classical case, this implies an (ad\`{e}lic) Poisson summation formula: for
any Schwartz function $f$ on $\mathbb{A}{}$ and its Fourier transform $\hat
{f}$%
\[
\sum\nolimits_{\gamma\in K}f(\gamma)=\sum\nolimits_{\gamma\in K}\hat{f}%
(\gamma){}.
\]

Let $\chi$ be a Hecke character of $K$, and let $\chi_{v}$ be the
quasicharacter $\chi\circ i_{v}$ on $K_{v}^{\times}$. Tate defines local
$L$-functions $L(\chi_{v})$ for each prime $v$ of $K$ (including the infinite
primes) as integrals over $K_{v}$, and proves functional equations for them.
He writes the global $L$-function as an integral over $J$, which then
naturally decomposes into a product of local $L$-functions. The functional
equation for the global $L$-function follows from the functional equations of
the local $L$-functions and the Poisson summation formula.

Although, the two proofs are superficially similar, in the details they are
quite different. Once Tate has developed the harmonic analysis of the local
fields and of the ad\`{e}le ring, including the Poisson summation formula,
\textquotedblleft an analytic continuation can be given at one stroke for all
of the generalized $\zeta$-functions, and an elegant functional equation can
be established for them \ldots\ without Hecke's complicated
theta-formulas\textquotedblright.\footnote{Tate 1967, pp. 305--306.}

One consequence of Tate's treating all primes equally, is that the $\Gamma
$-factors arise naturally as the local zeta functions of the infinite primes.
By contrast, in the classical treatment, their appearance is more mysterious.

As Kudla writes:\footnote{Kudla, S. In: An introduction to the Langlands
program. Edited by Joseph Bernstein and Stephen Gelbart. Birkh\"{a}user
Boston, Inc., Boston, MA, 2003, p.133.}

\begin{quote}
Tate provides an elegant and unified treatment of the analytic continuation
and functional equation of Hecke $L$-functions. The power of the methods of
abelian harmonic analysis in the setting of Chevalley's ad\`{e}les/id\`{e}les
provided a remarkable advance over the classical techniques used by Hecke.
\ldots\ In hindsight, Tate's work may be viewed as giving the theory of
automorphic representations and $L$-functions of the simplest connected
reductive group $G=\GL(1)$, and so it remains a fundamental reference and
starting point for anyone interested in the modern theory of automorphic representations.
\end{quote}

\noindent Tate's thesis completed the re-expression of the classical theory in
terms of id\`{e}les. In this way, it marked the end of one era, and the start
of a new.

\begin{nt}
Tate completed his thesis in May 1950. It was widely quoted\ long before its
publication in 1967. Iwasawa obtained similar results about the same time as
Tate, but published nothing except for the brief notes Iwasawa 1950,
1952.\footnote{Iwasawa, K., A Note on Functions, Proceedings of the
International Congress of Mathematicians, Cambridge, Mass., 1950, vol. 1,
p.322. Amer. Math. Soc., Providence, R. I., 1952; Letter to J. Dieudonn\'{e}.
Zeta functions in geometry, April 8, 1952, Adv. Stud. Pure Math., 21,
pp.445--450, Kinokuniya, Tokyo, 1992.}
\end{nt}

\subsubsection{Local constants}

Let $\chi$ be a Hecke character, and let $\Lambda(s,\chi)$ be its completed
$L$-series. The theorem of Hecke and Tate says that $\Lambda(s,\chi)$ admits a
meromorphic continuation to the whole complex plane, and satisfies a
functional equation%
\[
\Lambda(1-s,\chi)=W(\chi)\cdot\Lambda(s,\bar{\chi})
\]
with $W(\chi)$ a complex number of absolute value $1$. The number $W(\chi)$ is
called the root number or the epsilon factor. It is a very interesting number.
For example, for a Dirichlet character $\chi$ with conductor $f$, it equals
$\tau(\chi)/\sqrt{\pm f}$ where $\tau(\chi)$ is the Gauss sum $\sum_{a=1}%
^{f}\chi(a)\mathbf{e}(a/f)$. An importance consequence of Tate's description
of the global functional equation as a product of local functional equations
is that he obtains an expression%
\begin{equation}
W(\chi)=\prod\nolimits_{v}W(\chi_{v}) \label{e27}%
\end{equation}
of $W(\chi)$ as a product of (explicit) local root numbers $W(\chi_{v})$.

Langlands pointed out\footnote{See his \textquotedblleft Notes on Artin
$L$-functions\textquotedblright\ and the associated comments at
\url{http://publications.ias.edu/rpl/section/22}.} that his conjectural
correspondence between degree $n$ representations of the Galois groups of
number fields and automorphic representations of $\GL_{n}$ requires that there
be a similar decomposition for the root numbers of Artin $L$-series, or, more
generally, for the Artin-Hecke $L$-series that generalize both Artin and Hecke
$L$-series (see p.\pageref{ArtinHecke}). For a Hecke character, the required
decomposition is just that of Tate. Every expression (\ref{e26}),
p.\pageref{e26}, of an Artin character $\chi$ as a sum of monomial characters
gives a decomposition of its root number $W(\chi)$ as a product of local root
numbers --- the problem is to show that the decomposition is
\textit{independent} of the expression of $\chi$ as a sum.\footnote{In fact,
this is not quite true, but is true for \textquotedblleft virtual
representations\textquotedblright\ with \textquotedblleft virtual degree
$0$\textquotedblright. The decomposition of the root number of the character
$\chi$ of a Galois representation is obtained by writing it as $\chi
=(\chi-\dim\chi\cdot1)+\dim\chi\cdot1.$}

For an Artin character $\chi$, Dwork (1956)\footnote{Dwork, B., On the Artin
root number. Amer. J. Math. 78 (1956), 444--472. Based on his 1954 thesis as a
student of Tate.} proved that there exists a decomposition (\ref{e27}) of
$W(\chi)$ well-defined up to signs; more precisely, he proved that there
exists a well-defined decomposition for $\chi(-1)W(\chi)^{2}$. Langlands
completed Dwork's work and thereby found a local proof that there exists a
well-defined decomposition for $W(\chi)$. However, he abandoned the writing up
of his proof when Deligne (1973)\footnote{Deligne, P. Les constantes des
\'{e}quations fonctionnelles des fonctions $L$. Modular functions of one
variable, II (Proc. Internat. Summer School, Univ. Antwerp, Antwerp, 1972),
pp. 501--597. Lecture Notes in Math., Vol. 349, Springer, Berlin, 1973.} found
a simpler global proof.

Tate (1977b) gives an elegant exposition of these questions, including a proof
of (\ref{e27}) for Artin root numbers by a variant of Deligne's method, and a
proof of a theorem of Fr\"{o}hlich and Queyrut that $W(\chi)=1$ when $\chi$ is
the character of a representation that preserves a quadratic form.

\subsection{The cohomology of number fields\label{a3}}

\subsubsection{Tate cohomology}

With the action of a group $G$ on an abelian group $M$, there are associated
homology groups $H_{r}(G,M)$, $r\geq0$, and cohomology groups $H^{r}(G,M)$,
$r\geq0$. When $G$ is finite, the map $m\mapsto\sum_{\sigma\in G}\sigma m$
defines a homomorphism%
\[
H_{0}(G,M)\overset{\textup{{\tiny def}}}{=}M_{G}\overset{\mathrm{Nm}_{G}%
}{\longrightarrow}M^{G}\overset{\textup{{\tiny def}}}{=}H^{0}(G,M)\text{,}%
\]
and Tate defined cohomology groups $\hat{H}^{r}(G,M)$ for all integers $r$ by
setting%
\[
\renewcommand{\arraystretch}{1.3}\hat{H}^{r}(G,M)\overset{\textup{{\tiny def}%
}}{=}\left\{
\begin{array}
[c]{ll}%
H_{-r-1}(G,M) & r<-1\\
\Ker(\Nm_{G}) & r=-1\\
\Coker(\Nm_{G}) & r=0\\
H^{r}(G,M) & r>0.
\end{array}
\right.
\]
The diagram 
\[
\begin{tikzpicture}[>=angle 90]
\matrix(m)[matrix of math nodes,
row sep=2em, column sep=1em,,
text height=1.5ex, text depth=0.25ex]
{&&\hat{H}^{-1}(G,M')&\hat{H}^{-1}(G,M)&\hat{H}^{-1}(G,M'')\\
\cdots&H_1(G,M'')&H_0(G,M')&H_0(G,M)&H_0(G,M'')&0\\
&0&H^0(G,M')&H^0(G,M)&H^0(G,M'')&H^1(G,M')&\cdots\\
&&\hat{H}^0(G,M')&\hat{H}^0(G,M)&\hat{H}^0(G,M'')\\};
\path[->,font=\scriptsize] 
          (m-1-3) edge (m-2-3)
          (m-1-4) edge (m-2-4)
          (m-1-5) edge (m-2-5)
          (m-3-3) edge (m-4-3)
          (m-3-4) edge (m-4-4)
          (m-3-5) edge (m-4-5)
          (m-2-2) edge (m-2-3)
          (m-2-3) edge (m-2-4)
          (m-2-4) edge (m-2-5)
          (m-2-5) edge (m-2-6)
          (m-2-3) edge node[right]{$\mathrm{Nm}_G$} (m-3-3)
          (m-2-4) edge node[right]{$\mathrm{Nm}_G$} (m-3-4)
          (m-2-5) edge node[right]{$\mathrm{Nm}_G$} (m-3-5)
          (m-2-3) edge (m-3-3)
          (m-2-3) edge (m-3-3)
          (m-3-2) edge (m-3-3)
          (m-3-3) edge (m-3-4)
          (m-3-4) edge (m-3-5)
          (m-3-5) edge (m-3-6);
\path[->,font=\scriptsize,thick]
          (m-1-3) edge (m-1-4)
          (m-1-4) edge (m-1-5)
          (m-2-1) edge (m-2-2)
          (m-3-6) edge (m-3-7)
          (m-4-3) edge (m-4-4)
          (m-4-4) edge (m-4-5)
          (m-4-5) edge (m-3-6)
          (m-2-2) edge (m-1-3)
          (m-1-5) edge[out=0,in=180] (m-4-3);
\end{tikzpicture}
\]
shows that a short exact sequence of
$G$-modules gives an exact sequence of cohomology groups infinite in both
directions. These groups are now called the Tate cohomology groups%
\index{Tate cohomology groups}%
. Most of the usual constructions for cohomology groups (except the inflation
maps) extend to the Tate groups.

\begin{nt}
Tate's construction was included in Serre 1953,\footnote{Serre, Jean-Pierre,
Cohomologie et arithm\'{e}tique, S\'{e}minaire Bourbaki 1952/1953, no. 77.}
Cartan and Eilenberg 1956,\footnote{Cartan, Henri; Eilenberg, Samuel.
Homological algebra. Princeton University Press, Princeton, N. J., 1956} and
elsewhere. Farrell 1978\footnote{Farrell, F. Thomas, An extension of Tate
cohomology to a class of infinite groups. J. Pure Appl. Algebra 10 (1977/78),
no. 2, 153-161.} extended Tate's construction to infinite groups having finite
virtual cohomological dimension (Tate-Farrell cohomology%
\index{Tate-Farrell cohomology}%
), and others have defined an analogous extension of Hochschild cohomology
(Tate-Hochschild cohomology%
\index{Tate-Hochschild cohomology}%
).
\end{nt}

\subsubsection{The cohomology groups of algebraic number fields}

Let $G$ be a finite group, let $C$ be a $G$-module, and let $u$ be an element
of $H^{2}(G,C)$. Assume that $H^{1}(H,C)=0$ for all subgroups $H$ of $G$ and
that $H^{2}(H,C)$ is cyclic of order $(H\colon1)$ with generator the
restriction of $u$. Then Tate (1952c) showed that cup product with $u$ defines
an isomorphism%
\begin{equation}
x\mapsto x\cup u\colon\hat{H}^{r}(G,\mathbb{Z}{})\rightarrow\hat{H}^{r+2}(G,C)
\label{e32}%
\end{equation}
for all $r\in\mathbb{Z}{}$. He proves this by constructing an exact sequence%
\[
0\rightarrow C\rightarrow C(\varphi)\rightarrow\mathbb{Z}{}[G]\rightarrow
\mathbb{Z}{}\rightarrow0,
\]
depending on the choice of a $2$-cocycle $\varphi$ representing $u$, and
showing that
\[
\hat{H}^{r}(G,C(\varphi))=0=\hat{H}^{r}(G,\mathbb{Z}{}[G])
\]
for all $r\in\mathbb{Z}{}$. Now the double boundary map is an isomorphism
$\hat{H}^{r}(G,\mathbb{Z}{})\rightarrow\hat{H}^{r+2}(G,C)$.

On taking $G$ to be the Galois group of a finite extension $L/K$ of number
fields, $C$ to be the id\`{e}le class group of $L$, and $u$ the fundamental
class of $L/K$,\footnote{Which had been discovered by Nakayama and Weil, cf.
Artin and Tate, 2009, p189.} one obtains for $r=-2$ the inverse of the Artin
reciprocity map%
\[
G/[G,G]\overset{\simeq}{\longrightarrow}C^{G}\text{.}%
\]

Let $L/K$ be a finite Galois extension of global fields (e.g., number fields)
with Galois group $G$. There is an exact sequence of $G$-modules%
\begin{equation}
1\rightarrow L^{\times}\rightarrow J_{L}\rightarrow C_{L}\rightarrow1
\label{e34}%
\end{equation}
where $J_{L}$ is the group of id\`{e}les of $L$ and $C_{L}$ is the id\`{e}le
class group. Tate determined the cohomology groups of the terms in this
sequence by relating them to those in the much simpler sequence%
\begin{equation}
0\rightarrow X\rightarrow Y\rightarrow Z{}\rightarrow0. \label{e31}%
\end{equation}
Here $Y$ is the free abelian group on the set of primes of $L$ (including the
infinite primes) with $G$ acting through its action on the primes, and $Z$ is
just $\mathbb{Z}{}$ with $G$ acting trivially; the map $Y\rightarrow Z$ is
$\sum n_{P}P\mapsto\sum n_{P}$, and $X$ is its kernel. Tate proved that there
is a canonical isomorphism of doubly infinite exact sequences%
\begin{equation}
\minCDarrowwidth20pt\begin{CD} \cdots@>>> \hat{H}^{r}(G,X)@>>> \hat{H}^{r}(G,Y)@>>> \hat{H}^{r}(G,Z)@>>>\cdots\\ @.@VV{\simeq}V@VV{\simeq}V@VV{\simeq}V@.\\ \cdots@>>> \hat{H}^{r+2}(G,L^{\times})@>>> \hat{H}^{r+2}(G,J_L)@>>> \hat{H}^{r+2}(G,C_L)@>>>\cdots. \end{CD} \label{e33}%
\end{equation}
Tate announced this result in his Short Lecture at the 1954 International
Congress, but did not immediately publish the proof.

\subsubsection{The Tate-Nakayama theorem}

Nakayama (1957) generalized Tate's isomorphism (\ref{e32}) by weakening the
hypotheses --- it suffices to require them for Sylow subgroups --- and
strengthening the conclusion --- cup product with $u$ defines an isomorphism%
\[
x\mapsto x\cup u\colon\hat{H}^{r}(G,M{})\rightarrow\hat{H}^{r+2}(G,C\otimes
M)
\]
provided $C$ or $M$ is torsion-free. Building on this, Tate (1966c) proved
that the isomorphism (\ref{e33}) holds with each of the sequences (\ref{e34})
and (\ref{e31}) replaced by its tensor product with $M$. In other words, he
replaced the torus $\mathbb{G}_{m}$ implicit in (\ref{e33}) with an arbitrary
torus defined over $K$. He also proved the result for any \textquotedblleft
suitably large\textquotedblright\ set of primes $S$ --- the module $L^{\times
}$ is replaced with the group of $S$-units in $L$ and $J_{L}$ is replaced by
the group of id\`{e}les whose components are units outside $S$. This result is
usually referred to as the Tate-Nakayama theorem%
\index{Tate-Nakayama theorem}%
, and is widely used, for example, throughout the Langlands program including
in the proof of the fundamental lemma.

\subsubsection{Abstract class field theory: class formations}

Tate's theorem (see (\ref{e32}) above) shows that, in order to have a class
field theory over a field $k$, all one needs is, for each system of fields%
\[
k^{\mathrm{sep}}\supset L\supset K\supset k,\quad\lbrack L\colon
k]<\infty,\quad L/K\text{ Galois,}%
\]
a $G(L/K)$-module $C_{L}$ and a \textquotedblleft fundamental
class\textquotedblright\ $u_{L/K}\in H^{2}(G(L/K),C_{L})$ satisfying Tate's
hypotheses; the pairs $(C_{L},u_{L/K})$ should also satisfy certain natural
conditions when $K$ and $L$ vary. Then Tate's theorem then provides
\textquotedblleft reciprocity\textquotedblright\ isomorphisms%
\[
C_{L}^{G}\overset{\simeq}{\longrightarrow}G/[G,G],\quad G=G(L/K),
\]
Artin and Tate (1961, Chapter 14) formalized this by introducing the abstract
notion of a class formation.

For example, for any nonarchimedean local field $k$, there is a class
formation with $C_{L}=L^{\times}$ for any finite extension $L$ of $k$, and for
any global field, there is a class formation with $C_{L}=J_{L}/L^{\times}$. In
both cases, $u_{L/K}$ is the fundamental class.

Let $k$ be an algebraic function field in one variable with algebraically
closed constant field. Kawada and Tate (1955a) show that there is a class
formation for unramified extensions of $K$ with $C_{L}$ the dual of the group
of divisor classes of $L$. In this way they obtain a \textquotedblleft pseudo
class field theory\textquotedblright\ for $k$, which they examine in some
detail when $k=\mathbb{C}{}$.

\subsubsection{The Weil group}

Weil was the first to find a common generalization of Artin $L$-series and
Hecke $L$-series. For this he defined what is now known as the Weil group. The
Weil group of a finite Galois extension of number fields $L/K$ is an extension%
\[
1\rightarrow C_{L}\rightarrow W_{L/K}\rightarrow\Gal(L/K)\rightarrow1
\]
corresponding to the fundamental class in $H^{2}(G_{L/K},C_{L})$. Each
representation of $W_{L/K}$ has an $L$-series attached to it, and the
$L$-series arising in this way are called\label{ArtinHecke} Artin-Hecke
$L$-series. Weil (1951)\footnote{Weil, Andr\'{e}, Sur la th\'{e}orie du corps
de classes. J. Math. Soc. Japan 3, (1951). 1--35.} constructed these groups,
thereby discovering the fundamental class, and proved the fundamental
properties of Weil groups. Artin and Tate (1961, Chapter XV) developed the
theory of the Weil groups in the abstract setting of class formations, basing
their definition on the existence of a fundamental class. In the latest (2009)
edition of the work, Tate expanded their presentation and included a sketch of
Weil's original construction (pp. 185--189).

\subsubsection{Summary}

The first published exposition of class field theory in which full use of the
cohomology theory is made is Chevalley 1954.\footnote{Chevalley, C. Class
field theory. Nagoya University, Nagoya, 1954.} There Chevalley writes:

\begin{quote}
One of the most baffling features of classical class field theory was that it
appeared to say practically nothing about normal extensions that are not
abelian. It was discovered by A. Weil and, from a different point of view, T.
Nakayama that class field theory was actually much richer than hitherto
suspected; in fact, it can now be formulated in the form of statements about
normal extensions without any mention whatsoever of abelian extensions. Of
course, it is true that it is only in the abelian case that these statements
lead to laws of decomposition for prime ideals of the subfield and to the law
of reciprocity. Nevertheless, it is clear that, by now, we know something
about the arithmetic of non abelian extensions. In fact, since the work of J.
Tate, it may be said that we know almost everything that may be formulated in
terms of cohomology in the id\`{e}le class group, and generally a great deal
about everything that can be formulated in cohomological terms.
\end{quote}

\begin{nt}
Tate was not the first to make use of group cohomology in class field theory.
In a sense it had always been there, since crossed homomorphisms and factor
systems had long been used. Weil and Nakayama independently discovered the
fundamental class, Weil by constructing the Weil group, and Nakayama as a
consequence of his work (partly with Hochschild) to determine the cohomology
groups of number fields in degrees $1$ and $2$. Tate's contribution was to
give a remarkably simple description of all the basic cohomology groups of
number fields, and to construct a general isomorphism that, in the particular
case of an abelian extension and in degree $-2$, became the Artin reciprocity isomorphism.
\end{nt}

\subsection{The cohomology of profinite groups\label{a4}}

Krull (1928)\footnote{Krull, W., Galoissche Theorie der unendlichen
algebraischen Erweiterungen. Math. Ann. 100, 687-698 (1928).} showed that,
when the Galois group of an infinite Galois extension of fields $\Omega/F$ is
endowed with a natural topology, there is a Galois correspondence between the
intermediate fields of $\Omega/F$ and the \textit{closed} subgroups of the
Galois group. The topological groups that arise as Galois groups are exactly
the compact groups $G$ whose open normal subgroups $U$ form a fundamental
system $\mathcal{N}{}$ of neighbourhoods of $1$. Tate described such
topological groups as being \textquotedblleft of Galois-type\textquotedblright%
, but we now say they are \textquotedblleft profinite\textquotedblright.

For a profinite group $G$, Tate (1958d) considered the $G$-modules $M$ such
that $M=\bigcup\nolimits_{U\in\mathcal{N}{}}M^{U}$. These are the $G$-modules
$M$ for which the action is continuous relative to the discrete topology on
$M$. For such a module, Tate defined cohomology groups $H^{r}(G,M)$, $r\geq0$,
using continuous cochains, and he showed that%
\[
H^{r}(G,M)=\varinjlim\nolimits_{U\in\mathcal{N}{}}H^{r}(G/U,M^{U})
\]
where $H^{r}(G/U,M^{U})$ denotes the usual cohomology of the (discrete) finite
group $G/U$ acting on the abelian group $M^{U}$. In particular, $H^{r}(G,M)$
is torsion for $r>0$.

The cohomological dimension and strict cohomological dimension of a profinite
group $G$ relative to a prime number $p$ are defined by the conditions:%
\begin{align*}
\mathrm{cd}_{p}(G)  &  \leq n\iff H^{r}(G,M)(p)=0\text{ whenever }r>n\text{
and }M\text{ is torsion;}\\
\mathrm{scd}_{p}(G)  &  \leq n\iff H^{r}(G,M)(p)=0\text{ whenever }r>n\text{.}%
\end{align*}
Here $(p)$ denotes the $p$-primary component. The (strict) cohomological
dimension of a field is the (strict) cohomological dimension of its absolute
Galois group. Among Tate's theorems are the following statements:

\begin{enumerate}
\item A pro $p$-group $G$ is free if and only if $\mathrm{cd}_{p}(G)=1$. (A
pro $p$-group is a profinite group $G$ such that $G/U$ is a $p$-group for all
$U\in\mathcal{N}{}$; it is free if it is of the form $\varprojlim F/N$ where
$F$ is the free group on symbols $(a_{i})_{i\in I}$, say, and $N$ runs through
the normal subgroups of $G$ containing all but finitely many of the $a_{i}$
and such $G/N$ is a finite $p$-group).

\item If $k$ is a local field other than $\mathbb{R}{}$ or $\mathbb{C}{}$,
then $\mathrm{scd}_{p}(k)=2$ for all $p\neq\mathrm{char}(k)$.

\item Let $K\supset k$ be an extension of fields of transcendence degree $n$.
Then%
\[
\mathrm{cd}_{p}(K)\leq\mathrm{cd}_{p}(k)+n\text{,}%
\]
with equality if $K$ is finitely generated over $k$, $\mathrm{cd}%
_{p}(k)<\infty$, and $p\neq\mathrm{char}(k)$. In particular, if $k$ is
algebraically closed, then the $p$-cohomological dimension of a finitely
generated $K$ is equal to its transcendence degree over $k$ ($p\neq
\mathrm{char}(k)$).
\end{enumerate}

\noindent According to Tate 1958d, statement (c) \textquotedblleft
historically arose at [the theory's] beginning. Its conjecture and the sketch
of its proof are due to Grothendieck\textquotedblright. Indeed, from
Grothendieck's point of view, the cohomology of the absolute Galois group of a
field $k$ should be interpreted as the \'{e}tale cohomology of $\Spec k$, and
the last statement of (c) is suggested by the weak Lefschetz theorem in
\'{e}tale cohomology.

\begin{nt}
Tate explained the above theory in his 1958 seminar at Harvard.\footnote{See
Shatz, Math Reviews 0212073.} Douady reported on Tate's work in a Bourbaki
seminar in 1959,\footnote{Douady, Adrien, Cohomologie des groupes compacts
totalement discontinus (d'apr\`{e}s des notes de Serge Lang sur un article non
publie de Tate). S\'{e}minaire Bourbaki, Vol. 5, Exp. No. 189, 287--298, Soc.
Math. France, Paris, 1959.} and Lang included Tate's unpublished article 1958d
as Chapter VII of his 1967 book.\footnote{Lang, Serge, Rapport sur la
cohomologie des groupes. W. A. Benjamin, Inc., New York-Amsterdam 1967.} Serre
included the theory in his course at the Coll\`{e}ge de France, 1962--63; see
Serre 1964.\footnote{Serre, Jean-Pierre, Cohomologie Galoisienne. Cours au
Coll\`{e}ge de France, 1962-1963. Seconde \'{e}dition. Lecture Notes in
Mathematics 5 Springer-Verlag, Berlin-Heidelberg-New York.\label{serre1964}}
Tate himself published only the brief lectures Tate 2001.
\end{nt}

\subsection{Duality theorems\label{a5}}

In the early 1960s, Tate proved duality theorems for modules over the absolute
Galois groups of local and global fields that have become an indispensable
tool in Iwasawa theory, the theory of abelian varieties, and in other parts of
arithmetic geometry. The main global theorem was obtained independently by
Poitou, and is now referred to as the Poitou-Tate duality theorem%
\index{Poitou-Tate duality theorem}%
.

Throughout, $K$ is a field, $\bar{K}$ is a separable closure of $K$, and $G$
is the absolute Galois group $\Gal(\bar{K}/K)$. All $G$-modules are discrete
(i.e., the action is continuous for the discrete topology on the module). The
dual $M^{\prime}$ of such a module is $\Hom(M,\bar{K}^{\times})$.

\subsubsection{Local results}

Let $K$ be a nonarchimedean local field, i.e., a finite extension of
$\mathbb{Q}{}_{p}$ or $\mathbb{F}{}_{p}((t))$. Local class field theory
provides us with a canonical isomorphism $H^{2}(G,\bar{K}^{\times}%
)\simeq\mathbb{Q}{}/\mathbb{Z}{}$. Tate proved that, for every finite
$G$-module $M$ whose order $m$ is not divisible by characteristic of $K$, the
cup-product pairing%
\begin{equation}
H^{r}(G,M)\times H^{2-r}(G,M^{\prime})\rightarrow H^{2}(G,\bar{K}^{\times
})\simeq\mathbb{Q}{}/\mathbb{Z}{} \label{e29}%
\end{equation}
is a perfect duality of finite groups for all $r\in\mathbb{N}{}$ . In
particular, $H^{r}(G,M)=0$ for $r>2$. Moreover, the following holds for the
Euler-Poincar\'{e} characteristic of $M$:%
\[
\frac{\left\vert H^{0}(G,M)\right\vert \left\vert H^{2}(G,M)\right\vert
}{\left\vert H^{1}(G,M)\right\vert }=\frac{1}{(\mathcal{O}{}_{K}\colon
m\mathcal{O}{}_{K})}\text{.}%
\]
A $G$-module $M$ is said to be unramified if the inertia group $I$ in $G$ acts
trivially on $M$. When $M$ is unramified and its order is prime to the residue
characteristic, Tate proved that the submodules $H^{1}(G/I,M)$ and
$H^{1}(G/I,M^{\prime})$ of $H^{1}(G,M)$ and $H^{1}(G,M^{\prime})$ are exact
annihilators in the pairing (\ref{e29}).

Let $K=\mathbb{R}{}$. In this case, there is a canonical isomorphism
$H^{2}(G,\bar{K}^{\times})\simeq\frac{1}{2}\mathbb{Z}{}/\mathbb{Z}{}$. For any
finite $G$-module $M$, the cup-product pairing of Tate cohomology groups%
\[
\hat{H}^{r}(G,M)\times\hat{H}^{2-r}(G,M^{\prime})\rightarrow H^{2}(G,\bar
{K}^{\times})\simeq\tfrac{1}{2}\mathbb{Z}{}/\mathbb{Z}{}%
\]
is a perfect pairing for all $r\in\mathbb{Z}{}$. Moreover, $\hat{H}^{r}(G,M)$
is a finite group, killed by $2$, whose order is independent of $r$.

\subsubsection{Global results}

Let $K$ be a global field, and let $M$ be a finite $G$-module whose order is
not divisible by the characteristic of $K$. Let
\[
H^{r}(K_{v},M)=\left\{
\begin{array}
[c]{ll}%
H^{r}(G_{K_{v}},M)\quad & \text{if }v\text{ is nonarchmedean}\\
\hat{H}^{r}(G_{K_{v}},M) & \text{otherwise.}%
\end{array}
\right.
\]
The local duality results show that $\prod_{v}H^{0}(K_{v},M)$ is dual to
$\bigoplus_{v}H^{2}(K_{v},M^{\prime})$ and that $\prod\nolimits_{v}^{\prime
}H^{1}(K_{v},M)$ is dual to $\prod\nolimits_{v}^{\prime}H^{1}(K_{v},M^{\prime
})$ --- the $^{\prime}$ means that we are taking the restricted product with
respect to the subgroups $H^{1}(G_{K_{v}}/I_{v},M^{I_{v}})$.

In the table below, the homomorphisms at right are the duals of the
homomorphisms at left with $M$ replaced by $M^{\prime}$, i.e., $\beta
^{r}(M)=\alpha^{2-r}(M^{\prime})^{\ast}$ with $-^{\ast}=\Hom(-,\mathbb{Q}%
/\mathbb{Z}{})$:%
\[%
\begin{array}
[c]{lll}%
H^{0}(K,M)\overset{\alpha^{0}}{\longrightarrow}\prod_{v}H^{0}(K_{v},M) &  &
\bigoplus_{v}H^{2}(K_{v},M)\overset{\beta^{2}}{\longrightarrow}H^{0}%
(K,M^{\prime})^{\ast}\\
H^{1}(K,M)\overset{\alpha^{1}}{\longrightarrow}\prod_{v}^{\prime}H^{1}%
(K_{v},M) &  & \prod_{v}^{\prime}H^{1}(K_{v},M)\overset{\beta^{1}%
}{\longrightarrow}H^{1}(K,M^{\prime})^{\ast}\\
H^{2}(K,M)\overset{\alpha^{2}}{\longrightarrow}\bigoplus_{v}H^{2}(K_{v},M) &
& \prod\nolimits_{v}H^{0}(K_{v},M)\overset{\beta^{0}}{\longrightarrow}%
H^{2}(K,M^{\prime})^{\ast}.
\end{array}
\]

The Poitou-Tate duality theorem states that there is an exact sequence%
\begin{equation}
\minCDarrowwidth25pt\begin{CD} 0 @>>> H^{0}(K,M) @>{\alpha^{0}}>>\dstyle\prod\nolimits_{v}H^{0}(K_{v},M)@>{\beta^{0}}>> H^{2}(K,M^{\prime})^{\ast}\\ @.@.@.@VVV\\ @.H^{1}(K,M^{\prime})^{\ast}@<{\beta^{1}}<<\dstyle\prod\nolimits_{v}^{\prime}H^{1}(K_{v},M)@<{\alpha^{1}}<<H^{1}(K,M)\\ @.@VVV@.@.@.\\ @.H^{2}(K,M) @>{\alpha^{2}}>> \dstyle\bigoplus\nolimits_{v}H^{2}(K_{v},M)@>{\beta^{2}}>>H^{0}(K,M^{\prime})^{\ast}@>>>0\text{.} \end{CD} \label{e30}%
\end{equation}
(with explicit descriptions for the unnamed arrows). Moreover, for $r\geq3,$
the map%
\[
H^{r}(K,M)\rightarrow\prod\nolimits_{v\text{ real}}H^{r}(K_{v},M)
\]
is an isomorphism. In fact, the statement is more general in that one replaces
the set of all primes with a nonempty set $S$ containing the archimedean
primes in the number field case (there is then a restriction on the order of
$M$).

For the Euler-Poincar\'{e} characteristic, Tate proved that%
\begin{equation}
\frac{\left\vert H^{0}(G,M)\right\vert \left\vert H^{2}(G,M)\right\vert
}{\left\vert H^{1}(G,M)\right\vert }=\frac{1}{|M|^{r_{1}+2r_{2}}}%
\prod\nolimits_{v|\infty}\left\vert M^{G_{v}}\right\vert \label{e38}%
\end{equation}
where $r_{1}$ and $r_{2}$ are the numbers of real and complex primes.

\begin{nt}
Tate announced the above results (with brief indications of proof) in his talk
at the 1962 International Congress except for the last statement on the
Euler-Poincar\'{e} characteristic, which was announced in Tate 1966e. Tate's
proofs of the local statements were included in Serre
1964.$^{\text{\ref{serre1964}}}$ Later Tate (1966f) proved a duality theorem
for an abstract class formation, which included both the local and global
duality results, and in which the exact sequence (\ref{e30}) arises as a
sequence of $\Ext$s. This proof, as well as proofs of the formulas for the
Euler-Poincar\'{e} characteristics, are included in Milne
1986.\footnote{Milne, J. S. Arithmetic duality theorems. Perspectives in
Mathematics, 1. Academic Press, Inc., Boston, MA, 1986.\label{milne}}
\end{nt}

\subsection{Expositions\label{a6}}

The notes of the famous Artin-Tate seminar on class field theory have been a
standard reference on the topic since they first became available in 1961.
They have recently been republished in slightly revised form by the American
Mathematical Society. Tate made important contributions, both in his article
on global class field theory and in the exercises, to another classic
exposition of algebraic number theory, namely, the proceeding of the 1965
Brighton conference.\footnote{Algebraic number theory. Proceedings of an
instructional conference organized by the London Mathematical Society. Edited
by J. W. S. Cassels and A. Fr\"{o}hlich, Academic Press, London; Thompson Book
Co., Inc., Washington, D.C. 1967} His talk on Hilbert's ninth problem, which
asked for \textquotedblleft a proof of the most general reciprocity law in any
number field\textquotedblright, illuminates the problem and the work done on
it (Tate 1976a). Tate's contribution to the proceedings of the Corvallis
conference, gave a modern account of the Weil group and an explanation of the
hypothetical nonabelian reciprocity law in terms of the more general
Weil-Deligne group (Tate 1979).

\section{Abelian varieties and curves\label{b}}

In the course of proving the Riemann hypothesis for curves and abelian
varieties in the 1940s, Weil rewrote the foundations of algebraic geometry,
including the theory of abelian varieties. This made it possible to do
algebraic geometry in a rigorous fashion over arbitrary base fields. In the
late 1950s, Grothendieck rewrote the foundations again, developing the more
natural and flexible language of schemes.

\subsection{The Riemann hypothesis for curves\label{b1}}

After Hasse proved the Riemann hypothesis for elliptic curves over finite
fields in 1930, he and Deuring realized that, in order to extend the proof to
curves of higher genus, one should replace the endomorphisms of the elliptic
curve by correspondences. However, they regarded correspondences as objects in
a double field, and this approach didn't lead to a proof until Roquette
1953\footnote{Roquette, Peter, Arithmetischer Beweis der Riemannschen
Vermutung in Kongruenzfunktionenk\"{o}rpern beliebigen Geschlechts. J. Reine
Angew. Math. 191, (1953). 199--252.} (Roquette was a student of Hasse). In the
meantime Weil had realized that everything needed for the proof could be found
already in the work of the Italian geometers on correspondences, at least in
characteristic zero. In order to give a rigorous proof, he laid the
foundations for algebraic geometry over arbitrary fields,\footnote{The main
lacunae at the time were a rigorous intersection theory taking account of the
phenomenon of pure inseparability and the construction of the Jacobian variety
in nonzero characteristic.} and completed the proof of the Riemann hypothesis
for all curves over finite fields in 1945.\footnote{Weil, Andr\'{e}. Sur les
courbes alg\'{e}briques et les vari\'{e}t\'{e}s qui s'en d\'{e}duisent. Publ.
Inst. Math. Univ. Strasbourg 7 (1945).}\footnote{Much has been written on
these events. I've found the following particularly useful: Schappacher,
Norbert, The Bourbaki Congress at El Escorial and other mathematical
(non)events of 1936. The Mathematical Intelligencer, Special issue
International Congress of Mathematicians Madrid August 2006, 8-15.}

The key point of Weil's proof is that the inequality of Castelnuovo-Severi
continues to hold in characteristic $p$, i.e., for a divisor $D$ on the
product of two complete nonsingular curves $C$ and $C^{\prime}$ over an
algebraically closed field,
\begin{equation}
\lbrack D\cdot D]\leq2dd^{\prime} \label{e8}%
\end{equation}
where $d=[D\cdot(P\times C^{\prime})]$ and $d^{\prime}=[D\cdot(C\times
P^{\prime})]$ are the degrees of $D$ over $C$ and $C^{\prime}$ respectively.
Mattuck and Tate (1958a) showed that it is possible to derive (\ref{e8})
directly and easily from the Riemann-Roch theorem for surfaces, for which they
were able to appeal to Zariski 1952\footnote{Zariski, Oscar, Complete linear
systems on normal varieties and a generalization of a lemma of
Enriques-Severi. Ann. of Math. (2) 55, (1952). 552--592.} or to a
sheaf-theoretic proof of Serre which is sketched in Zariski
1956.\footnote{Zariski, Oscar, Scientific report on the Second Summer
Institute, III Algebraic Sheaf Theory, Bull. Amer. Math. Soc. 62 (1956), 117--141.}

The Mattuck-Tate proof%
\index{Mattuck-Tate proof}
is the most attractive geometric proof of Weil's theorem.
Grothendieck\footnote{Grothendieck, A. Sur une note de Mattuck-Tate. J. Reine
Angew. Math. 200 1958 208--215.} simplified it further by showing that the
Castelnuovo-Severi inequality can most naturally be derived from the Hodge
index theorem for surfaces, which itself can be derived directly from the
Riemann-Roch theorem.

Hodge proved his index theorem for smooth projective varieties over
$\mathbb{C}{}$. That it should hold for such varieties in nonzero
characteristic is known as Grothendieck's \textquotedblleft Hodge standard
conjecture\textquotedblright, whose proof Grothendieck calls one of the
\textquotedblleft most urgent tasks in algebraic geometry\textquotedblright%
.\footnote{Grothendieck, A. Standard conjectures on algebraic cycles. 1969
Algebraic Geometry (Internat. Colloq., Tata Inst. Fund. Res., Bombay, 1968)
pp. 193--199 Oxford Univ. Press, London.} In the more than forty years since
Grothendieck formulated the conjecture, almost no progress has been made
towards its proof --- even in characteristic zero, there exists no algebraic
proof in dimensions greater than $2$.

\subsubsection{The Tate module of an abelian variety}

Let $A$ be an abelian variety over a field $k$. For a prime $l$, let
$A(l)=\bigcup A(k^{\mathrm{sep}})_{l^{n}}$ where $A(k^{\mathrm{sep}})_{l^{n}%
}=\Ker(A(k^{\mathrm{sep}})\overset{l^{n}}{\longrightarrow}A(k^{\mathrm{sep}%
}))$. Then $A\rightsquigarrow A(l)$ is a functor from abelian varieties over
$k$ to $l$-divisible groups equipped with an action of $\Gal(k^{\mathrm{sep}%
}/k)$. When $l\neq\mathrm{char}(k)$, $A(l)\simeq(\mathbb{Q}{}_{\ell
}/\mathbb{Z}{}_{\ell})^{2\dim A}$, and Weil used $A(l)$ to study the
endomorphisms of $A$. Tate observed that it is more convenient to work with%
\[
T_{l}A=\varprojlim A(k^{\mathrm{sep}})_{l^{n}},
\]
which is a free $\mathbb{Z}{}_{l}$-module of rank $2\dim A$ when
$l\neq\mathrm{char}(k)$ --- this is now called the Tate module%
\index{Tate module}
of $A$.\label{Tatemodule}

\subsection{Heights on abelian varieties\label{b2}}

\subsubsection{The N\'{e}ron-Tate (canonical) height}

Let $K$ be a number field, and normalize the absolute values $|\cdot|_{v}$ of
$K$ so that the product formula holds:%
\[
\prod\nolimits_{v}|a|_{v}=1\text{ for all }v\in K^{\times}\text{.}%
\]
The logarithmic height of a point $P=(a_{0}\colon\ldots\colon a_{n})$ of
$\mathbb{P}{}^{n}(K)$ is defined to be%
\[
h(P)=\log\left(  \prod\nolimits_{v}\max\{|a_{0}|_{v},\ldots,|a_{n}%
|_{v}\}\right)  .
\]
The product formula shows that this is independent of the representation of
$P$.

Let $X$ be a projective variety. A morphism $f\colon X\rightarrow\mathbb{P}%
{}^{n}$ from $X$ into projective space defines a height function
$h_{f}(P)=h(f(P))$ on $X$. In a Short Communication at the 1958 International
Congress, N\'{e}ron conjectured that, in certain cases, the height is given by
a quadratic form.\footnote{Only the title, \textit{Valeur asymptotique du
nombre des points rationnels de hauteur bornd\'{e}e sur une courbe
elliptique}, of N\'{e}ron's communication is included in the Proceedings. The
sentence paraphrases one from: Lang, Serge. Les formes bilin\'{e}aires de
N\'{e}ron et Tate. S\'{e}minaire Bourbaki, 1963/64, Fasc. 3, Expos\'{e}
274.\label{lang}} Tate proved this for abelian varieties by a simple direct argument.

Let $A$ be an abelian variety over a number field $K$. A nonconstant map
$f\colon A\rightarrow\mathbb{P}{}^{n}$ of $A$ into projective space is said to
be symmetric if the inverse image $D$ of a hyperplane is linearly equivalent
to $(-1)^{\ast}D$. For a symmetric embedding $f$, Tate proved that there
exists a unique \textit{quadratic} map $\hat{h}\colon A(K)\rightarrow
\mathbb{R}{}$ such that $\hat{h}(P)-h_{f}(P)$ is bounded on $A(K)$. To say
that $\hat{h}$ is quadratic means that $\hat{h}(2P)=4\hat{h}(P)$ and that the
function
\begin{equation}
P,Q\mapsto\frac{1}{2}\left(  \hat{h}(P+Q)-\hat{h}(P)-\hat{h}(Q)\right)
\label{e2}%
\end{equation}
is bi-additive on $A(K)\times A(K)$.

Note first that there exists at most one function $\hat{h}\colon
A(K)\rightarrow\mathbb{R}{}$ such that (a) $\hat{h}(P)-h_{f}(P)$ is bounded on
$A(K)$, and (b) $\hat{h}(2P)=4\hat{h}(P)$ for all $P\in A(K)$. \noindent
Indeed, if $\hat{h}$ satisfies (a) with bound $B$, then%
\[
\left\vert \hat{h}(2^{n}P)-h_{f}(2^{n}P)\right\vert \leq B
\]
for all $P\in A(K)$ and all $n\geq0$. If in addition it satisfies (b), then%
\[
\left\vert \hat{h}(P)-\frac{h_{f}(2^{n}P)}{4^{n}}\right\vert \leq\frac
{B}{4^{n}}%
\]
for all $n$, and so%
\begin{equation}
\hat{h}(P)=\lim_{n\rightarrow\infty}\frac{h_{f}(2^{n}P)}{4^{n}}. \label{e1}%
\end{equation}
Tate used the equation (\ref{e1}) to define $\hat{h}$, and applied results of
Weil on abelian varieties to verify that it is quadratic.

Let $A^{\prime}$ be the dual abelian variety to $A$. For a map $f\colon
A\rightarrow\mathbb{P}^{n}$ corresponding to a divisor $D$, let $\varphi
_{f}\colon A(K)\rightarrow A^{\prime}(K)$ be the map sending $P$ to the point
on $A$ represented by the divisor $(D+P)-D$. Tate showed that there is a
unique bi-additive pairing
\begin{equation}
\langle\,\,,\,\,\rangle\colon A^{\prime}(K)\times A(K)\rightarrow\mathbb{R}{}
\label{e4}%
\end{equation}
such that, for every symmetric $f$, the function $\langle\varphi
_{f}(P),P\rangle+2h_{f}(P)$ is bounded on $A(K)$.

N\'{e}ron (1965)\footnote{N\'{e}ron, A., Quasi-fonctions et hauteurs sur les
vari\'{e}t\'{e}s ab\'{e}liennes. Ann. of Math. (2) 82 1965 249--331.} found
his own construction of $\hat{h}$, which is much longer than Tate's, but which
has the advantage of expressing $\hat{h}$ as a sum of local heights. The
height function $\hat{h}$, is now called the N\'{e}ron-Tate, or canonical,
height%
\index{N\'{e}ron-Tate height}%
. It plays a fundamental role in arithmetic geometry.

\begin{nt}
Tate explained his construction in his course on abelian varieties at Harvard
in the fall of 1962, but did not publish it. However, it was soon published by
others.\footnote{Lang, S., see footnote \ref{lang}, and Diophantine
approximations on toruses. Amer. J. Math. 86 1964 521--533; Manin, Ju. I., The
Tate height of points on an abelian variety, its variants and applications.
(Russian) Izv. Akad. Nauk SSSR Ser. Mat. 28 1964 1363--1390.}
\end{nt}

\subsubsection{Variation of the canonical height of a point depending on a
parameter}

Let $T$ be an algebraic curve over $\mathbb{Q}{}^{\mathrm{al}}$, and let
$E\rightarrow T$ be an algebraic family of elliptic curves parametrised by
$T$. Let $P\colon T\rightarrow E$ be a section of $E/T$, and let $\hat{h}_{t}$
be the N\'{e}ron-Tate height on the fibre $E_{t}$ of $E/T$ over a closed point
$t\ $of $T$. Tate (1983a) proves that the map $t\mapsto\hat{h}_{t}(P_{t})$ is
a height function on the curve $T$ for a certain divisor class $q(P)$ on $T$;
moreover, the degree of $q(P)$ is the N\'{e}ron-Tate height of $P$ regarded as
a point on the generic fibre of $E/T$.

As Tate noted \textquotedblleft The main obstacle to extending the theorem in
this paper to abelian varieties seems to be the lack of a canonical
compactification of the N\'{e}ron model in higher
dimensions.\textquotedblright\ After Faltings compactified the moduli stack of
abelian varieties,\footnote{Faltings, G. Arithmetische Kompaktifizierung des
Modulraums der abelschen Variet\"{a}ten. Workshop Bonn 1984 (Bonn, 1984),
321--383, Lecture Notes in Math., 1111, Springer, Berlin, 1985.} one of his
students, William Green, extended Tate's theorem to abelian
varieties.\footnote{Green, William, Heights in families of abelian varieties.
Duke Math. J. 58 (1989), no. 3, 617--632.}

\subsubsection{Height pairings via biextensions}

Let $A$ be an abelian variety over a number field $K$, and let $A^{\prime}$ be
its dual. The classical N\'{e}ron-Tate height pairing is a pairing%
\[
A(K)\times A^{\prime}(K)\rightarrow\mathbb{R}%
\]
whose kernels are precisely the torsion subgroups of $A(K)$ and $A^{\prime
}(K)$. In order, for example, to state a $p$-adic version of the conjecture of
Birch and Swinnerton-Dyer, it is necessary to define a $\mathbb{Q}{}_{p}%
$-valued height pairing,%
\[
A(K)\times A^{\prime}(K)\rightarrow\mathbb{Q}{}_{p}\text{.}%
\]
When $A$ has good ordinary or multiplicative reduction at the $p$-adic primes,
Mazur and Tate (1983b) use the expression of the duality between $A$ and
$A^{\prime}$ in terms of biextensions, and exploit the local splittings of
these biextensions, to define such pairings. They compare their definition
with other suggested definitions. It is not known whether the pairings are
nondegenerate modulo torsion.

\subsection{The cohomology of abelian varieties\label{b3}}

\subsubsection{The local duality for abelian varieties}

Let $A$ be an abelian variety over a field $k$. A principal homogeneous space
over $A$ is a variety $V$ over $k$ together with a regular map $A\times
V\rightarrow V$ such that, for every field $K$ containing $k$ for which $V(K)$
is nonempty, the pairing $A(K)\times V(K)\rightarrow V(K)$ makes $V(K)$ into a
principal homogeneous space for $A(K)$ in the usual sense. The isomorphism
classes of principal homogeneous spaces form a group, which Tate (1958b) named
the Weil-Ch\^{a}telet group, and denoted $\mathrm{WC}(A/k)$.

For a finite extension $k$ of $\mathbb{Q}{}_{p}$, local class field theory
provides a canonical isomorphism $H^{2}(k,\mathbb{G}_{m})\simeq\mathbb{Q}%
/\mathbb{Z}{}$. Tate (ibid.) defines an \textquotedblleft
augmented\textquotedblright\ cup-product pairing\label{ldt}%
\begin{equation}
H^{r}(k,A)\times H^{1-r}(k,A^{\prime})\rightarrow H^{2}(k,\mathbb{G}%
_{m})\simeq\mathbb{Q}{}/\mathbb{Z}{}, \label{e39}%
\end{equation}
and proves that it is a perfect duality for $r=1$. In other words, the
discrete group $\mathrm{WC}(A/k)$ is canonically dual to the compact group
$A^{\prime}(k)$. Later, he showed that (\ref{e39}) is a perfect duality for
all $r$. In the case $k=\mathbb{R}{}$, he proved that $H^{1}(\mathbb{R}{},A)$
is canonically dual to $A^{\prime}(\mathbb{R}{})/A^{\prime}(\mathbb{R}%
{})^{\circ}=\pi_{0}(A^{\prime}(\mathbb{R}{}))$.

\begin{nt}
The above results are proved in Tate 1958b, 1959b, or 1962d. The analogous
statements for local fields of characteristic $p$ are proved in Milne
1970.\footnote{Milne, J. S. Weil-Ch\^{a}telet groups over local fields. Ann.
Sci. \'{E}cole Norm. Sup. (4) 3 1970 273--284; ibid. 5 (1972), 261-264.}
\end{nt}

\subsubsection{Principal homogeneous spaces over abelian varieties}

Lang and Tate (1958c) explain the relation between the set $\mathrm{WC}(A/k)$
of isomorphism classes of principal homogeneous spaces over a group variety
$A$ and the Galois cohomology group $H^{1}(k,A)$. Briefly, there is a
canonical injective map $\mathrm{WC}(A/k)\rightarrow H^{1}(k,A)$ which Weil's
descent theorems show to be surjective. This generalizes results of Ch\^{a}telet.

Let $K$ be a field complete with respect to a discrete valuation with residue
field $k$, and let $A$ be an abelian variety over $K$ with good reduction to
an abelian variety $\bar{A}$ over $k$. Then, for any integer $m$ prime to the
characteristic of $k$, Lang and Tate (ibid.) prove that there is a canonical
exact sequence%
\[
0\rightarrow H^{1}(k,\bar{A})_{m}\rightarrow H^{1}(k,A)_{m}\rightarrow
\Hom(\mu_{m}(k),\bar{A}(k^{\mathrm{sep}})_{m})\rightarrow0.
\]
In the final section of the article, they study abelian varieties over global
fields. In particular, they prove the weak Mordell-Weil theorem.

As Cassels wrote,\footnote{Math. Reviews 0138625.} the article Tate 1962b
provides \textquotedblleft A laconic but useful review of the existing state
of knowledge [on principal homogeneous spaces for abelian varieties] for
different types of groundfield.\textquotedblright\ 

\subsubsection{The conjecture of Birch and Swinnerton-Dyer}

For an elliptic curve $A$ over $\mathbb{Q}{}$, Mordell showed that the group
$A(\mathbb{Q}{})$ is finitely generated. It is easy to compute the torsion
subgroup of $A(\mathbb{Q}{})$, but there is at present no proven algorithm for
computing its rank $r(A)$. Computations led Birch and Swinnerton-Dyer to
conjecture that $r(A)$ is equal to the order of the zero at $1$ of the
$L$-series of $A$, and further work led to a more precise conjecture. Tate
(1966e) formulated the analogues of their conjectures for an abelian variety
$A$ over a global field $K$.

Let $v$ be a nonarchimedean prime of $K$, and let $\kappa(v)$ be the
corresponding residue field. If $A$ has good reduction at $v$, then it gives
rise to an abelian variety $A(v)$ over $\kappa(v)$. The characteristic
polynomial of the Frobenius endomorphism of $A(v)$ is a polynomial $P_{v}(T)$
of degree $2d$ with coefficients in $\mathbb{Z}$ such that, when we factor it
as $P_{v}(T)=\prod_{i}(1-a_{i}T)$, then $\prod_{i}(1-a_{i}^{m})$ is the number
of points on $A(v)$ with coordinates in the finite field of degree $m$ over
$\kappa(v)$. For any finite set $S$ of primes of $K$ including the archimedean
primes and those where $A$ has bad reduction, we define the $L$-series
$L_{S}(s,A)$ by the formula
\[
L_{S}(A,s)={\textstyle\prod_{v\notin S}}P_{v}(A,Nv^{-s})^{-1}%
\]
where $Nv=[\kappa(v)]$. The product converges for $\Re(s)>3/2$, and it is
conjectured that $L_{S}(A,s)$ can be analytically continued to a meromorphic
function on the whole complex plane. This is known in the function field case,
and over $\mathbb{Q}{}$ for elliptic curves. The analogue of the first
conjecture of Birch and Swinnerton-Dyer for $A$ is that%
\begin{equation}
L_{S}(A,s)\text{ has a zero of order }r(A)\text{ at }s=1\text{.} \label{e3}%
\end{equation}

Let $\omega$ be a nonzero global differential $d$-form on $A$. As
$\Gamma(A,\Omega_{A}^{d})$ has dimension $1$, $\omega$ is uniquely determined
up to multiplication by an element of $K^{\times}$. For each nonarchimedean
prime $v$ of $K$, let $\mu_{v}$ be the Haar measure on $K_{v}$ for which
$\mathcal{O}{}_{v}$ has measure $1$, and for each archimedean prime, take
$\mu_{v}$ to be the usual Lebesgue measure on $K_{v}$. Define%
\[
\mu_{v}(A,\omega)=\int_{A(K_{v})}|\omega|_{v}\mu_{v}^{d}%
\]
Let $\mu$ be the measure ${\textstyle\prod}\mu_{v}$ on the ad\`{e}le ring
$\mathbb{A}{}_{K}$ of $K$, and set $|\mu|=\int_{\mathbb{A}{}_{K}/K}\mu$. For
any finite set $S$ of primes of $K$ including all archimedean primes and those
nonarchimedean primes for which $A$ has bad reduction or such that $\omega$
does not reduce to a nonzero differential $d$-form on $A(v)$, we define%
\[
L_{S}^{\ast}(s,A)=L_{S}(s,A)\frac{|\mu|^{d}}{%
{\textstyle\prod_{v\in S}}
\mu_{v}(A,\omega)}\text{.}%
\]
The product formula shows that this is independent of the choice of $\omega$.
The asymptotic behaviour of $L_{S}^{\ast}(s,A)$ as $s\rightarrow1$, which is
all we are interested in, doesn't depend on $S$. The analogue of the second
conjecture of Birch and Swinnerton-Dyer is that%
\begin{equation}
\lim_{s\rightarrow1}\frac{L_{S}^{\ast}(s,A)}{(s-1)^{r(A)}}=\frac
{[\TS(A)]\cdot\left\vert D\right\vert }{[A^{\prime}(K)_{\mathrm{tors}}%
]\cdot\lbrack A(K)_{\mathrm{tors}}]} \label{e5}%
\end{equation}
where $\TS(A)$ is the Tate-Shafarevich group%
\index{Tate-Shafarevich group}
of $A$,%
\[
\TS(A)\overset{\textup{{\tiny def}}}{=}\Ker\left(  H^{1}(K,A)\rightarrow
\prod\nolimits_{v}H^{1}(K_{v},A)\right)  ,
\]
which is conjectured to be finite, and $D$ is the discriminant of the height
pairing (\ref{e4}), which is known to be nonzero.

\subsubsection{Global duality}

In his talk at the 1962 International Congress, Tate stated the local duality
theorems reviewed above (p.\pageref{ldt}), and he announced some global
theorems which we now discuss.

In their computations, Birch and Swinnerton-Dyer found that the order of the
Tate-Shafarevich group predicted by (\ref{e5}) is always a square. Cassels and
Tate conjectured independently that the explanation for this is that there
exists an alternating pairing%
\begin{equation}
\TS(A)\times\TS(A)\rightarrow\mathbb{Q}{}/\mathbb{Z}{} \label{e35}%
\end{equation}
that annihilates only the divisible subgroup of $\TS(A)$. Cassels proved this
for an elliptic curve over a number field.\footnote{Cassels, J. W. S.
Arithmetic on curves of genus $1$. IV. Proof of the Hauptvermutung. J. Reine
Angew. Math. 211 1962 95--112.} For an abelian variety $A$ and its dual
abelian variety $A^{\prime}$, Tate proved that there exists a canonical
pairing%
\begin{equation}
\TS(A)\times\TS(A^{\prime})\rightarrow\mathbb{Q}{}/\mathbb{Z}{} \label{e36}%
\end{equation}
that annihilates only the divisible subgroups; moreover, for a divisor $D$ on
$A$ and the homomorphism $\varphi_{D}\colon A\rightarrow A^{\prime}$,
$a\mapsto\lbrack D_{a}-D]$, it defines, the pair $(\alpha,\varphi_{D}%
(\alpha))$ maps to zero under (\ref{e36}) for all $\alpha\in\TS(A)$. The
pairing (\ref{e36}), or one of its several variants, is now called the
Cassels-Tate pairing.%
\index{Cassels-Tate pairing}%

For an elliptic curve $A$ over a number field $k$ such that $\TS(A)$ is
finite, Cassels determined the Pontryagin dual of the exact sequence%
\begin{equation}
0\rightarrow\TS(A)\rightarrow H^{1}(k,A)\rightarrow\bigoplus\nolimits_{v}%
H^{1}(k_{v},A)\rightarrow\TB(A)\rightarrow0 \label{e6}%
\end{equation}
(regarded as a sequence of discrete groups). Assume that $\TS(A)$ is finite.
Using Tate's local duality theorem (see p.\pageref{ldt}) for an elliptic
curve, Cassels (1964)\footnote{Cassels, J. W. S., Arithmetic on curves of
genus 1. VII. The dual exact sequence. J. Reine Angew. Math. 216 1964
150--158.} showed that the dual of (\ref{e6}) takes the form%
\begin{equation}
0\leftarrow\text{ }\TS\leftarrow\Theta\leftarrow\prod\nolimits_{v}%
A(k_{v})^{\prime}\leftarrow\widetilde{A(k)}\leftarrow0 \label{e37}%
\end{equation}
for a certain explicit $\Theta$ and with $\widetilde{A(k)}$ equal to the
closure of $A(k)$ in $\prod\nolimits_{v}A(k_{v})^{\prime}$. Tate proved the
same statement for abelian varieties over number fields, except that, in
(\ref{e37}), it is necessary to replace $A$ with its dual $A^{\prime}$. So
modified, the sequence (\ref{e37}) is now called the Cassels-Tate dual exact
sequence.%
\index{Cassels-Tate sequence}%

Let $A$ and $B$ be isogenous elliptic curves over a number field. Then
$L_{S}(s,A)=L_{S}(s,B)$ and $r(A)=r(B)$, and so the first conjecture of Birch
and Swinnerton-Dyer is true for $A$ if and only if it is true for $B$. Cassels
proved the same statement for the second conjecture.\footnote{Cassels, J. W.
S. Arithmetic on curves of genus 1. VIII. On conjectures of Birch and
Swinnerton-Dyer. J. Reine Angew. Math. 217 1965 180--199.} This amounts to
showing that a certain product of terms doesn't change in passing from $A$ to
$B$ (even though the individual terms may change). Using his duality theorems
and the formula (\ref{e38}), p.\pageref{e38}, for the Euler-Poincar\'{e}
characteristic, Tate (1966e, 2.1) proved the same result for abelian varieties
over number fields.

Tate's global duality theorems were widely used, even before there were
published proofs. Since 1994, the duality theorems have been used in cryptography.

\begin{nt}
Tate's results are more general and complete than stated above; in particular,
he works with a nonempty set $S$ of primes of $k$ (not necessarily the
complete set). Proofs of the theorems of Tate in this subsection can be found
in Milne 1986.$^{\text{\ref{milne}}}$
\end{nt}

\subsection{Serre-Tate liftings of abelian varieties\label{b4}}

In a talk at the 1964 Woods Hole conference, Tate discussed some results of
his and Serre on the lifting of abelian varieties from characteristic $p$.

For an abelian scheme $A$ over a ring $R$, let $A_{n}$ denote the kernel of
$A\overset{p^{n}}{\longrightarrow}A$ regarded as a finite group scheme over
$R$, and let $A(p)$ denote the direct system%
\[
A_{1}\hookrightarrow A_{2}\hookrightarrow\cdots\hookrightarrow A_{n}%
\hookrightarrow\cdots
\]
of finite group schemes. Let $R$ be an artinian local ring with residue field
$k$ of characteristic $p\neq0$. An abelian scheme $A$ over $R$ defines an
abelian variety $\bar{A}$ over $k$ and a system of finite group schemes $A(p)$
over $R$. Serre and Tate prove that the functor%
\[
A\rightsquigarrow(\bar{A},A(p))
\]
is an equivalence of categories (Serre-Tate theorem%
\index{Serre-Tate theorem}%
). In particular, to lift an abelian variety $A$ from $k$ to $R$ amounts to
lifting the system of finite group schemes $A(p)$.

This has many important consequences.

\begin{itemize}
\item Let $A$ and $B$ be abelian schemes over a complete local noetherian ring
$R$ with residue field a field $k$ of characteristic $p\neq0$. A homomorphism
$f\colon\bar{A}\rightarrow\bar{B}$ of abelian varieties over $k$ lifts to a
homomorphism $A\rightarrow B$ of abelian schemes over $R$ if and only if
$f(p)\colon\bar{A}(p)\rightarrow\bar{B}(p)$ lifts to $R$. For the artinian
quotients of $R$, this is part of the above statement, and the statement for
$R$ follows by passing to the limit over the artinian quotients of $R$ and
applying a theorem of Grothendieck.

\item Let $A$ be an abelian variety over a perfect field $k$ of characteristic
$p\neq0$. If $A$ is ordinary, then%
\[
A_{n}\approx(\mathbb{Z}{}/p^{n}\mathbb{Z}{})^{\dim A}\times(\mu_{p^{n}})^{\dim
A},
\]
and so each $A_{n}$ has a canonical lifting to a finite group scheme over the
ring $W(k)$ of Witt vectors of $k$. Thus $A$ has a canonical lifting to an
abelian scheme over $W(k)$ (at least formally, but the existence of
polarizations implies that the formal abelian scheme is an abelian scheme).
Deligne has used this to give a \textquotedblleft linear
algebra\textquotedblright\ description of the category of ordinary abelian
varieties over a finite field similar to the classical description of abelian
varieties over $\mathbb{C}{}$.\footnote{Deligne, Pierre. Vari\'{e}t\'{e}s
ab\'{e}liennes ordinaires sur un corps fini. Invent. Math. 8 1969 238--243.}

\item Over a ring $R$ in which $p$ is nilpotent, the infinitesimal deformation
theory of $A$ is equivalent to the infinitesimal deformation theory of $A(p)$.
For example, when $A$ is ordinary, this implies that the local deformation
space of an ordinary abelian variety $A$ over $k$ has a natural structure of a
formal torus over $W(k)$ of relative dimension $\dim(A)^{2}$.
\end{itemize}

\begin{nt}
Lifting results were known to Hasse and Deuring for elliptic curves. The
canonical lifting of an ordinary abelian variety was found by Serre, prompting
Tate to prove the general result. Lubin, Serre, and Tate 1964b contains a
sketch of the proofs. The liftings of $A$ obtained from liftings of the system
$A(p)$ are sometimes called Serre-Tate liftings%
\index{Serre-Tate liftings}%
, especially in the ordinary case. Messing 1972\footnote{Messing, William. The
crystals associated to Barsotti-Tate groups: with applications to abelian
schemes. Lecture Notes in Mathematics, Vol. 264. Springer-Verlag, Berlin-New
York, 1972.} includes a proof of the Serre-Tate theorem.
\end{nt}

\subsection{Mumford-Tate groups and the Mumford-Tate conjecture \label{b5}}

In 1965, Mumford gave a talk at the AMS Summer Institute\footnote{Mumford,
David. Families of abelian varieties. 1966 Algebraic Groups and Discontinuous
Subgroups (Proc. Sympos. Pure Math., Boulder, Colo., 1965) pp. 347--351 Amer.
Math. Soc., Providence, R.I.} whose results he described as being
\textquotedblleft partly joint work with J. Tate\textquotedblright. In it, he
attached a reductive group to an abelian variety, and stated a conjecture. The
first is now called the Mumford-Tate group%
\index{Mumford-Tate group}%
, and the second is the Mumford-Tate conjecture.%
\index{Mumford-Tate conjecture}%

Let $A$ be a complex abelian variety of dimension $g$. Then $V\overset
{\textup{{\tiny def}}}{=}H_{1}(A,\mathbb{Q}{})$ is a $\mathbb{Q}{}$-vector
space of dimension $2g$ whose tensor product with $\mathbb{R}{}$ acquires a
complex structure through the canonical isomorphism
\[
H_{1}(A,\mathbb{Q}{})_{\mathbb{R}{}}{}\simeq\Tgt_{0}(A).
\]
Let $u\colon U^{1}\rightarrow\GL(V_{\mathbb{R}{}})$ be the homomorphism
describing this complex structure, where $U^{1}=\{z\in\mathbb{C}{}%
\mid\left\vert z\right\vert =1\}$. The Mumford-Tate group of $A$ is defined to
be the smallest algebraic subgroup $H$ of $\GL_{V}$ such that $H(\mathbb{R}%
{})$ contains $u(U^{1})$.\footnote{Better, it should be thought of as the pair
$(H,u)$.} Then $H$ is a reductive algebraic group over $\mathbb{Q}{}$, which
acts on $H^{\ast}(A^{r},\mathbb{Q}{})$, $r\in\mathbb{N}{}$, through the
isomorphisms%
\begin{align*}
H^{\ast}(A^{r},\mathbb{Q}{})  &  \simeq\bigwedge\nolimits^{\ast}H^{1}%
(A^{r},\mathbb{Q}{}),\quad\\
H^{1}(A^{r},\mathbb{Q}{})  &  \simeq rH^{1}(A,\mathbb{Q}{}),\quad\\
H^{1}(A,\mathbb{Q}{})  &  \simeq\Hom(V,\mathbb{Q}{}).
\end{align*}
It can be characterized as the algebraic subgroup of $\GL_{V}$ that fixes
exactly the Hodge tensors in the spaces $H^{\ast}(A^{r},\mathbb{Q}{})$, i.e.,
the elements of the $\mathbb{Q}{}$-spaces%
\[
H^{2p}(A^{r},\mathbb{Q}{})\cap\bigoplus H^{p,p}(A^{r})\text{.}%
\]

Let $A$ be an abelian variety over a number field $k$, and, for a prime number
$l$, let%
\[
V_{l}A=\mathbb{Q}{}_{l}\otimes_{\mathbb{Z}{}_{l}}T_{l}A
\]
where $T_{l}A$ is the Tate module of $A$ (p.\pageref{Tatemodule}). Then%
\[
V_{l}A\simeq\mathbb{Q}_{l}\otimes_{\mathbb{Q}{}}{}H_{1}(A_{\mathbb{C}{}%
},\mathbb{Q}{})\text{.}%
\]
The Galois group $G(k^{\mathrm{al}}/k)$ acts on $A(k^{\mathrm{al}})$, and
hence there is a representation%
\[
\rho_{l}\colon G(k^{\mathrm{al}}/k)\rightarrow\GL(V_{l}A)
\]
The Zariski closure $H_{l}(A)$ of $\rho_{l}(G(k^{\mathrm{al}}/k))$ is an
algebraic group in $\GL_{V_{l}A}$. Although $H_{l}(A)$ may change when $k$ is
replaced by a finite extension, its identity component $H_{l}(A)^{\circ}$ does
not and can be thought of as determining the image of $\rho_{l}$ up to finite
groups. The Mumford-Tate conjecture states that
\[
H_{l}(A)^{\circ}=(\text{Mumford-Tate group of }A_{\mathbb{C}{}})_{\mathbb{Q}%
{}_{l}}\text{ inside }\GL_{V_{l}A}\simeq(\GL_{H_{1}(A_{\mathbb{C}{}%
},\mathbb{Q}{})})_{\mathbb{Q}{}_{l}}\text{.}%
\]
In particular, it posits that the $\mathbb{Q}{}_{l}$-algebraic groups
$H_{l}(A)^{\circ}$ are independent of $l$ in the sense that they all arise by
base change from a single algebraic group over $\mathbb{Q}{}$. In the presence
of the Mumford-Tate conjecture, the Hodge and Tate conjectures for $A$ are
equivalent. Much is known about the Mumford-Tate conjecture.

Let $H$ be the Mumford-Tate group of an abelian variety $A$, and let $u\colon
U^{1}\rightarrow H(\mathbb{R}{})$ be the above homomorphism. The centralizer
$K$ of $u$ in $H(\mathbb{R}{})$ is a maximal compact subgroup of
$H(\mathbb{R}{})$, and the quotient manifold $X=H(\mathbb{R}{})/K$ has a
unique complex structure for which $u(z)$ acts on the tangent space at the
origin as multiplication by $z$. With this structure $X$ is isomorphic to a
bounded symmetric domain, and it supports a family of abelian varieties whose
Mumford-Tate groups \textquotedblleft refine\textquotedblright\ that of $A$.
The quotients of $X$ by congruence subgroups of $H(\mathbb{Q)}$ are connected
Shimura varieties.

The notion of a Mumford-Tate group has a natural generalization to an
arbitrary polarizable rational Hodge structure. In this case the quotient
space $X$ is a homogeneous complex manifold, but it is not necessarily a
bounded symmetric domain. The complex manifolds arising in this way were
called Mumford-Tate domains%
\index{Mumford-Tate domains}
by Green, Griffiths, and Kerr.\footnote{Green, Mark; Griffiths, Phillip; Kerr,
Matt, Mumford-Tate domains. Boll. Unione Mat. Ital. (9) 3 (2010), no. 2,
281--307.} As these authors say: \textquotedblleft Mumford-Tate groups have
emerged as the principal symmetry groups in Hodge theory.\textquotedblright

\subsection{Abelian varieties over finite fields (Weil, Tate, Honda
theory)\label{b6}}

Consider the category whose objects are the abelian varieties over a field $k$
and whose morphisms are given by%
\[
\Hom^{0}(A,B)\overset{\textup{{\tiny def}}}{=}\Hom(A,B)\otimes\mathbb{Q}{}.
\]
Weil's results\footnote{Weil, Andr\'{e}. Vari\'{e}t\'{e}s ab\'{e}liennes et
courbes alg\'{e}briques. Actualit\'{e}s Sci. Ind., no. 1064 = Publ. Inst.
Math. Univ. Strasbourg 8 (1946). Hermann \& Cie., Paris, 1948.} imply that
this is a semisimple abelian category whose endomorphism algebras are finite
dimensional $\mathbb{Q}{}$-algebras. Thus, to describe the category up to
equivalence, it suffices to list the isomorphism classes of simple objects
and, for each class, describe the endomorphism algebra of an object in the
class. This the theory of Weil, Tate, and Honda does when $k$ is finite.
Briefly: Weil showed that there is a well-defined map from isogeny classes of
simple abelian varieties to conjugacy classes of Weil numbers, Tate proved
that the map is injective and determined the endomorphism algebra of each
simple class, and Honda used the theory of Shimura and Taniyama to prove that
the map is surjective.

In more detail, let $k$ be a field with $q=p^{a}$ elements. Each abelian
variety $A$ over $k$ admits a Frobenius endomorphism $\pi_{A}$, which acts on
the $k^{\mathrm{al}}$-points of $A$ as $(a_{0}\colon a_{1}\colon\ldots
)\mapsto(a_{0}^{q}\colon a_{1}^{q}\colon\ldots)$. Weil proved that the image
of $\pi_{A}$ in $\mathbb{C}{}$ under any homomorphism $\mathbb{Q}{}[\pi
_{A}]\rightarrow\mathbb{C}{}$ is a Weil $q$-integer, i.e., it is an algebraic
integer with absolute value $q^{\frac{1}{2}}$ (this is the Riemann
hypothesis). Thus, attached to every simple abelian variety $A$ over $k$,
there is a conjugacy class of Weil $q$-integers. Isogenous simple abelian
varieties give the same conjugacy class.

Tate (1966b) proved that a simple abelian $A$ is determined up to isogeny by
the conjugacy class of $\pi_{A}$, and moreover, that $\mathbb{Q}{}[\pi_{A}]$
is the centre of $\End^{0}(A)$. Since $\End^{0}(A)$ is a division algebra with
centre the field $\mathbb{Q}{}[\pi_{A}]$, class field theory shows that its
isomorphism class is determined by its invariants at the primes $v$ of
$\mathbb{Q}{}[\pi_{A}]$. These Tate determined as follows:%
\[
\inv_{v}(\End^{0}(A))=\left\{
\begin{array}
[c]{ll}%
\frac{1}{2} & \text{if }v\text{ is real,}\\
\dstyle\frac{\ord_{v}(\pi_{A})}{\ord_{v}(q)}[\mathbb{Q}{}[\pi_{A}]_{v}%
\colon\mathbb{Q}{}_{p}] & \text{if }v|p\text{,}\\
0 & \text{otherwise.}%
\end{array}
\right.
\]
Moreover,%
\[
2\dim A=[\End^{0}(A)\colon\mathbb{Q}{}[\pi_{A}]]^{\frac{1}{2}}\cdot
\lbrack\mathbb{Q}{}[\pi_{A}]\colon\mathbb{Q}{}]\text{.}%
\]

The abelian varieties of CM-type over $\mathbb{C}{}$ are classified up to
isogeny by their CM-types, and every such abelian variety has a model over
$\mathbb{Q}{}^{\mathrm{al}}$. When we choose a $p$-adic prime of $\mathbb{Q}%
{}^{\mathrm{al}}$, an abelian variety $A$ of CM-type over $\mathbb{Q}%
{}^{\mathrm{al}}$ specializes to an abelian variety $\bar{A}$ over a finite
field of characteristic $p$. The Shimura-Taniyama formula determines
$\pi_{\bar{A}}$ up to a root of $1$ in terms of the CM-type of $A$. Using
this, Honda proved that every Weil $q$-number arises from an abelian variety,
possibly after a finite extension of the base field.\footnote{Honda, Taira,
Isogeny classes of abelian varieties over finite fields. J. Math. Soc. Japan
20 1968 83--95.} An application of Weil restriction of scalars completes the proof.

\subsection{Good reduction of Abelian Varieties\label{b7}}

The language of Weil's foundations of algebraic geometry is ill-suited to the
study of algebraic varieties in mixed characteristic. For example, it makes it
cumbersome to prove even that an algebraic variety over a number field has
good reduction at almost all primes of the field.\footnote{See the proof of
Theorem 26 of Shimura, Goro, Reduction of algebraic varieties with respect to
a discrete valuation of the basic field. Amer. J. Math. 77, (1955). 134--176.}
Serre and Tate (1968a) use schemes and N\'{e}ron's theory of minimal
models\footnote{N\'{e}ron, Andr\'{e}. Mod\`{e}les minimaux des
vari\'{e}t\'{e}s ab\'{e}liennes sur les corps locaux et globaux. Inst. Hautes
\'{E}tudes Sci. Publ.Math. No. 21 1964 128 pp.} to simplify and sharpen known
results for abelian varieties, and to extend some statements from elliptic
curves to abelian varieties.

Let $R$ be a discrete valuation ring with field of fractions $K$ and perfect
residue field $k$. For an abelian variety $A$ over $K$, N\'{e}ron proved that
the functor sending a smooth $R$-scheme $X$ to $\Hom(X_{K},A)$ is represented
by a smooth group scheme $\tilde{A}$ of finite type over $R$. Using this,
Serre and Tate prove the following criterion:

\begin{quote}
If $A$ has good reduction, then the $\Gal(K^{\mathrm{sep}}/K)$-module
$A(K^{\mathrm{sep}})_{m}$ is unramified for all integers $m$ prime to
$\mathrm{\mathrm{char}}(k)$; conversely, if $A(K^{\mathrm{sep}})_{m}$ is
unramified for infinitely many $m$ prime to $\mathrm{char}(k)$, then $A$ has
good reduction.
\end{quote}

\noindent The necessity was known earlier, and the sufficiency was known to
Ogg and Shafarevich. in the case of elliptic curves. Because it is a direct
consequence of the existence of N\'{e}ron's models, Serre and Tate call it the
\textquotedblleft N\'{e}ron-Ogg-Shafarevich criterion\textquotedblright. It is
of fundamental importance.

Serre and Tate say that an abelian variety has \emph{potential good reduction}
if it acquires good reduction after a finite extension of the base field, and
they prove a number of results about such varieties. For example, when $R$ is
strictly henselian, there is a smallest extension $L$ of $K$ in
$K^{\mathrm{al}}$ over which such an abelian variety $A$ has good reduction,
namely, the extension of $K$ generated by the coordinates of the points of
order $m$ for any $m\geq3$ prime to $\mathrm{char}(k)$. Moreover, just as for
elliptic curves, the notion of the conductor of an abelian variety is well defined.

Let $(A,i)$ be an abelian variety over a number field $K$ with complex
multiplication by $E$. By this we mean that $E$ is a CM field of degree $2\dim
A$ over $\mathbb{Q}{}$, and that $i$ is a homomorphism of $\mathbb{Q}{}%
$-algebras $E\rightarrow\End^{0}(A)$. Serre and Tate apply their earlier
results to show that such an abelian variety $A$ acquires good reduction
everywhere over a cyclic extension $L$ of $K$; moreover, $L$ can be chosen to
have degree $m$ or $2m$ where $m$ is the least common multiple of the images
of the inertia groups acting on the torsion points of $A$.

Let $(A,i)$ and $E$ be as in the last paragraph, and let $C_{K}$ be the
id\`{e}le class group of $K$. Shimura and Taniyama (1961,
18.3)\footnote{Shimura, Goro; Taniyama, Yutaka Complex multiplication of
abelian varieties and its applications to number theory. Publications of the
Mathematical Society of Japan, 6 The Mathematical Society of Japan, Tokyo
1961.} show there exists a (unique) homomorphism $\rho\colon C_{K}%
\rightarrow(\mathbb{R}{}\otimes_{\mathbb{Q}{}}E)^{\times}$ with the following
property: for each $\sigma\colon E\rightarrow\mathbb{C}{}$, let $\chi_{\sigma
}$ be the Hecke character%
\[
\chi_{\sigma}\colon C_{K}\overset{\rho}{\longrightarrow}(\mathbb{R}{}%
\otimes_{\mathbb{Q}{}}E)^{\times}\overset{1\otimes\sigma}{\longrightarrow
}\mathbb{C}{}^{\times};
\]
then the $L$-series $L(s,A)$ coincides with the product $\prod
\nolimits_{\sigma}L(s,\chi_{\sigma})$ of the $L$-series of the $\chi_{\sigma}%
$, except possibly for the factors corresponding to a finite number of primes
of $K$.\footnote{Rather, this is Serre and Tate's interpretation of what they
prove; Shimura and Taniyama express their results in terms of ideals.} Serre
and Tate make this more precise by showing that the conductor of $A$ is the
product of the conductors of the $\chi_{\sigma}$ (which each equals the
conductor of $\rho$). In particular, the support of the conductor of each
$\chi_{\sigma}$ equals the set of primes where $A$ has bad reduction, from
which it follows that $L(s,A)$ and $\prod\nolimits_{\sigma}L(s,\chi_{\sigma})$
coincide exactly.

\subsection{CM abelian varieties and Hilbert's twelfth problem\label{b8}}

A CM-type on a CM field $E$ is a subset $\Phi$ of $\Hom(E,\mathbb{C}{})$ such
that $\Phi\sqcup\bar{\Phi}=\Hom(E,\mathbb{C}{})$. For $\sigma\in
\Aut(\mathbb{C}{})$, let $\sigma\Phi=\{\sigma\circ\varphi\mid\varphi\in\Phi
\}$. Then $\sigma\Phi$ is also a CM-type on $E$. The reflex field of
$(E,\Phi)$ is the subfield $F$ of $\mathbb{C}{}$ such that an automorphism
$\sigma$ of $\mathbb{C}{}$ fixes $F$ if and only if $\sigma\Phi=\Phi$. It is
easy to see that $F$ is a CM-subfield of $\mathbb{Q}{}^{\mathrm{al}}%
\subset\mathbb{C}{}$.

Let $(A,i)$ be an abelian variety over $\mathbb{C}{}$ with complex
multiplication by $E$. Then $E$ acts on the tangent space of $A$ at $0$
through a CM-type $\Phi$, and $(A,i)$ is said to be of CM-type $(E,\Phi)$. For
$\sigma\in\Aut(\mathbb{C}{}/\mathbb{Q}{})$, $\sigma(A,i)$ is of CM-type
$\sigma\Phi$, and it follows that $\sigma(A,i)$ is isogenous to $(A,i)$ if and
only if $\sigma$ fixes the reflex field $F$. Fix a polarization $\lambda$ of
$A$ whose Rosati involution acts as complex multiplication on $E$. For an
integer $m\geq1$, let $\mathcal{S}{}(m)$ be the set of isomorphism classes of
quadruples $(A^{\prime},\lambda^{\prime},i^{\prime},\eta)$ such that
$(A^{\prime},\lambda^{\prime},i^{\prime})$ is isogenous to $(A,\lambda,i)$ and
$\eta$ is a level $m$-structure on $(A^{\prime},i^{\prime})$. According to the
preceding observation, $\Aut(\mathbb{C}{}/F)$ acts on the set $\mathcal{S}%
{}(m)$. Shimura and Taniyama prove that this action factors through
$\Aut(F^{\mathrm{ab}}/F)$, and they describe it explicitly. In this way, they
generalized the theory of complex multiplication from elliptic curves to
abelian varieties, and they provided a partial solution to Hilbert's twelfth
problem for $F$.

In one respect the result of Shimura and Taniyama falls short of generalizing
the elliptic curve case: for an elliptic curve, the reflex field $F$ is a
complex quadratic extension of $\mathbb{Q}{}$; since one knows how complex
conjugation acts on CM elliptic curves and their torsion points, the elliptic
curve case provides a description of how the \emph{full} group
$\Aut(\mathbb{C}{}/\mathbb{Q}{})$ acts on CM elliptic curves and their torsion
points. Shimura asked whether there was a similar result for abelian
varieties, but concluded rather pessimistically that \textquotedblleft In the
higher-dimensional case, however, no such general answer seems
possible.\textquotedblright\footnote{Shimura, Goro. On abelian varieties with
complex multiplication. Proc. London Math. Soc. (3) 34 (1977), no. 1, 65--86.}

Grothendieck's theory of motives suggests the framework for an answer. The
Hodge conjecture implies the existence of Tannakian category of CM-motives
over $\mathbb{Q}$, whose motivic Galois group is an extension%
\[
1\rightarrow S\rightarrow T\rightarrow\Gal(\mathbb{Q}{}^{\mathrm{al}%
}/\mathbb{Q}{})\rightarrow1
\]
of $\Gal(\mathbb{Q}{}^{\mathrm{al}}/\mathbb{Q}{})$ (regarded as a pro-constant
group scheme) by the Serre group $S$ (a certain pro-torus). \'{E}tale
cohomology defines a section $\lambda$ of $T\rightarrow\Gal(\mathbb{Q}%
{}^{\mathrm{al}}/\mathbb{Q}{})$ over the finite ad\`{e}les. The pair
$(T,\lambda)$ (tautologically) describes the action of $\Aut(\mathbb{C}%
{}/\mathbb{Q}{})$ on the CM abelian varieties and their torsion points.
Deligne's theorem on Hodge classes on abelian varieties allows one to
construct the pair $(T,\lambda)$ without assuming the Hodge conjecture. To
answer Shimura's question, it remains to give a direct explicit description of
$(T,\lambda)$.

Langlands's work on the zeta functions of Shimura varieties led him to define
a certain explicit cocycle,\footnote{See \S 5 of Langlands, R. P. Automorphic
representations, Shimura varieties, and motives. Ein M\"{a}rchen. Automorphic
forms, representations and $L$-functions (Proc. Sympos. Pure Math., Oregon
State Univ., Corvallis, Ore., 1977), Part 2, pp. 205--246, Proc. Sympos. Pure
Math., XXXIII, Amer. Math. Soc., Providence, R.I., 1979} which Deligne
recognized as conjecturally being that describing the pair $(T,\lambda)$.

Tate was inspired by this to commence his own investigation of Shimura's
question. He gave a simple direct construction of a map $f$ that he
conjectured describes how $\Aut(\mathbb{C}{}/\mathbb{Q}{})$ acts on the CM
abelian varieties and their torsion points, and proved this up to signs. More
precisely, he proved it up to a map $e$ with values in an ad\`{e}lic group
such that $e^{2}=1$. See Tate 1981c.

It was soon checked that Langlands's and Tate's conjectural descriptions of
how $\Aut(\mathbb{C}{}/\mathbb{Q}{})$ acts on the CM abelian varieties and
their torsion points coincided, and a few months later Deligne proved that
their conjectural descriptions are indeed correct.\footnote{Deligne, P.,
Motifs et groupe de Taniyama, pp.261--279 in Hodge cycles, motives, and
Shimura varieties. Lecture Notes in Mathematics, 900. Springer-Verlag,
Berlin-New York, 1982.}

\section{Rigid analytic spaces\label{c}}

After Hensel introduced the $p$-adic number field $\mathbb{Q}{}_{p}$ in the
1890s, there were attempts to develop a theory of analytic functions over
$\mathbb{Q}{}_{p}$, the most prominent being that of Krasner. The problem is
that every disk $D$ in $\mathbb{Q}{}_{p}$ can be written as a disjoint union
of arbitrarily many open-closed smaller disks, and so there are too many
functions on $D$ that can be represented locally by power series. Outside a
small group of mathematicians, $p$-adic analysis attracted little attention
until the work of Dwork and Tate in late 1950s. In February, 1958, Tate sent
Dwork a letter in which he stated a result concerning elliptic curves, and
challenged Dwork to find a proof using $p$-adic analysis. In answering the
letter, Dwork found \textquotedblleft the first suggestion of a connection
between $p$-adic analysis and the theory of zeta functions.\textquotedblright%
\footnote{Katz and Tate, 1999, p.343; Dwork, Bernard. A deformation theory for
the zeta function of a hypersurface. 1963 Proc. Internat. Congr.
Mathematicians (Stockholm, 1962) pp. 247--259 Inst. Mittag-Leffler,
Djursholm.} By November, 1959, Dwork had found his famous proof of the
rationality of the zeta function $Z(V,T)$ of an algebraic variety $V$ over a
finite field, a key point of which is to express $Z(V,T)$, which initially is
a power series with integer coefficients, as a quotient of two $p$-adically
entire functions.\footnote{Dwork, Bernard. On the rationality of the zeta
function of an algebraic variety. Amer. J. Math. 82 1960 631--648.}

In 1959 also, Tate discovered that, suitably normalized, certain classical
formulas allow one to express many elliptic curves $E$ over a nonarchimedean
local field $K$ as a quotient $E(K)=K^{\times}/q^{\mathbb{Z}{}}$. This
persuaded him that there should exist a category in which $E$ itself, not just
its points, is a quotient; in other words, that there exists a category in
which $E$, as an \textquotedblleft analytic space\textquotedblright, is the
quotient of $K^{\times}$, as an \textquotedblleft analytic
space\textquotedblright, by the discrete group $q^{\mathbb{Z}{}}$. Two years
later, Tate constructed the correct category of \textquotedblleft rigid
analytic spaces\textquotedblright, thereby founding a new subject in
mathematics (with its own Math. Reviews number 14G22).

\subsection{The Tate curve\label{c1}}

Let $E$ be an elliptic curve over $\mathbb{C}{}$. The choice of a differential
$\omega$ realizes $E(\mathbb{C}{})$ as the quotient $\mathbb{C}{}%
/\Lambda\simeq E(\mathbb{C}{})$ of $\mathbb{C}{}$ by the lattice of periods of
$\omega$. More precisely, it realizes the complex analytic manifold
$E^{\text{an}}$ as the quotient of the complex analytic manifold $\mathbb{C}%
{}$ by the action of the discrete group $\Lambda$.

For an elliptic curve $E$ over a $p{}$-adic field $K$, there is no similar
description of $E(K)$ because there are no nonzero discrete subgroups of $K$
(if $\lambda\in K$, then $p^{n}\lambda\rightarrow0$ as $n\rightarrow\infty$).
However, there is an alternative uniformization of elliptic curves over
$\mathbb{C}{}$. Let $\Lambda$ be the lattice $\mathbb{Z}{}+\mathbb{Z}{}\tau$
in $\mathbb{C}{}$. Then the exponential map $\underline{e}\colon\mathbb{C}%
{}\rightarrow\mathbb{C}{}^{\times}$ sends $\mathbb{C}{}/\Lambda$
isomorphically onto $\mathbb{C}{}^{\times}/q^{\mathbb{Z}{}}$ where
$q=\underline{e}(\tau)$, and so $\mathbb{C}{}^{\times}/q^{\mathbb{Z}{}}\simeq
E^{\text{an}}$ (as analytic spaces). If $\im(\tau)>0$, then $|q|<1$, and the
elliptic curve $E_{q}$ is given by the equation%
\begin{equation}
Y^{2}Z+XYZ=X^{3}-b_{2}XZ^{2}-b_{3}Z^{3}\text{,} \label{e15}%
\end{equation}
where
\begin{equation}
\left\{
\begin{aligned} b_{2} & =5\sum\nolimits_{n=1}^{\infty}\frac{n^{3}q^{n}}{1-q^{n}}=5q+45q^{2}+140q^{3}+\cdots\\ b_{3} & =\sum\nolimits_{n=1}^{\infty}\frac{7n^{5}+5n^{3}}{12}\frac{q^{n}}{1-q^{n}}=q+23q^{2}+154q^{3}+\cdots\label{e16} \end{aligned}\right.
\end{equation}
are power series with integer coefficients. The discriminant and modular
invariant of $E_{q}$ are given by the usual formulas%
\begin{align}
\Delta &  =q\prod\nolimits_{n\geq1}(1-q^{n})^{24}\label{e17}\\
j(E_{q})  &  =\frac{(1+48b_{2})^{3}}{q\prod\nolimits_{n\geq1}(1-q^{n})^{24}%
}=\frac{1}{q}+744+196884q+\cdots. \label{e20}%
\end{align}

Now let $K$ be a field complete with respect to a nontrivial nonarchimedean
valuation with residue field of characteristic $p\neq0$, and let $q$ be an
element of $K^{\times}$ with $|q|<1$. The series (\ref{e16}) converge in $K$,
and Tate discovered\footnote{\textquotedblleft I still remember the thrill and
amazement I felt when it occurred to me that the classical formulas for such
an isomorphism over $\mathbb{C}{}$ made sense $p$-adically when properly
normalized.\textquotedblright\ Tate 2008.} that (\ref{e15}) is an elliptic
curve $E_{q}$ such that $K^{\prime}/q^{\mathbb{Z}{}}\simeq E_{q}(K^{\prime})$
for all finite extension $K^{\prime}$ of $K$. It follows from certain power
series identities, valid over $\mathbb{Z}{}$, that the discriminant and
modular invariant of $E_{q}$ are given by (\ref{e17}) and (\ref{e20}). Every
$j\in K^{\times}$ with $|j|<1$ arises from a $q$ (determined by (\ref{e20}),
which allows $q$ to be expressed as a power series in $1/j$ with integer
coefficients). The function field $K(E_{q})$ of $E_{q}$ consists of the
quotients $F/G$ of Laurent series%
\[
F=\sum\nolimits_{-\infty}^{\infty}a_{n}z^{n},\quad G=\sum\nolimits_{-\infty
}^{\infty}b_{n}z^{n},\quad a_{n},b_{n}\in K,
\]
converging for all nonzero $z$ in $\mathbb{C}{}_{p}$, such that the $F/G$ is
invariant under $q^{\mathbb{Z}{}}$:
\[
F(qz)/G(qz)=F(z)/G(z).
\]
The elliptic curves $E$ over $K$ with $|j(E)|<1$ that arise in this way are
exactly those whose reduced curve has a node with tangents that are rational
over the base field. They are now called Tate (elliptic) curves%
\index{Tate (elliptic) curve}%
.

Tate's results were contained in a 1959 manuscript, which he did not publish
until 1995, but there soon appeared several summaries of his results in the
literature, and Roquette\footnote{Roquette, Peter, Analytic theory of elliptic
functions over local fields. Hamburger Mathematische Einzelschriften (N.F.),
Heft 1 Vandenhoeck \& Ruprecht, G\"{o}ttingen 1970} gave a very detailed
account of the theory. The Tate curve has found many applications, for
example, to Tate's isogeny conjecture (Serre 1968;\footnote{Serre,
Jean-Pierre. Abelian $l$-adic representations and elliptic curves. W. A.
Benjamin, Inc., New York-Amsterdam 1968.} Tate 1995, p.180) and to the study
of elliptic modular curves near a cusp (Deligne and Rapoport\footnote{Deligne,
P.; Rapoport, M. Les sch\'{e}mas de modules de courbes elliptiques. Modular
functions of one variable, II (Proc. Internat. Summer School, Univ. Antwerp,
Antwerp, 1972), pp. 143--316. Lecture Notes in Math., Vol. 349, Springer,
Berlin, 1973.}). Mumford\footnote{Mumford, David. An analytic construction of
degenerating curves over complete local rings. Compositio Math. 24 (1972),
129--174.} generalized Tate's construction to curves of higher genus, and
McCabe\footnote{McCabe, John. $p$-adic theta functions. Ph.D. thesis, Harvard,
1968, 222 pages.} and Raynaud\footnote{Raynaud, Michel. Vari\'{e}t\'{e}s
ab\'{e}liennes et g\'{e}om\'{e}trie rigide. Actes du Congr\`{e}s International
des Math\'{e}maticiens (Nice, 1970), Tome 1, pp. 473--477. Gauthier-Villars,
Paris, 1971.} generalized it to abelian varieties of higher dimension.

\subsection{Rigid analytic spaces\label{c2}}

Tate's idea that his $p$-adic uniformization of elliptic curves indicated the
existence of a general theory of $p$-adic analytic spaces was radically new.
For example, Grothendieck was initially very
negative.\footnote{\textquotedblleft Tate has written to me about his elliptic
curve stuff, and has asked me if I had any ideas for a global definition of
analytic varieties over complete valuation fields. I must admit that I have
absolutely not understood why his results might suggest the existence of such
a definition, and I remain skeptical. Nor do I have the impression of having
understood his theorem at all; it does nothing more than exhibit, via brute
formulas, a certain isomorphism of analytic groups.\textquotedblright%
\ Grothendieck, letter to Serre, August 18, 1959.} However, when Tate began to
work out his theory in the fall of $1961$, Grothendieck, who was visiting
Harvard at the time, became very optimistic,\footnote{\textquotedblleft Sooner
or later it will be necessary to subsume ordinary analytic spaces, rigid
analytic spaces, formal schemes, and maybe even schemes themselves into a
single kind of structure for which all these usual theorems will
hold.\textquotedblright\ Grothendieck, letter to Serre, October 19, 1961.} and
was very supportive.

Let $K$ be a field complete with respect to a nontrivial nonarchmedean
valuation, and let $\bar{K}$ be its algebraic closure. Tate began by
introducing a new class of $K$-algebras. The Tate algebra%
\index{Tate algebra}
$T_{n}=K\{X_{1},\ldots,X_{n}\}$ consists of the formal power series in
$K[[X_{1},\ldots,X_{n}]]$ that are convergent on the unit ball,
\[
B^{n}=\left\{  (c_{i})_{1\leq i\leq n}\in\bar{K}\mid\left\vert c_{i}%
\right\vert \leq1\right\}  .
\]
Thus the elements of $T_{n}$ are the power series%
\[
f=\sum a_{i_{1}\cdots i_{n}}X_{1}^{i_{1}}\cdots X_{n}^{i_{n}},\quad
a_{i_{1}\cdots i_{n}}\in K,\quad\text{such that }a_{i_{1}\cdots i_{n}%
}\rightarrow0\text{ as }(i_{1},\ldots,i_{n})\rightarrow\infty.
\]
Tate (1962c) shows that $T_{n}$ is a Banach algebra for the norm $\left\Vert
f\right\Vert =\sup\left\vert a_{i_{1}\cdots i_{n}}\right\vert $, and that the
ideals $\mathfrak{a}{}$ of $T_{n}$ are closed and finitely generated. A
quotient $T_{n}/\mathfrak{a}{}$ of $T_{n}$ is a Banach algebra whose topology
is independent of its presentation (because every homomorphism of such
algebras is continuous). Such quotients are called affinoid (or Tate)
$K$-algebras,%
\index{Tate algebra}
and the category of affine rigid analytic spaces is the opposite of the
category of affinoid $K$-algebras.

We need a geometric interpretation of this category. Tate showed that $T_{n}$
is Jacobson (i.e., every prime ideal is an intersection of maximal ideals),
and that the map $A\mapsto\mathrm{max}(A)$ sending an affinoid algebra to its
set of maximal ideals is a functor: a homomorphism $\varphi\colon A\rightarrow
B$ of affinoid algebras defines a map $\varphi^{\circ}\colon\max
(B)\rightarrow\max(A)$. The set $\mathrm{max}(A)$ has the Zariski topology,
which is very coarse, and a canonical topology induced from that of $K$. When
$K$ is algebraically closed, $\mathrm{max}(T_{n})\simeq B^{n}$, and, by
definition, $\mathrm{max}(A)$ can be realized as a closed subset of
$\mathrm{max}(T_{n})$ for some $n$.

Let $X=\mathrm{max}(A)$. One would like to define a sheaf $\mathcal{O}{}_{X}$
on $X$ such that, for every open subset $U$ isomorphic to $B^{n}$,
$\mathcal{O}{}_{X}(U)\simeq T_{n}$. As noted at the start of this section,
this is impossible. However, Tate's realized that it is possible to achieve
something like this by allowing only certain \textquotedblleft
admissible\textquotedblright\ open subsets and certain \textquotedblleft
admissible\textquotedblright\ coverings. He defined an \emph{affine subset} of
$X$ to be a subset $Y$ such that the functor of affinoid $K$-algebras%
\[
B\rightsquigarrow\left\{  \varphi\colon A\rightarrow B\mid\varphi^{\circ
}(\mathrm{max}(B))\subset Y\right\}
\]
is representable (say, by $A\rightarrow A(Y)$). A subset $Y$ of $X$ is a
\emph{special affine subset} of $X$ if there exist two finite families
$(f_{i})$ and $(g_{j})$ of elements of $A$ such that%
\[
Y=\left\{  x\in X\,\middle|\,\left\vert f_{i}(x)\right\vert \leq
1,\quad\left\vert g_{j}(x)\right\vert \geq1,\quad\text{all }i,j\right\}  .
\]
Every special affine subset is affine. Tate's acyclicity theorem%
\index{Tate's acyclicity theorem}
(Tate 1962c, 8.2) says that, for every finite covering $(X_{i})_{i\in I}$ of
$X$ by special affines, the \u{C}ech complex of the presheaf $Y\mapsto A(Y)$,
\[
0\rightarrow A\rightarrow\prod_{i_{0}}A(X_{i_{0}})\rightarrow\prod
_{i_{0}<i_{1}}A(X_{i_{0}}\cap X_{i_{1}})\rightarrow\cdots\rightarrow
\prod_{i_{0}<\cdots<i_{p}}A(X_{i_{0}}\cap\cdots\cap X_{i_{p}})\rightarrow
\cdots,
\]
is exact. In particular, $Y\mapsto A(Y)$ satisfies the sheaf condition on such coverings.

Using Tate's acyclicity theorem it is possible to define a collection of
admissible open subsets of $X=\mathrm{max(}A)$ and admissible coverings of
them for which there exists a functor $\mathcal{O}{}_{X}$ satisfying the sheaf
conditions and such that $\mathcal{O}{}_{X}(Y)=A(Y)$ for any affine subset.
Although the admissible open subsets and coverings don't form a topology in
the usual sense, they satisfy the conditions necessary for them to support a
sheaf theory --- in fact, they form a Grothendieck topology. So, in this
sense, Tate recovers analytic continuation.

For the final step, extending the category of affine rigid analytic spaces to
a category of global rigid analytic spaces, Tate followed suggestions of
Grothendieck. This step has since been clarified and simplified; see, for
example, Bosch 2005,\footnote{Bosch, S., Lectures on Formal and Rigid
Geometry, Preprint 378 of the SFB Geometrische Strukturen in der Mathematik,
M\'{u}nster, 2005.} especially 1.12.

Tate reported on his work in a series of letters to Serre, who had them typed
by IHES as the notes Tate 1962c. These notes were distributed to a number of
mathematicians and libraries. They soon attracted the attention of the German
school of complex analytic geometers, who were able to transfer many of their
arguments and results to the new setting (e.g., Kiehl 1967\footnote{Kiehl,
Reinhardt. Der Endlichkeitssatz f\"{u}r eigentliche Abbildungen in der
nichtarchimedischen Funktionentheorie. Invent. Math. 2 1967 191--214; Theorem
A und Theorem B in der nichtarchimedischen Funktionentheorie. Ibid.
256--273.}). Already by 1984 to give a comprehensive account of the theory
required a book of over 400 pages (Bosch, G\"{u}ntzer, Remmert
1984\footnote{Bosch, S.; G\"{u}ntzer, U.; Remmert, R. Non-Archimedean
analysis. A systematic approach to rigid analytic geometry. Grundlehren der
Mathematischen Wissenschaften, 261. Springer-Verlag, Berlin, 1984.}). Tate did
not publish his work, but eventually the editors of \textquotedblleft
Mir\textquotedblright\ published a Russian translation of his notes (Tate
1969a), and the editors of \textquotedblleft Inventiones\textquotedblright%
\ published the original (Tate 1971).

There have been a number of extensions of Tate's theory. For example,
following a suggestion of Grothendieck, Raynaud showed that it is possible to
realize a rigid analytic space over a field $K$ as the \textquotedblleft
generic fibre\textquotedblright\ of a formal scheme over the valuation ring of
$K$. One problem with rigid analytic spaces is that, while they are adequate
for the study of coherent sheaves, they have too few points for the study of
locally constant sheaves --- for example, there exist nonzero such sheaves
whose stalks are all zero. Berkovich found a solution to this problem by
enlarging the underlying set of a rigid analytic space without altering the
sheaf of functions so that the spaces now support an \'{e}tale cohomology
theory (Berkovich 1990;\footnote{Berkovich, Vladimir G. Spectral theory and
analytic geometry over non-Archimedean fields. Mathematical Surveys and
Monographs, 33. American Mathematical Society, Providence, RI, 1990.}
1993\footnote{Berkovich, Vladimir G. Etale cohomology for non-Archimedean
analytic spaces. Inst. Hautes \'{E}tudes Sci. Publ. Math. No. 78 (1993),
5--161.}).

Rigid analytic spaces are now part of the landscape of arithmetic geometry:
just as it is natural to regard the $\mathbb{R}{}$-points of a $\mathbb{Q}{}%
$-variety as a real analytic space, it has become natural to regard the
$\mathbb{Q}{}_{p}$-points of the variety as a rigid analytic space. They have
found numerous applications, for example, in the solution by Harbater and
Raynaud of Abyhankar's conjecture on the \'{e}tale fundamental groups of
curves, and in the Langlands program (see \ref{ee1}).

\section{The Tate conjecture\label{d}}

\hfill\begin{minipage}{2.5in}
\textit{This stuff is too beautiful not to be true}\\
\hspace*{2in}Tate\footnotemark
\end{minipage}\footnotetext{As a thesis topic, Tate gave me the problem of
proving a formula that he and Mike Artin had conjectured concerning algebraic
surfaces over finite fields (Conjecture C below). One day he ran into me in
the corridors of 2 Divinity Avenue and asked how it was going.
\textquotedblleft Not well\textquotedblright\ I said, \textquotedblleft In one
example, I computed the left hand side and got $p^{13}$; for the other side, I
got $p^{17}$; $13$ is not equal to $17$, and so the conjecture is
false.\textquotedblright\ For a moment, Tate was taken aback, but then he
broke into a grin and said \textquotedblleft That's great! That's really
great! Mike and I must have overlooked some small factor which you have
discovered.\textquotedblright\ He took me off to his office to show him. In
writing it out in front of him, I discovered a mistake in my work, which in
fact proved that the conjecture is correct in the example I considered. So I
apologized to Tate for my carelessness. But Tate responded: \textquotedblleft
Your error was not that you made a mistake --- we all make mistakes. Your
error was not realizing that you must have made a mistake. This stuff is too
beautiful not to be true.\textquotedblright\ Benedict Gross tells of a similar
experience, but as he writes: \textquotedblleft John was so encouraging,
saying that everyone made mistakes, and the key was to understand them and to
keep thinking about the problem. I felt that one of his greatest talents as an
advisor was to make his students feel like we were partners in a great
enterprise, modern number theory.\textquotedblright
\par
{}}\bigskip

The Hodge conjecture says that a rational cohomology class on a nonsingular
projective variety over $\mathbb{C}{}$ is algebraic if it is of type $(p,p)$.
The Tate conjecture says that an $\ell$-adic cohomology class on a nonsingular
projective variety over a finitely generated field $k$ is in the span of the
algebraic classes if it is fixed by the Galois group. (A field is finitely
generated if it is finitely generated as a field over its prime field.)

\subsection{Beginnings \label{d0}}

In the last section of his talk at the 1962 International Congress, Tate
states several conjectures.

\begin{E}
\label{d0a}For every abelian variety $A$ over a global field $k$ and prime
$\ell\neq\mathrm{char}(k)$, $\TS(A/k)(\ell)$ is finite.
\end{E}

\noindent\begin{minipage}{4in}
Let $k$ be a global function field, and let $k_{0}$ be its finite field of constants,
so that $k=k_{0}(C)$ for a complete nonsingular curve $C$ over $k_{0}$. An
elliptic curve $A$ over $k$ is the generic fibre of a map $X\rightarrow C$
with $X$ a complete nonsingular surface over $k_{0}$, which may be taken to
be a minimal. Tate showed that, in this case,
(\ref{d0a}) is equivalent to the following conjecture.
\end{minipage}\hfill\begin{minipage}{1.4in}
\begin{tikzpicture}
\matrix(m)[matrix of math nodes, row sep=3em, column sep=2.5em,
text height=1.5ex, text depth=0.25ex]
{X&A\\
C&\Spec(k)\\};
\path[->,>=angle 90]
(m-1-2) edge  (m-1-1)
(m-1-1) edge  (m-2-1)
(m-2-2) edge  (m-2-1)
(m-1-2) edge  (m-2-2);
\end{tikzpicture}
\end{minipage}

\begin{E}
\label{d0b}Let $q=\left\vert k_{0}\right\vert $. The $\mathbb{Z}{}_{\ell}%
$-submodule of $H_{\mathrm{et}}^{2}(X_{k_{0}^{\mathrm{al}}},\mathbb{Z}{}%
_{\ell})$ on which the Frobenius map acts as multiplication by $q$ is exactly
the submodule generated by the algebraic classes.
\end{E}

As Tate notes, (\ref{d0b}) makes sense for any complete nonsingular surface
over $k_{0}$, and that, so generalized, it is equivalent to the following
statement.\footnote{Assuming the Weil conjectures, which weren't proved until
1973.}

\begin{E}
\label{d0c}Let $X$ be a complete nonsingular surface over a finite field. The
order of the pole of $\zeta(X,s)$ at $s=1$ is equal to the number of
algebraically independent divisors on $X$.
\end{E}

Mumford pointed out that (\ref{d0c}) implies that elliptic curves over a
finite field are isogenous if and only if they have the same zeta function,
and he proved this using results of Deuring\footnote{Deuring, Max. Die Typen
der Multiplikatorenringe elliptischer Funktionenk\"{o}rper. Abh. Math. Sem.
Hansischen Univ. 14, (1941). 197--272.} on the lifting to characteristic $0$
of the Frobenius automorphism.

In his talk at the 1964 Woods Hole conference, Tate vastly generalized these conjectures.

\subsection{Statement of the Tate conjecture\label{d1}}

For a connected nonsingular projective variety $V$ over a field $k$, we let
$\mathcal{Z}{}^{r}(V)$ denote the $\mathbb{Q}{}$-vector space of algebraic
cycles on $V$ of codimension $r$, i.e., the $\mathbb{Q}{}$-vector space with
basis the irreducible closed subsets of $V$ of dimension $\dim V-r$. We let
$H_{\mathrm{et}}^{r}(V,\mathbb{Q}{}_{\ell}(s))$ denote the \'{e}tale
cohomology group of $V$ with coefficients in the \textquotedblleft Tate
twist\textquotedblright\ $\mathbb{Q}{}_{\ell}(s)$ of $\mathbb{Q}{}_{\ell}$.
There are cycle maps%
\[
c^{r}\colon\mathcal{Z}{}^{r}(V)\rightarrow H_{\mathrm{et}}^{2r}(V,\mathbb{Q}%
{}_{\ell}(r)).
\]

Assume that $\ell\neq\mathrm{char}(k)$. Let $\bar{k}$ be an algebraically
closed field containing $k$, and let $G(\bar{k}/k)$ be the group of
automorphisms of $\bar{k}$ fixing $k$. Then $G(\bar{k}/k)$ acts on
$H_{\mathrm{et}}^{2r}(V_{\bar{k}},\mathbb{Q}{}_{\ell}(r))$, and the \emph{Tate
conjecture}\footnote{In the literature, a number of variants of $T^{r}(V)$,
not obviously equivalent to it, are also called the Tate conjecture. It is not
always easy to discern what an author means by the \textquotedblleft Tate
conjecture\textquotedblright.} (Tate 1964a, Conjecture 1) is the following
statement:\footnote{Since Atiyah and Hirzebruch had already found their
counterexample to an integral Hodge conjecture, Tate was not tempted to state
his conjecture integrally.}

\begin{quote}
$T^{r}(V)$: When $k$ is finitely generated, the $\mathbb{Q}{}_{\ell}$-space
spanned by $c^{r}(\mathcal{Z}{}^{r}(V_{\bar{k}}))$ consists of the elements of
$H_{\mathrm{et}}^{2r}(V_{\bar{k}},\mathbb{Q}{}_{\ell}(r))$ fixed by some open
subgroup of $G(\bar{k}/k)$.
\end{quote}

\noindent Suppose for simplicity that $\bar{k}$ is an algebraic closure of
$k$. For any finite extension $k^{\prime}$ of $k$ in $\bar{k}$, there is a
commutative diagram%
\[
\begin{CD}
\mathcal{Z}^{r}(V_{\bar{k}}) @>{c^{r}}>> H_{\mathrm{et}}^{2r}(V_{\bar{k}},\mathbb{Q}%
_{\ell}(r)){}\\
@AAA@AAA\\
\mathcal{Z}^{r}(V_{k^{\prime}}) @>{c^{r}}>> H_{\mathrm{et}}^{2r}(V_{k^{\prime}},\mathbb{Q}%
{}_{\ell}(r)),
\end{CD}
\]
and the image of the right hand map is $H_{\mathrm{et}}^{2r}(V_{\bar{k}%
},\mathbb{Q}_{\ell}(r))^{G(\bar{k}/k^{\prime})}$. As $\mathcal{Z}{}%
^{r}(V_{\bar{k}})=\bigcup\nolimits_{k^{\prime}}\mathcal{Z}{}^{r}(V_{k^{\prime
}})$, we see that%
\[
c^{r}(\mathcal{Z}{}^{r}(V_{\bar{k}}))\subset\bigcup\nolimits_{k^{\prime}%
}H_{\mathrm{et}}^{2r}(V_{\bar{k}},\mathbb{Q}_{\ell}(r))^{G(\bar{k}/k^{\prime
})}\text{.}%
\]
The content of the Tate conjecture is that the first set spans the space on
the right. If an element of $\mathcal{Z}{}^{r}(V_{\bar{k}})$ is fixed by
$G(\bar{k}/k^{\prime})$, then it lies in $\mathcal{Z}{}^{r}(V_{k^{\prime}})$,
and so $T^{r}(V)$ implies that%
\begin{equation}
c^{r}(\mathcal{Z}^{r}(V_{k^{\prime}}))\mathbb{Q}{}_{\ell}=H_{\mathrm{et}%
}^{2r}(V_{\bar{k}},\mathbb{Q}_{\ell}(r))^{G(\bar{k}/k^{\prime})}; \label{e14}%
\end{equation}
conversely, if (\ref{e14}) holds for all (sufficiently large) $k^{\prime}$,
then $T^{r}(V)$ is true.\medskip

When asked about the origin of the Tate conjecture, Tate responded (Tate 2011):

\begin{quote}
Early on I somehow had the idea that the special case about endomorphisms of
abelian varieties over finite fields might be true. A bit later I realized
that a generalization fit perfectly with the function field version of the
Birch and Swinnerton-Dyer conjecture. Also it was true in various particular
examples which I looked at and gave a heuristic reason for the Sato-Tate
distribution. So it seemed a reasonable conjecture.
\end{quote}

\noindent I discuss each of these motivations in turn.

\subsection{Homomorphisms of abelian varieties\label{d2}}

Let $A$ be an abelian variety over a field $k$, let $\bar{k}$ be an
algebraically closed field containing $k$, and let $G(\bar{k}/k)$ denote the
group of automorphisms of $\bar{k}$ over $k$. For $\ell\neq\mathrm{char}(k)$,
\[
A\rightsquigarrow T_{\ell}A=\varprojlim A(\bar{k})_{\ell^{n}}%
\]
is a functor from abelian varieties over $k$ to $\mathbb{Z}_{\ell}$-modules
equipped with an action of $G(\bar{k}/k)$. The (Tate) \textit{ isogeny
conjecture}%
\index{Tate isogeny conjecture}
is the following statement:

\begin{quote}
$H(A,B)$: For abelian varieties $A,B$ over a finitely generated field $k$, the
canonical map%
\[
\mathbb{Z}_{\ell}\otimes\Hom(A,B)\rightarrow\Hom(T_{\ell}A,T_{\ell}%
B)^{G(\bar{k}/k)}%
\]
is an isomorphism.
\end{quote}

\noindent It follows from Weil's theory of correspondences and the
interpretation of divisorial correspondences as homomorphisms, that, for
varieties $V$ and $W$,%
\begin{equation}
\mathrm{NS}(V\times W)\simeq\mathrm{NS}(V)\oplus\mathrm{NS}(W)\oplus\Hom(A,B)
\label{e41}%
\end{equation}
where $A$ is the Albanese variety of $V$, $B$ is the Picard variety of $W$,
and $\mathrm{NS}$ denotes the N\'{e}ron-Severi group. On comparing (\ref{e41})
with the decomposition of $H^{2}(V\times W,\mathbb{Q}{}(1))$ given by the
K\"{u}nneth formula, we find that, for varieties $V$ and $W$ over a finitely
generated field $k$,%
\begin{equation}
T^{1}(V\times W)\iff T^{1}(V)+T^{1}(W)+H(A,B)\text{.} \label{e43}%
\end{equation}
When $V$ is a curve, $T^{1}(V)$ is obviously true, and so, for elliptic curves
$E$ and $E^{\prime}$,%
\[
T^{1}(E\times E^{\prime})\iff H(E,E^{\prime}).
\]

At the time Tate made his conjecture, $H(E,E^{\prime})$ was known for elliptic
curves over a finite field as a consequence of work of Deuring (see above),
and $H(E,E)$ was known for elliptic curves over number fields with at least
one real prime (Serre 1964).\footnote{Serre, Jean-Pierre, Groupes de Lie
l-adiques attach\'{e}s aux courbes elliptiques. Colloque de Clermond-Ferrand
1964, Les Tendances G\'{e}om. en Alg\'{e}bre et Th\'{e}orie des Nombres pp.
239--256 \'{E}ditions du Centre National de la Recherche Scientifique, Paris,
1966.}

Tate (1966b) proved $H(A,B)$ for all abelian varieties over finite fields (see
below). As we discussed in (\ref{b6}), this has implication for the
classification of abelian varieties over finite fields (and even cryptography).

Zarhin extended Tate's result to fields finitely generated over $\mathbb{F}%
{}_{p}$, and Faltings proved $H(A,B)$ for all abelian varieties over number
fields in the same article in which he proved Mordell's conjecture. In fact,
$H(A,B)$ has now been proved in all generality (Faltings et al.
1994).\footnote{Rational points. Papers from the seminar held at the
Max-Planck-Institut f\"{u}r Mathematik, Bonn, 1983/1984. Edited by Gerd
Faltings and Gisbert W\"{u}stholz. Aspects of Mathematics, E6. Friedr. Vieweg
\& Sohn, Braunschweig, 1984.}

Tate's theorem proves that $T^{1}$ is true for surfaces over finite fields
that are a product of curves (by (\ref{e43})). When Artin and Swinnerton-Dyer
(1973)\footnote{Artin, M.; Swinnerton-Dyer, H. P. F. , The Shafarevich-Tate
conjecture for pencils of elliptic curves on K3 surfaces. Invent. Math. 20
(1973), 249--266.} proved $T^{1}$ for elliptic $K3$ surfaces over finite
fields, there was considerable optimism that $T^{1}$ would soon be proved for
all surfaces over finite fields. However, there has been little progress in
the years since then. By contrast, the Hodge conjecture is easily proved for divisors.

\subsubsection{Tate's proof of $H(A,B)$ over a finite field}

It suffices to prove the statement with $A=B$. As the map
\[
\mathbb{Z}_{\ell}\otimes\End(A)\rightarrow\End(T_{\ell}A)^{G(\bar{k}/k)}%
\]
is injective, the problem is to construct enough endomorphisms of $A$. I
briefly outline Tate's proof.

(a) If $H(A,A)$ is true for one prime $\ell\neq\mathrm{char}(k)$, then it is
true for all. This allows Tate to choose an $\ell$ that is well adapted to his arguments.

(b) A polarization on $A$ defines a skew-symmetric pairing $V_{\ell}A\times
V_{\ell}A\rightarrow\mathbb{Q}_{\ell}$. Let $W$ be a maximal isotropic
subspace of $V_{\ell}A$ that is stable under $G(\bar{k}/k)$, and let%
\[
X_{n}=(T_{\ell}A\cap W)+\ell^{n}T_{\ell}A.
\]
There is an infinite sequence of isogenies%
\[
\cdots\rightarrow B_{n}\rightarrow B_{n-1}\rightarrow\cdots\rightarrow
B_{1}\rightarrow B_{0}=A
\]
such that the image of $T_{\ell}B_{n}$ in $T_{\ell}A$ is $X_{n}$. Using a
theorem of Weil, Tate shows that each $B_{n}$ has a polarization \emph{of the
same degree} as the original polarization on $A$. As $k$ is finite, this
implies that the $B_{n}$ fall into finitely many isomorphism classes. An
isomorphism $B_{n}\rightarrow B_{n^{\prime}},$ $n\neq n^{\prime}$, gives a
nontrivial isogeny $A\rightarrow A$.

(c) Having constructed one endomorphism of $A$ not in $\mathbb{Z}$, Tate makes
adroit use of the semisimplicity of the rings involved (and his choice of
$\ell$) to complete the proof.

\subsection{Relation to the conjectures of Birch and Swinnerton-Dyer\label{d3}%
}

The original conjectures of Birch and Swinnerton-Dyer were stated for elliptic
curves over $\mathbb{Q}{}$. Tate re-stated them more generally (see \ref{b3}).

\begin{quote}
(A) For an abelian variety $A$ over a global field $K$, the function $L(s,A)$
has a zero of order $r=\mathrm{\rank}A(K)$ at $s=1$.

(B) Moreover,%
\[
L^{\ast}(s,A)\quad\sim\quad\frac{|\TS(A)|\cdot\left\vert D\right\vert
}{\left\vert A(K)_{\mathrm{tors}}\right\vert \cdot\left\vert A^{\prime
}(K)_{\mathrm{tors}}\right\vert }(s-1)^{r}\text{ as }s\rightarrow1.
\]

\end{quote}

Let $f\colon V\rightarrow C$ be a proper map with fibres of dimension $1$,
where $V$ (resp. $C$) is a nonsingular projective surface (resp. curve) over a
finite field $k$. The generic fibre of $f$ is a curve over the global field
$k(C)$, and we let $A(f)$ denote its Jacobian variety (an abelian variety over
$k(C)$). A comparison of the invariants of $V$ with the invariants of $A(f)$
yields the following statement:%
\[
\text{Conjecture }T^{1}\text{ holds for }V\text{ }\iff\text{ Conjecture (A)
holds for }A(f)\text{.}%
\]
\noindent\noindent In examining the situation further, Artin and Tate (Tate
1966e) were led to make the following (Artin-Tate) conjecture:%
\index{Artin-Tate conjecture}%

\begin{quote}
(C) For a projective smooth geometrically-connected surface $V$ over a finite
field $k$, the Brauer group $\Br(V)$ of $V$ is finite, and%
\[
P_{2}(q^{-s})\quad\sim\quad\frac{\left\vert \Br(V)\right\vert \cdot\left\vert
D\right\vert }{q^{\alpha(X)}\left\vert \mathrm{NS}(V)_{\mathrm{tors}%
}\right\vert ^{2}}(1-q^{1-s})^{\rho(V)}\quad\text{as}\quad s\rightarrow1
\]
where $P_{2}(T)$ is the characteristic polynomial of the Frobenius
automorphism acting on $H^{2}(V_{k^{\mathrm{al}}},\mathbb{Q}_{\ell}){}$, $D$
is the discriminant of the intersection pairing on $\mathrm{NS}(V)$, $\rho(V)$
is the rank of \negthinspace\textrm{NS}$(V)$, and $\alpha(V)=\chi
(X,\mathcal{O}{}_{X})-1+\dim($\textrm{PicVar}$(V)).$
\end{quote}

\noindent Naturally, they also conjectured:

\begin{quote}
(d) Let $f\colon V\rightarrow C$ be a proper map, as above, and assume that
$f$ has connected geometric fibres and a smooth generic fibre. Then Conjecture
(B) holds for $A(f)$ over $k(C)$ if and only if Conjecture (C) holds for $V$
over $k$.
\end{quote}

\noindent Tate explains that he gave this conjecture \textquotedblleft only a
small letter (d) as label, because it is of a more elementary nature than (B)
and (C)\textquotedblright, and indeed, it has been proved. Artin and Tate
checked it directly when $f$ is smooth and has a section, Gordon
(1979)\footnote{Gordon, W. J. Linking the conjectures of Artin-Tate and
Birch--Swinnerton-Dyer. Compositio Math. 38 (1979), no. 2, 163--199.} checked
it when the generic fibre has a rational cycle of degree $1$, and Milne
(1982)\footnote{Comparison of the Brauer group with the Tate-Shafarevich
group, J. Fac. Sci. Univ. Tokyo (Shintani Memorial Volume) IA 28 (1982),
735-743.} checked it when this condition holds only locally. However,
ultimately the proof of (d) came from a different direction, by combining the
following two statements:

\begin{itemize}
\item Conjecture C holds for a surface $V$ over a finite field if and only if
$\Br(V)(\ell)$ is finite for some prime $\ell$ (Tate 1966e ignoring the $p$
part; Milne 1975\footnote{Milne, J. S., On a conjecture of Artin and Tate.
Ann. of Math. (2) 102 (1975), no. 3, 517--533.} complete statement);

\item Conjecture A holds for an abelian variety $A$ over a global field of
nonzero characteristic if and only if $\TS(A)(\ell)$ is finite for some $\ell$
(Kato and Trihan 2003\footnote{Kato, Kazuya; Trihan, Fabien, On the
conjectures of Birch and Swinnerton-Dyer in characteristic $p>0$. Invent.
Math. 153 (2003), no. 3, 537--592.}).
\end{itemize}

\noindent In the situation of (d), $\Br(V)(\ell)$ is finite if and only if
$\TS(A(f))(\ell)$ is finite.\footnote{Cf. Liu, Qing; Lorenzini, Dino; Raynaud,
Michel, On the Brauer group of a surface. Invent. Math. 159 (2005), no. 3,
673--676.}

The known cases of Conjecture B over function fields have proved useful in the
construction of lattice packings.\footnote{\textquotedblleft One of the most
exciting developments has been Elkies' (sic) and Shioda's construction of
lattice packings from Mordell-Weil groups of elliptic curves over function
fields. Such lattices have a greater density than any previously known in
dimensions from about 54 to 4096.\textquotedblright\ Preface to Conway, J. H.;
Sloane, N. J. A. Sphere packings, lattices and groups. Second edition.
Springer-Verlag, New York, 1993.}

\subsection{Poles of zeta functions\label{d4}}

Throughout this subsection, $V$ is a nonsingular projective variety over a
field $k$. We regard algebraic cycles on $V$ as elements of the $\mathbb{Q}%
$-vector spaces $\mathcal{Z}{}^{r}(V)$.

Algebraic cycles $D$ and $D^{\prime}$ are said to be numerically equivalent if
$D\cdot E=D^{\prime}\cdot E$ for all algebraic cycles $E$ on $V$ of
complementary dimension, and they are $\ell$-homologically equivalent if they
have the same class in $H^{2r}(V_{k^{\mathrm{al}}},\mathbb{Q}{}_{\ell}(r))$.
In his Woods Hole talk, Tate asked whether the following statement is always true:

\begin{quote}
$E^{r}(V)$: Numerical equivalence coincides with $\ell$-homological
equivalence for algebraic cycles on $V$ of codimension $r$.
\end{quote}

\noindent This is now generally regarded as a folklore conjecture (it is also
a consequence of Grothendieck's standard conjectures). Note that, like the
Tate conjecture, $E^{r}(V)$ is an existence statement for algebraic cycles:
for an algebraic cycle $D$, it says that there exists an algebraic cycle $E$
of complementary dimension such that $D\cdot E\neq0$ if there exists a
cohomological cycle with this property.

Let $\mathcal{A}{}^{r}$ denote the image of $\mathcal{Z}{}^{r}(V)$ in
$H^{2r}(V,\mathbb{Q}{}_{\ell}(r))$, and let $\mathcal{N}{}^{r}$ denote the
subspace of classes numerically equivalent to zero. Thus, $\mathcal{A}{}^{r}$
(resp. $\mathcal{A}{}^{r}/\mathcal{N}{}^{r}$) is the $\mathbb{Q}{}$-space of
algebraic classes of codimension $r$ modulo homological equivalence (resp.
modulo numerical equivalence). In particular, $\mathcal{A}{}^{r}/\mathcal{N}%
{}^{r}$ is independent of $\ell$.

Now assume that $k$ is finitely generated. We need to consider also the
following statement:

\begin{quote}
$S^{r}(V)$: The map $H^{2r}(V_{k^{\mathrm{al}}},\mathbb{Q}{}_{\ell
}(r))^{G(k^{\mathrm{al}}/k)}\rightarrow H^{2r}(V_{k^{\mathrm{al}}}%
,\mathbb{Q}{}_{\ell}(r))_{G(k^{\mathrm{al}}/k)}$ induced by the identity map
is bijective.
\end{quote}

\noindent When $k$ is finite, this means that $1$ occurs semisimply (if at
all) as an eigenvalue of the Frobenius map acting on $H^{2r}(V_{k^{\mathrm{al}%
}},\mathbb{Q}{}_{\ell}(r))$.

An elementary argument suffices to prove that the following three statements
are equivalent (for a fixed variety $V$, integer $r$, and prime $\ell$):

\begin{enumerate}
\item $T^{r}+E^{r};$

\item $T^{r}+T^{\dim V-r}+S^{r};$

\item $\dim_{\mathbb{Q}{}}(\mathcal{A}{}^{r}/\mathcal{N}{}^{r})=\dim
_{\mathbb{Q}_{l}}H^{2r}(V_{k^{\mathrm{al}}},\mathbb{Q}{}_{\ell}%
(r))^{G(k^{\mathrm{al}}/k)}.$
\end{enumerate}

\noindent When $k$ is finite, each statement is equivalent to:

\begin{enumerate}
\item[(d)] the order of the pole of $\zeta(V,s)$ at $s=r$ is equal to
$\dim_{\mathbb{Q}{}}(\mathcal{A}{}^{r}/\mathcal{N}^{r})$.
\end{enumerate}

\noindent See Tate 1979, 2.9. Note that (d) is independent of $\ell$.

Tate (1964a, Conjecture 2) conjectured the following general version of (d):

\begin{quote}
$P^{r}(V)$: Let $V$ be a nonsingular projective variety over a finitely
generated field $k$. Let $d$ be the transcendence degree of $k$ over the prime
field, augmented by $1$ if the prime field is $\mathbb{Q}{}$. Then the $2r$th
component $\zeta^{2r}(V,s)$ of the zeta function of $V$ has a pole of order
$\dim_{\mathbb{Q}{}}(\mathcal{A}{}^{r}(V))$ at the point $s=d+r$.
\end{quote}

\noindent This is also known as the Tate conjecture%
\index{Tate conjecture}%
. For a discussion of the known cases of $P^{r}(V)$, see Tate 1964a, 1994a.

When $k$ is a global function field, statements (a), (b), (c) are independent
of $\ell$, and are equivalent to the statement that $\zeta_{S}^{2r}(V,s)$ has
a pole of order $\dim_{\mathbb{Q}{}}(\mathcal{A}{}^{r}(V))$ at the point
$s=d+r$; here $\zeta_{S}^{2r}(V,s)$ omits the factors at a suitably large
finite set $S$ of primes. This follows from Lafforgue's proof of the global
Langlands correspondence and other results in Langlands program by an argument
that will, in principle, also work over number fields. --- see Lyons
2009.\footnote{Lyons, Christopher. A rank inequality for the Tate conjecture
over global function fields. Expo. Math. 27 (2009), no. 2, 93--108.}

In the presence of $E^{r}$, Conjecture $T^{r}(V)$ is equvalent to $P^{r}(V)$
if and only if the order of the pole of $\zeta^{2r}(V,s)$ at $s=r$ is
$\dim_{\mathbb{Q}{}_{l}}H^{2r}(V_{k^{\mathrm{al}}},\mathbb{Q}{}_{\ell
}(r))^{G(k^{\mathrm{al}}/k)}$. This is known for some Shimura varieties.

\subsubsection{The Sato-Tate conjecture}

Let $A$ be an elliptic curve over $\mathbb{Q}{}$. For a prime $p$ of good
reduction, the number $N_{p}$ of points on $A$ mod $p$ can be written%
\begin{align*}
N_{p}  &  =p+1-a_{p}\\
a_{p}  &  =2\sqrt{p}\cos\theta_{p},\quad0\leq\theta_{p}\leq\pi.
\end{align*}
When $A$ has complex multiplication over $\mathbb{C}{}$, it is easily proved
that the $\theta_{p}$ are uniformly distributed in the interval $0\leq
\theta\leq\pi$ as $p\rightarrow\infty$. In the opposite case, Mikio Sato found
computationally that the $\theta_{p}$ appeared to have a density distribution
$\frac{2}{\pi}\sin^{2}\theta$.

Tate proved that, for a power of an elliptic curve, the $\mathbb{Q}{}$-algebra
of algebraic cycles is generated modulo homological equivalence by divisor
classes. Using this, he computed that, for an elliptic curve $A$ over
$\mathbb{Q}{}$ without complex multiplication,%
\[
\rank(\mathcal{A}{}{}^{i}(A^{m}))=\binom{m}{i}^{2}-\binom{m}{i-1}\binom
{m}{i+1}\text{,}%
\]
from which he deduced that Sato's distribution is the only symmetric density
distribution for which the zeta functions of the powers of $A$ have their
zeros and poles in agreement with the Conjecture $P^{r}(V)$.

The conjecture that, for an elliptic curve over $\mathbb{Q}$ without complex
multiplication, the $\theta_{p}$ are distributed as $\frac{2}{\pi}\sin
^{2}\theta$ is known as the Sato-Tate conjecture.%
\index{Sato-Tate conjecture}
It has been proved only recently, as the fruit of a long collaboration
(Richard Taylor, Michael Harris, Laurent Clozel, Nicholas Shepherd-Barron,
Thomas Barnet-Lamb, David Geraghty). As did Tate, they approach the conjecture
through the analytic properties of the zeta functions of the powers of
$A$.\footnote{The proof was completed in:
\par
T.Barnet-Lamb, D.Geraghty, M.Harris and R.Taylor, A family of Calabi-Yau
varieties and potential automorphy II. P.R.I.M.S. 47 (2011), 29-98.
\par
For expository accounts, see:
\par
Carayol, Henri La conjecture de Sato-Tate (d'apr\`{e}s Clozel, Harris,
Shepherd-Barron, Taylor). S\'{e}minaire Bourbaki. Vol. 2006/2007.
Ast\'{e}risque No. 317 (2008), Exp. No. 977, ix, 345--391.
\par
Clozel, L. The Sato-Tate conjecture. Current developments in mathematics,
2006, 1--34, Int. Press, Somerville, MA, 2008.}

Needless to say, the Sato-Tate conjecture has been generalized to motives.
Langlands has pointed out that his functoriality conjecture contains a very
general form of the Sato-Tate conjecture.\footnote{Langlands, Robert P.
Reflexions on receiving the Shaw Prize. On certain L-functions, 297--308, Clay
Math. Proc., 13, Amer. Math. Soc., Providence, RI, 2011.}

\subsection{Relation to the Hodge conjecture\label{d5}}

For a variety $V$ over $\mathbb{C}{}$, there is a well-defined cycle map%
\[
c^{r}\colon\mathcal{Z}{}^{r}(V)\rightarrow H^{2r}(V,\mathbb{Q})
\]
(cohomology with respect to the complex topology). Hodge proved that there is
a decomposition%
\[
H^{2r}(V,\mathbb{Q})_{\mathbb{C}{}}=\bigoplus\nolimits_{p+q=2r}H^{p,q}%
,\quad\overline{H^{p,q}}=H^{q,p}.
\]
In Hodge 1950,\footnote{Hodge, W. V. D., The topological invariants of algebraic varieties.  Proceedings of the International Congress of Mathematicians, Cambridge, Mass., 1950, vol. 1,  pp. 182--192. Amer. Math. Soc., Providence, R. I., 1952.} he observed that the image of $c^{r}$ is contained in%
\[
H^{2r}(V,\mathbb{Q})\cap V^{r,r}%
\]
and asked whether this $\mathbb{Q}$-module is exactly the image of $c^{r}$.
This has become known as the Hodge conjecture.\footnote{Hodge actually asked
the question with $\mathbb{Z}{}$-coefficients.}

In his original article (Tate 1964), Tate wrote:

\begin{quote}
I can see no direct logical connection between [the Tate conjecture] and
Hodge's conjecture that a rational cohomology class of type $(p,p)$ is
algebraic\ldots. However, the two conjectures have an air of compatibility.
\end{quote}

\noindent Pohlmann (1968)\footnote{Pohlmann, Henry. Algebraic cycles on
abelian varieties of complex multiplication type. Ann. of Math. (2) 88 1968
161--180.} proved that the Hodge and Tate conjectures are equivalent for CM
abelian varieties, Piatetski-Shapiro (1971)\footnote{Piatetski-Shapiro, I. I.,
Interrelations between the Tate and Hodge hypotheses for abelian varieties.
(Russian) Mat. Sb. (N.S.) 85(127) (1971), 610--620.} proved that the Tate
conjecture for abelian varieties in characteristic zero implies the Hodge
conjecture for abelian varieties, and Milne (1999)\footnote{Milne, J. S.,
Lefschetz motives and the Tate conjecture. Compositio Math. 117 (1999), no. 1,
45--76.} proved that the Hodge conjecture for CM abelian varieties implies the
Tate conjecture for abelian varieties over finite fields.

The relation between the two conjectures has been greatly clarified by the
work of Deligne. He defines the notion of an absolute Hodge class on a
(complete smooth) variety over a field of characteristic zero, and conjectures
that every Hodge class on a variety over $\mathbb{C}$ is absolutely Hodge. The
Tate conjecture for a variety implies that all absolute Hodge classes on the
variety are algebraic. Therefore, in the presence of Deligne's conjecture, the
Tate conjecture implies the Hodge conjecture. As Deligne has proved his
conjecture for abelian varieties, this gives another proof of
Piatetski-Shapiro's theorem.\medskip

The twin conjectures of Hodge and Tate have a status in algebraic and
arithmetic geometry similar to that of the Riemann hypothesis in analytic
number theory. A proof of either one for any significantly large class of
varieties would be a major breakthrough. On the other hand, whether or not the
Hodge conjecture is true, it is known that Hodge classes behave in many ways
as if they were algebraic.\footnote{Deligne, P., Hodge cycles on abelian
varieties (notes by J.S. Milne). In: Hodge Cycles, Motives, and Shimura
Varieties, Lecture Notes in Math. 900, Springer-Verlag, 1982, pp. 9--100.
Cattani, Eduardo; Deligne, Pierre; Kaplan, Aroldo. On the locus of Hodge
classes. J. Amer. Math. Soc. 8 (1995), no. 2, 483--506.} There is some
fragmentary evidence that the same is true for Tate classes in nonzero
characteristic.\footnote{Milne, J. S. Rational Tate classes. Mosc. Math. J. 9
(2009), no. 1, 111--141.}

\section{Lubin-Tate theory and Barsotti-Tate group schemes\label{ee}}

\subsection{Formal group laws and applications\label{ee1}}

Let $R$ be a commutative ring. By a formal group law over $R$, we shall always
mean a one-parameter commutative formal group law, i.e., a formal power series
$F\in R[[X,Y]$] such that

\begin{itemize}
\item $F(X,Y)=X+Y+$terms of higher degree,

\item $F(F(X,Y),Z)=F(X,F(Y,Z))$

\item $F(X,Y)=F(Y,X)$.
\end{itemize}

\noindent These conditions imply that there exists a unique $i_{F}(X)\in
X\cdot R[[X]]$ such that $F(X,i_{F}(X))=0$. A homomorphism $F\rightarrow G$ of
formal group laws is a formal power series $f\in XR[[X]]$ such that
$f(F(X,Y))=G(f(X),f(Y))$.

The formal group laws form a $\mathbb{Z}$-linear category. Let $c(f)$ be the
first-degree coefficient of an endomorphism $f$ of $F$. If $R$ is an integral
domain of characteristic zero, then $f\mapsto c(f)$ is an injective
homomorphism of rings $\End_{R}(F)\rightarrow R$. See Lubin
1964.\footnote{Lubin, Jonathan, One-parameter formal Lie groups over p-adic
integer rings. Ann. of Math. (2) 80 1964 464--484.}

Let $F$ be a formal group law over a field $k$ of characteristic $p\neq0$. A
nonzero endomorphism $f$ of $F$ has the form
\[
f=aX^{p^{h}}+\text{terms of higher degree,}\quad a\neq0,
\]
where $h$ is a nonnegative integer, called the height of $f$. The height of
the multiplication-by-$p$ map is called the height of $F$.

\subsubsection{Lubin-Tate formal group laws and local class field theory}

Let $K$ be a nonarchimedean local field, i.e., a finite extension of
$\mathbb{Q}{}_{p}$ or $\mathbb{F}{}_{p}((t))$. Local class field theory
provides us with a homomorphism (the local reciprocity map)%
\[
\rec_{K}\colon K^{\times}\rightarrow\Gal(K^{\mathrm{ab}}/K)
\]
such that, for every finite abelian extension $L$ of $K$ in $K^{\mathrm{ab}}$,
$\rec_{K}$ induces an isomorphism%
\[
(-,L/K)\colon K^{\times}/\Nm L^{\times}\rightarrow\Gal(L/K);
\]
moreover, every open subgroup $K^{\times}$ of finite index arises as the norm
group of a (unique) finite abelian extension. This statement shows that the
finite abelian extensions of $K$ are classified by the open subgroups of
$K^{\times}$ of finite index, but leaves open the following problem:

\begin{quote}
Let $L/K$ be the abelian extension corresponding to an open subgroup $H$ of
$K^{\times}$ of finite index; construct generators for $L$ and describe how
$K^{\times}/H$ acts on them.
\end{quote}

\noindent Lubin and Tate (1965a) found an elegantly simple solution to this problem.

The choice of a prime element $\pi$ determines a decomposition $K^{\times
}=\mathcal{O}{}_{K}^{\times}\times\langle\pi\rangle$, and hence (by local
class field theory) a decomposition $K^{\mathrm{ab}}=K_{\pi}\cdot
K^{\mathrm{un}}$. Here $K_{\pi}$ is a totally ramified extension of $K$ with
the property that $\pi$ is a norm from every finite subextension. Since
$K^{\mathrm{un}}$ is well understood, the problem then is to find generators
for the subfields of $K_{\pi}$ and to describe the isomorphism%
\[
(-,K_{\pi}/K)\colon\mathcal{O}{}_{K}^{\times}\rightarrow\Gal(K_{\pi}/K)
\]
given by reciprocity map. Let $\mathcal{O}{}=\mathcal{O}_{K}$, let
$\mathfrak{p}{}=(\pi)$ be the maximal ideal in $\mathcal{O}$, and let $q=(\mathcal{O}\colon\mathfrak{p})$.

Let $f\in\mathcal{O}{}[[T]]$ be a formal power series such that%
\[
\left\{
\begin{array}
[c]{lll}%
f(T) & = & \pi T+\text{ terms of higher degree}\\
f(T) & \equiv & T^{q}\text{ modulo }\pi\mathcal{O}{}[[T]],
\end{array}
\right.
\]
for example, $f=\pi T+T^{q}$ is such a power series. An elementary argument
shows that, for each $a\in\mathcal{O}{}$, there is a unique formal power
series $[a]_{f}\in\mathcal{O}{}_{K}[[T]]$ such that%
\[
\left\{
\begin{array}
[c]{lll}%
\lbrack a]_{f}(T) & = & aT+\text{ terms of higher degree}\\
\lbrack a]_{f}\circ f & = & f\circ\lbrack a]_{f}\text{.}%
\end{array}
\right.
\]
Let
\[
X_{m}=\{x\in K^{\mathrm{al}}\mid\left\vert x\right\vert <1\text{, }%
\overbrace{(f\circ\cdots\circ f)}^{m}(x)=0\}.
\]

\noindent Then Lubin and Tate (1965a) prove:

\begin{enumerate}
\item the field $K[X_{m}]$ is the totally ramified abelian extension of $K$
with norm group $U_{m}\times\langle\pi\rangle$ where $U_{m}=1+\mathfrak{p}%
{}_{K}^{m}$;

\item the map%
\begin{align*}
\mathcal{O}{}^{\times}/U_{m}  &  \rightarrow\Gal(K[X_{m}]/K)\\
u  &  \mapsto\left(  x\mapsto\lbrack u^{-1}]_{f}(x)\right)
\end{align*}
is an isomorphism.
\end{enumerate}

\noindent For example, if $K=\mathbb{Q}{}_{p}$, $\pi=p$, and $f(T)=(1+T)^{p}%
-1$, then%
\[
X_{m}=\{\zeta-1\mid\zeta^{p^{m}}=1\}\simeq\mu_{p^{m}}(K^{\mathrm{al}})
\]
and $[u]_{f}(\zeta-1)=\zeta^{u^{-1}}-1$.

Lubin and Tate (1965a) show that, for each $f$ as above, there is a unique
formal group law $F_{f}$ admitting $f$ as an endomorphism. Then $X_{m}$ can be
realized as a group of \textquotedblleft torsion points\textquotedblright\ on
$F_{f}$, which endows it with the structure of an $\mathcal{O}{}$-module for
which it is isomorphic to $\mathcal{O}{}/\mathfrak{p}{}^{m}$. From this the
statements follow in a straightforward way.

The proof of the above results does not use local class field theory. Using
the Hasse-Arf theorem, one can show that $K_{\pi}\cdot K^{\mathrm{un}%
}=K^{\mathrm{ab}}$, and deduce local class field theory. Alternatively, using
local class field theory, one can show that $K_{\pi}\cdot K^{\mathrm{un}%
}=K^{\mathrm{ab}}$, and deduce the Hasse-Arf theorem. In either case, one
finds that the isomorphism in (b) is the local reciprocity map.

The $F_{f}$ are called \emph{Lubin-Tate formal group laws}%
\index{Lubin-Tate formal group law}%
, and the above theory is called \emph{Lubin-Tate theory}%
\index{Lubin-Tate theory}%
.

\subsubsection{Deformations of formal group laws (Lubin-Tate spaces)}

Let $F$ be a formal group law of height $h$ over a perfect field $k$ of
characteristic $p\neq0$. We consider local artinian $k$-algebras $A$ with
residue field $k$. A deformation of $F$ over such an $A$ is a formal group law
$F_{A}$ over $A$ such that $F_{A}\equiv F$ mod $\mathfrak{m}{}_{A}$. An
isomorphism of deformations is an isomorphism of $\varphi\colon F_{A}%
\rightarrow G_{A}$ of formal group laws such that $\varphi(T)\equiv T$ mod
$\mathfrak{m}{}_{A}$.

Let $W$ denote the ring of Witt vectors with residue field $k$. Lubin and Tate
(1966d) prove that there exists a formal group $\mathcal{F}(t_{1}%
,\ldots,t_{h-1})(X,Y){}$ over $W[[t_{1},\ldots,t_{h-1}]]$ with the following properties:

\begin{itemize}
\item $\mathcal{F}{}(0,\ldots,0)\equiv F$ mod $\mathfrak{m}{}_{W}$;

\item for any deformation $F_{A}$ of $F$, there is a unique homomorphism
$W[[t_{1},\ldots,t_{h-1}]]\rightarrow A$ sending $\mathcal{F}{}$ to $F_{A}$.
\end{itemize}

\noindent The results of Lubin and Tate are actually stronger, but this seems
to be the form in which they are most commonly used.

In particular, the above result identifies the space of deformations of $F$
with the formal scheme $\mathrm{Spf}(W[[t_{1},\ldots,t_{h-1}]).$ These spaces
are now called Lubin-Tate deformation spaces%
\index{Lubin-Tate deformation space}%
. Drinfeld showed that, by adding (Drinfeld) level structures, it is possible
to construct towers of deformation spaces, called Lubin-Tate towers%
\index{Lubin-Tate tower}%
.\footnote{Drinfeld, V. G. Elliptic modules. (Russian) Mat. Sb. (N.S.) 94(136)
(1974), 594--627, 656.} These play an important role in the study of the
moduli varieties of abelian varieties with PEL-structure and in the Langlands
program. For example, it is known that both the Jacquet-Langlands
correspondence and the local Langlands correspondence for $\GL_{n}$ can be
realized in the \'{e}tale cohomology of a Lubin-Tate tower (or, more
precisely, in the \'{e}tale cohomology of the Berkovich space that is attached to the
rigid analytic space which is the generic fibre of the Lubin-Tate
tower).\footnote{Carayol, H. Nonabelian Lubin-Tate theory. Automorphic forms,
Shimura varieties, and $L$-functions, Vol. II (Ann Arbor, MI, 1988), 15--39,
Perspect. Math., 11, Academic Press, Boston, MA, 1990.$\quad$Boyer, P.
Mauvaise r\'{e}duction des vari\'{e}t\'{e}s de Drinfeld et correspondance de
Langlands locale. Invent. Math. 138 (1999), no. 3, 573--629.$\quad$Harris,
Michael; Taylor, Richard The geometry and cohomology of some simple Shimura
varieties. With an appendix by Vladimir G. Berkovich. Annals of Mathematics
Studies, 151. Princeton University Press, Princeton, NJ, 2001.}

\subsection{Finite flat group schemes\label{ee2}}

A group scheme $G$ over a scheme $S$ is finite and flat if $G=\Spec(A)$ with
$A$ locally free of finite rank as a sheaf of $\mathcal{O}{}_{S}$-modules.
When $A$ has constant rank $r$, $G$ is said to have order $r$. A finite flat
group scheme of prime order is necessarily commutative.

In his course on Formal Groups at Harvard in the fall of 1966, Tate discussed
the following classification problem:

\begin{quote}
let $R$ be a local noetherian ring with residue field of characteristic
$p\neq0$; assume that $R$ contains the $(p-1)$st roots of $1$, i.e., that
$R^{\times}$ contains a cyclic subgroup of order $p-1$; determine the finite
flat group schemes of order $p$ over $R$.
\end{quote}

\noindent When $R$ is complete, Tate found that such group schemes correspond
to pairs of elements $(a,c)$ of $R$ such that $ac\in p\cdot R^{\times}$; two
pairs $(a,c)$ and $(a_{1},c_{1})$ correspond to isomorphic groups if and only
if $a_{1}=u^{p-1}a$ and $c_{1}=u^{1-p}c$ for some $u\in R^{\times}$.

These results were extended and completed in Tate and Oort 1970a. Let
\[
\Lambda_{p}=\mathbb{Z}{}[\zeta,\frac{1}{p(p-1)}]\cap\mathbb{Z}{}_{p}%
\]
where $\zeta$ is a primitive $(p-1)$st root of $1$ in $\mathbb{Z}{}_{p}$. Tate
and Oort define a sequence $w_{1}=1,w_{2},\ldots,w_{p}$ of elements of
$\Lambda_{p}$ in which $w_{1},\ldots,w_{p-1}$ are units and $w_{p}=pw_{p-1}$.
Then, given an invertible $\mathcal{O}{}_{S}$-module $L$ and sections $a$ of
$L^{\otimes(p-1)}$ and $b$ of $L^{\otimes(1-p)}$ such that $a\otimes b=w_{p}$,
they show that there is a group scheme $G_{a,b}^{L}$ such that, for every
$S$-scheme $T$,%
\[
G_{a,b}^{L}(T)=\{x\in\Gamma(T,L\otimes_{\mathcal{O}{}_{S}}\mathcal{O}{}%
_{T})\mid x^{\otimes p}=a\otimes x\},
\]
and the multiplication on $G_{a,b}^{L}(T)$ is given by%
\[
x_{1}\cdot x_{2}=x_{1}+x_{2}+\frac{b}{w_{p-1}}\otimes D_{p}(x_{1}%
\otimes1,1\otimes x_{2}),
\]
where%
\[
D_{p}(X_{1},X_{2})=\frac{w_{p-1}}{1-p}\sum_{i=1}^{p-1}\frac{X_{1}^{i}}{w_{i}%
}\frac{X_{2}^{p-i}}{w_{p-i}}\in\Lambda_{p}[X_{1},X_{2}].
\]
Every finite flat group scheme of order $p$ over $S$ is of the form
$G_{a,b}^{L}$ for some triple $(L,a,b)$, and $G_{a,b}^{L}$ is isomorphic to
$G_{a_{1},b_{1}}^{L_{1}}$ if and only if there exists an isomorphism from $L$
to $L_{1}$ carrying $a$ to $a_{1}$ and $b$ to $b_{1}$. The Cartier dual of
$G_{a,b}^{L}$ is $G_{b,a}^{L^{-1}}$. The proofs of these statements make
ingenious use of the action of $\mathbb{F}{}_{p}^{\times}$ on $\mathcal{O}%
{}_{G}$.

Tate and Oort apply their result to give a classification of finite flat group
schemes of order $p$ over the ring of integers in a number field in terms of
id\`{e}le class characters. In particular, they show that the only such group
schemes over $\mathbb{Z}{}$ are the constant group scheme $\mathbb{Z}%
{}/p\mathbb{Z}{}$ and its dual $\mu_{p}$.

In the years since Tate and Oort wrote their article, the classification of
finite flat commutative group schemes over various bases has been intensively
studied, and some of the results were used in the proof of the modularity
conjecture for elliptic curves (hence of Fermat's last theorem). Tate (1997a)
has given a beautiful exposition of the basic theory of finite flat group
schemes, including the results of Raynaud 1974\footnote{Raynaud, Michel,
Sch\'{e}mas en groupes de type (p,\ldots,p). Bull. Soc. Math. France 102
(1974), 241--280.} extending the above theory to group schemes of type
$(p,\ldots,p)$.

\subsection{Barsotti-Tate groups $p$-divisible
groups)\label{ee3}}

Let $A$ be an abelian variety over a field $k$. In his study of abelian
varieties and their zeta functions, Weil used the $\ell$-primary component
$A(\ell)$ of the group $A(k^{\mathrm{sep}})$ for $\ell$ a prime distinct from
$\mathrm{char}(k)$. This is an $\ell$-divisible group isomorphic to
$(\mathbb{Q}{}_{\ell}/\mathbb{Z}{}_{\ell})^{2\dim A}$ equipped with an action
of $\Gal(k^{\mathrm{sep}}/k)$. For $p=\mathrm{char}(k)$, it is natural to replace
$A(\ell)$ with the direct system%
\[
A(p)\colon\quad A_{1}\hookrightarrow A_{2}\hookrightarrow\cdots\hookrightarrow
A_{\nu}\hookrightarrow A_{\nu+1}\hookrightarrow\cdots
\]
where $A_{\nu}$ is the finite group \emph{scheme} $\Ker(p^{\nu}\colon
A\rightarrow A)$.

In the mid sixties, Serre and Tate\footnote{Usually this is credited to Tate
alone, but Tate writes: \textquotedblleft We were both contemplating them. I
think it was probably Serre who first saw clearly the simple general
definition and its relation to formal groups of finite
height.\textquotedblright\ The dual of a $p$-divisible group is often called the
Serre dual.
\par
{}} defined a $p$\emph{-divisible group} of height $h$ over a ring $R$ to be a
direct system $G=(G_{\nu},i_{\nu})_{\nu\in\mathbb{N}{}}$ where, for each
$\nu\geq0$, $G_{\nu}$ is a finite group scheme over $R$ of order $p^{\nu h}$
and the sequence%
\[
0\longrightarrow G_{\nu}\overset{i_{\nu}}{\longrightarrow}G_{\nu+1}%
\overset{p^{\nu}}{\longrightarrow}G_{\nu+1}%
\]
is exact. An abelian scheme $A$ over $R$ defines a $p$-divisible group $A(p)$
over $R$ of height $2\dim(A)$.

The dual of a $p$-divisible group $G=(G_{\nu},i_{\nu})$ is the system
$G^{\prime}=(G_{\nu}^{\prime},i_{\nu}^{\prime})$ where $G_{\nu}^{\prime}$ is
the Cartier dual of $G_{\nu}$ and $i_{\nu}^{\prime}$ is the Cartier dual of
the map $G_{\nu+1}\overset{p}{\longrightarrow}G_{\nu}$. It is again a
$p$-divisible group.

Tate developed the basic theory $p$-divisible group in his article for the
proceedings of the 1966 Driebergen conference (Tate 1967c) and in a series of
ten lectures at the Coll\`{e}ge de France in 1965-1966.\footnote{See Serre,
OEuvres, II pp.321--324.} He showed that $p$-divisible groups generalize formal
Lie groups in the following sense: let $R$ be a complete noetherian local ring
with residue field $k$ of characteristic $p>0$; an $n$-dimensional commutative formal Lie
group $\Gamma$ over $R$ can be defined to be a family $f(Y,Z)=(f_{i}%
(Y,Z))_{1\leq i\leq n}$ of $n$ power series in $2n$ variables satisfying the
conditions in the first paragraph of this section; if such a group $\Gamma$ is
divisible (i.e., $p\colon\Gamma\rightarrow\Gamma$ is an isogeny), then one can
define the kernel $G_{\nu}$ of $p^{\nu}\colon\Gamma\rightarrow\Gamma$ as a
group scheme over $R$; Tate shows that $\Gamma(p)=(G_{\nu})_{\nu\geq1}$ is a
$p$-divisible group $\Gamma(p)$, and that the functor $\Gamma\rightsquigarrow
\Gamma(p)$ is an equivalence from the category of divisible commutative formal
Lie groups over $R$ to the category of connected $p$-divisible groups over $R$.

The main theorem of Tate 1967c states the following:

\begin{quote}
Let $R$ be an integrally closed, noetherian, integral domain whose field of
fractions $K$ is of characteristic zero, and let $G$ and $H$ be $p$-divisible
groups over $R$. Then every homomorphism $G_{K}\rightarrow H_{K}$ of the
generic fibres extends uniquely to a homomorphism $G\rightarrow H$.
\end{quote}

\noindent In other words, the functor $G\rightsquigarrow G_{K}$ is fully
faithful. This was extended to rings $R$ of characteristic $p\neq0$ by de Jong
in 1998.\footnote{de Jong, A. J., Homomorphisms of Barsotti-Tate groups and
crystals in positive characteristic. Invent. Math. 134 (1998), no. 2,
301--333. Erratum: ibid. 138 (1999), no. 1, 225.}

Since their introduction, $p$-divisible groups have become an essential tool
in the study of abelian schemes. We have already seen in (\ref{b4}) one
application of $p$-divisible groups to the problem of lifting abelian
varieties. Another application was to the proof of the Mordell conjecture
(Faltings 1983\footnote{Faltings, G. Endlichkeitss\"{a}tze f\"{u}r abelsche
Variet\"{a}ten \"{u}ber Zahlk\"{o}rpern. Invent. Math. 73 (1983), no. 3,
349--366. Erratum: Ibid. (1984), no. 2, 381.}).

\ In his talk at the 1970 International Congress, Grothendieck renamed
$p$-divisible groups \textquotedblleft Barsotti-Tate groups\textquotedblright%
\index{Barsotti-Tate group}%
. Today, both terms are used.

\subsection{Hodge-Tate decompositions\label{ee4}}

Now let $R$ be a complete discrete valuation ring of unequal characteristics,
and let $K$ be its field of fractions. Let $K^{\mathrm{al}}$ be an algebraic
closure of $K$, and let $C$ be the completion of $K^{\mathrm{al}}$. As Serre
noted:\footnote{\OE \ uvres, II, p.322.}

\begin{quote}
One of the most surprising results of Tate's theory is the fact that \emph{the
properties of }$p$\emph{-divisible groups are intimately related to the
structure of} $C$ \emph{as a Galois module over} $\Gal(K^{\mathrm{al}}/K)$.
\end{quote}

\noindent Let $\mathcal{G}{}{}=\Gal(K^{\mathrm{al}}/K)$, and let $V$ be a
$C$-vector space on which $\mathcal{G}{}$ acts semi-linearly. The Tate twist%
\index{Tate twist}
$V(i)$, $i\in\mathbb{Z}{}$, is $V$ with $\mathcal{G}{}$ acting by
\[
(\sigma,v)\mapsto\chi(\sigma)^{i}\cdot\sigma v,\quad\chi\text{ the
}p\text{-adic cyclotomic character.}%
\]
\noindent Tate proved that $H^{0}(\mathcal{G}{},C)=K$ and $H^{1}(\mathcal{G}%
{},C)\approx K$, and that $H^{q}(\mathcal{G}{},C(i))=0$ for $q=0,1$ and
$i\neq0$.\footnote{Sen and Ax simplified and generalized Tate's proof that
$C^{G}=K$, and the result is now known as the Ax-Sen-Tate theorem%
\index{Ax-Sen-Tate theorem}%
.} Using these statements, he proved that, for a $p$-divisible group $G$ over
$R$, there is a canonical isomorphism%
\begin{equation}
\Hom(TG,C)\simeq t_{G^{\prime}}(C)\oplus t_{G}^{\vee}(C)(-1). \label{e21}%
\end{equation}
where $TG=\varprojlim_{\nu}G_{\nu}(K^{\mathrm{al}})$ is the Tate module%
\index{Tate module}
of $G$, $t_{G}$ is the tangent space to $G$ at zero, and $G^{\prime}$ is the
dual $p$-divisible group. In particular $TG$ determines the dimension of $G$,
a fact that is used in the proof of the main theorem in the last subsection.

When $G$ is the $p$-divisible group of an abelian scheme $A$ over $R$,
(\ref{e21}) can be written as:%
\[
H_{\mathrm{et}}^{1}(A_{C},\mathbb{Q}{}_{p})\otimes C\simeq H^{1}(A_{C}%
,\Omega_{A_{C}/C}^{0})\oplus H^{0}(A_{C},\Omega_{A_{C}/C}^{1})(-1).
\]
This result led Tate to make the following (Hodge-Tate%
\index{Hodge-Tate conjecture}%
) conjecture:\footnote{Tate 1967c, p.180; see also Serre's summary of Tate's
lectures, \OE uvres, II, p.324.}

\begin{quote}
For every nonsingular projective variety $X$ over $K$, there exists a
canonical (Hodge-Tate%
\index{Hodge-Tate decomposition}%
) decomposition%
\begin{equation}
H_{\mathrm{et}}^{n}(X_{C},\mathbb{Q}{}_{p})\otimes_{\mathbb{Q}{}_{p}}%
C\simeq\bigoplus\nolimits_{p+q=n}H^{p,q}(X_{C})(-p) \label{e22}%
\end{equation}
where $H^{p,q}(X_{C})=H^{q}(X,\Omega_{X/K}^{p})\otimes_{K}C$. This
decomposition is compatible with the action of $\Gal(K^{\mathrm{al}}/K$).
\end{quote}

Tate's conjecture launched a new subject in mathematics, called $p$-adic Hodge
theory. The isomorphism (\ref{e22}) can be regarded as a statement about the
\'{e}tale cohomology of $X_{C}$ regarded as a module over $\Gal(K^{\mathrm{al}%
}/K)$. About 1980, Fontaine stated a series of successively stronger
conjectures, beginning with the Hodge-Tate conjecture, that describe the
structure of these Galois modules.\footnote{See Fontaine, Jean-Marc. Sur
certains types de repr\'{e}sentations $p$-adiques du groupe de Galois d'un
corps local; construction d'un anneau de Barsotti-Tate. Ann. of Math. (2) 115
(1982), no. 3, 529--577, and many other articles.} Most of Fontaine's
conjectures have now been proved. The Hodge-Tate conjecture itself was proved
by Faltings in 1988.\footnote{Faltings, Gerd. p-adic Hodge theory. J. Amer.
Math. Soc. 1 (1988), no. 1, 255--299.}

\section{Elliptic curves\label{f}}

Although elliptic curves are just abelian varieties of dimension one, their
study is quite different. Throughout his career, Tate has returned to the
study of elliptic curves.

\subsection{Ranks of elliptic curves over global fields\label{f1}}

Mordell proved that, for an elliptic curve $E$ over $\mathbb{Q}{}$, the group
$E(\mathbb{Q}{})$ is finite generated. At one time, it was widely conjectured
that the rank of $E(\mathbb{Q}{})$ is bounded, but, as Cassels 1966 pointed
out, this is implausible.\footnote{\textquotedblleft It has been widely
conjectured that there is an upper bound for the rank depending only on the
groundfield. This seems to me implausible because the theory makes it clear
that an abelian variety can only have high rank if it is defined by equations
with very large coefficients.\textquotedblright\ p.257 of Cassels, J. W. S.,
Diophantine equations with special reference to elliptic curves. J. London
Math. Soc. 41 1966 193--291.} Tate and Shafarevich (1967d)\footnote{For a long
time I was puzzled as to how this article came to be written, because I was
not aware that Shafarevich had been allowed to travel to the West, but Tate
writes: \textquotedblleft sometime during the year 1965--66, which I spent in
Paris, Shafarevich appeared. There must have been a brief period when the
Soviets relaxed their no-travel policy\ldots. Shafarevich was in Paris for a
month or so, and the paper grew out of some discussion we had. We both liked
the idea of our having a joint paper, and I was happy to have it in
Russian.\textquotedblright} made it even less plausible by proving that, for
elliptic curves over the global field $\mathbb{F}{}_{p}(t)$, the ranks are
unbounded. Their examples are quadratic twists of a supersingular elliptic
curve with coefficients in $\mathbb{F}{}_{p}$; in particular, they are
isotrivial (i.e., have $j\in\mathbb{F}{}_{p}$). More recently, it has been
shown that the ranks are unbounded even among the nonisotrivial elliptic
curves over $\mathbb{F}{}_{p}(t)$.\footnote{Ulmer, Douglas Elliptic curves
with large rank over function fields. Ann. of Math. (2) 155 (2002), no. 1,
295--315.} Meanwhile, the largest known rank for an elliptic curve over
$\mathbb{Q}{}$ is $28$.\footnote{Elkies, see
\url{http://web.math.hr/~duje/tors/tors.html}.}

\subsection{Torsion points on elliptic curves over
$\mathbb{Q}$\label{f2}}

Beppo Levi constructed elliptic curves $E$ over $\mathbb{Q}{}$ having each of
the groups%
\[
{}%
\begin{array}
[c]{rcl}%
\mathbb{Z}{}/n\mathbb{Z} & \quad & n=1,2,\ldots,10,12,\\
\mathbb{Z}{}/2\mathbb{Z}{}\times\mathbb{Z}{}/n\mathbb{Z}{} &  & n=2,4,6,8,
\end{array}
\]
as the torsion subgroup of $E(\mathbb{Q}{})$, and he conjectured that this
exhausted the list of possible such groups.\footnote{Beppo Levi,
Sull'equazione indeterminata del 3${}^{\circ}$ ordine, Rom. 4. Math. Kongr. 2,
173-177 1909. (Talk at the 1908 ICM.)}

Consistent with this, Mazur and Tate (1973c) show that there is no elliptic
curve over $\mathbb{Q}{}$ with a rational point of order $13$, or,
equivalently, that the curve $X_{1}(13)$ that classifies the elliptic curves
with a chosen point of order $13$ has no rational points (except for its
cusps). Ogg found a rational point of order $19$ on the Jacobian $J$ of
$X_{1}(13)$, and Mazur and Tate show that $J$ has exactly $19$ rational
points. They then deduce that $X_{1}(13)$ has no noncuspidal rational point by
examining how it sits in its Jacobian.

The interest of their article is more in its methods than in the result
itself.\footnote{About the same time, J. Blass found a more elementary proof
of the same result.} The ring $\mathbb{Z[}\sqrt[3]{1}]$ acts on $J$, and Mazur
and Tate perform a 19-descent by studying the flat cohomology of the exact
sequence of group schemes%
\[
0\rightarrow F\rightarrow J\overset{\pi}{\longrightarrow}J\rightarrow0
\]
on $\Spec\mathbb{Z}{}\smallsetminus(13)$, where $\pi$ is one of the factors of
$19$ in $\mathbb{Z[}\sqrt[3]{1}]$. In a major work, Mazur developed these
methods further, and completely proved Levi's conjecture.\footnote{Mazur, B.
Modular curves and the Eisenstein ideal. Inst. Hautes \'{E}tudes Sci. Publ.
Math. No. 47 (1977), 33--186.}

The similar problem for an arbitrary number field $K$ is probably beyond
reach, but, following work of Kamienny, Merel (1996\footnote{Merel, Lo\"{\i}c.
Bornes pour la torsion des courbes elliptiques sur les corps de nombres.
Invent. Math. 124 (1996), no. 1-3, 437--449.}) proved that, for a fixed number
field $K$, the order of the torsion subgroup of $E(K)$ for $E$ an elliptic
curve over $K$ is bounded by a constant that depends only the degree of $K$
over $\mathbb{Q}{}$.

\subsection{Explicit formulas and algorithms\label{f3}}

The usual Weierstrass form of the equation of an elliptic curve is valid only
in characteristics $\neq2,3.$ About 1965 Tate wrote out the complete form,
valid in all characteristics, and even over $\mathbb{Z}{}$. For an elliptic
curve over a nonarchimedean local field with perfect residue field, he wrote
out an explicit algorithm (known as Tate's algorithm%
\index{Tate's algorithm}%
) for computing the minimal model of the curve and determining the Kodaira
type of the special fibre. Ogg's formula then gives the conductor of the
curve. The handwritten manuscript containing these formulas was invaluable to
people working in the field. A copy, which had been sent to Cassels, was
included, essentially verbatim, in the proceedings of the Antwerp conference
(Tate 1975b).\footnote{Tate writes: \textquotedblleft Early in that summer
[1965], Weil had told me of the idea that all elliptic curves over
$\mathbb{Q}{}$ are modular [and that the conductor of the elliptic curve
equals the conductor of the corresponding modular form]. That motivated
Swinnerton-Dyer to make a big computer search for elliptic curves over
$\mathbb{Q}{}$ with not too big discriminant, in order to test Weil's idea.
But of course it was necessary to be able to compute the conductor to do that
test. That was my main motivation.\textquotedblright}

\subsection{Analogues at $p$ of the conjecture of Birch
and Swinnerton-Dyer\label{f4}}

For an elliptic curve $E$ over $\mathbb{Q}{}$, the conjecture of Birch and
Swinnerton-Dyer states that%
\[
L(s,E)\quad\sim\quad\Omega\prod_{p\text{ bad}}c_{p}\frac{\left\vert
\TS(E)\right\vert \cdot R}{\left\vert E(\mathbb{Q}{})_{\mathrm{tors}%
}\right\vert ^{2}}(s-1)^{r}\quad\text{as}\quad s\rightarrow1
\]
where $r$ is the rank of $E(\mathbb{Q}{})$, $\TS(E)$ is the Tate-Shafarevich
group, $R$ is the discriminant of the height pairing on $E(\mathbb{Q}{})$,
$\Omega$ is the real period of $E$, and the $c_{p}=(E(\mathbb{Q}{}_{p})\colon
E^{0}(\mathbb{Q}{}_{p})$). When $E$ has good ordinary or multiplicative
reduction at $p$, there is a $p$-adic zeta function $L_{p}(s,E)$, and Mazur,
Tate, and Teitelbaum (1986) investigated whether the behaviour of $L_{p}(s,E)$
near $s=1$ is similarly related to the arithmetic invariants of $E$%
.\footnote{The authors assume the $E$ is modular --- at the time, it was not
known that all elliptic curves over $\mathbb{Q}{}$ are modular.} They found it
is, but with one major surprise: there is an \textquotedblleft
exceptional\textquotedblright\ case in which $L_{p}(s,E)$ is related to an
extended version of $E(\mathbb{Q}{})$ rather than $E(\mathbb{Q}{})$ itself.
Supported by numerical evidence, they conjectured:

\begin{quote}
BSD$(p)$. When $E$ has good ordinary or nonsplit multiplicative reduction at a
prime $p$, the function $L_{p}(s,E)$ has a zero at $s=1$ of order at least
$r=\rank E(\mathbb{Q}{})$, and $L_{p}^{(r)}(1,E)$ is equal to a certain
expression involving $|\TS(E)|$ and a $p$-adic regulator $R_{p}(E)$. When $E$
has split multiplicative reduction, it is necessary to replace $r$ with $r+1$.
\end{quote}

\noindent\noindent The $L$-function $L_{p}(s,E)$ is the $p$-adic Mellin
transform of a $p$-adic measure obtained from modular symbols. The $p$-adic
regulator is the discriminant of the canonical $p$-adic height pairing
(augmented in the exceptional case). Much more is known about BSD$(p)$ than
the original conjecture of Birch and Swinnerton-Dyer. For example,
Kato\footnote{Kato, Kazuya. $p$-adic Hodge theory and values of zeta functions
of modular forms. Cohomologies $p$-adiques et applications arithm\'{e}tiques.
III. Ast\'{e}risque No. 295 (2004), ix, 117--290.} has proved the following statement:

\begin{quote}
The function $L_{p}(s,E)$ has a zero at $s=1$ of order at least $r$ (at least
$r+1$ in the exceptional case). When the order of the zero equals its
conjectured value, then the $p$-primary component of $\TS(E)$ is finite and
$R_{p}(E)\neq0$.
\end{quote}

\noindent In the exceptional case, $E_{\mathbb{Q}{}_{p}}$ is a Tate elliptic
curve and $L_{p}(1,E)=0$. On comparing their conjecture in this case with the
original conjecture of Birch and Swinnerton-Dyer, the authors were led to the
conjecture%
\[
L_{p}^{(1)}(1,E)=\frac{\log_{p}(q)}{\ord_{p}(q)}\frac{L(1,E)}{\Omega}%
\]
where $q$ is the period of the Tate curve $E_{\mathbb{Q}{}_{p}}$ and $\Omega$
is the real period of $E$. This became known as the Mazur-Tate-Teitelbaum
conjecture.%
\index{Mazur-Tate-Teitelbaum conjecture}
It was proved by Greenberg and Stevens in 1993\footnote{Greenberg, Ralph;
Stevens, Glenn, $p$-adic $L$-functions and $p$-adic periods of modular forms.
Invent. Math. 111 (1993), no. 2, 407--447.} for $p\neq2,3$, and by several
authors in general.

Mazur and Tate (1987) state \textquotedblleft refined\textquotedblright%
\ conjectures that avoid any mention of $p$-adic $L$-functions and, in
particular, avoid the problem of constructing these functions. Let $E$ be an
elliptic curve over $\mathbb{Q}{}$. For a fixed integer $M>0$, they use
modular symbols to construct an element $\theta$ in the group algebra
$\mathbb{Q}{}[(\mathbb{Z}{}/M\mathbb{Z}{})^{\times}/\{\pm1\}]$. Let $R$ be a
subring of $\mathbb{Q}{}$ containing the coefficients of $\theta$ and such
that the order the torsion subgroup of $E(\mathbb{Q}{})$ is invertible in $R$.
The analogue of an $L$-function having a zero of order $r$ at $s=1$ is that
$\theta$ lie in the $r$th power of the augmentation ideal $I$ of the group
algebra $R[(\mathbb{Z}{}/M\mathbb{Z}{})^{\times}/\{\pm1\}]$. Assume that $M$
is not divisible by $p^{2}$ for any prime $p$ at which $E$ has split
multiplicative reduction. Then Mazur and Tate conjectured:

\begin{quote}
Let $r$ be the rank of $E(\mathbb{Q}{})$, and let $r^{\prime}$ be the number
of primes dividing $M$ at which $E$ has split multiplicative reduction. Then
$\theta\in I^{r+r^{\prime}}$, and there is a formula (involving the order of
$\TS (E)$) for the image of $\theta$ in $I^{r+r^{\prime}}/I^{r+r^{\prime}+1}$.
\end{quote}

\noindent Again, the authors provide numerical evidence for their conjecture.
Tate's student, Ki-Seng Tan, restated the Mazur-Tate conjecture for an
elliptic curve over a global field, and he proved that part of the new
conjecture is implied by the conjecture of Birch and
Swinnerton-Dyer.\footnote{Tan, Ki-Seng, Refined theorems of the Birch and
Swinnerton-Dyer type. Ann. Inst. Fourier (Grenoble) 45 (1995), no. 2,
317--374.}

In the first article discussed above, Mazur, Tate, Teitelbaum gave explicit
formulas relating the canonical $p$-adic height pairings to a $p$-adic sigma
function, and used the sigma function to study the height pairings. Mazur and
Tate (1991), present a detailed construction of the $p$-adic sigma function
for an elliptic curve $E$ with good ordinary reduction over a $p$-adic field
$K$, and they prove the properties used in the earlier article. In contrast to
the classical sigma function, which is defined on the universal covering of
$E$, the $p$-adic sigma function is defined on the formal group of $E$.

Finally, the article Mazur, Stein, and Tate 2006 studies the problem of
efficiently computing of $p$-adic heights for an elliptic curve $E$ over a
global field $K$. This amounts to efficiently computing the sigma function,
which in turn amounts to efficiently computing the $p$-adic modular form
$E_{2}$.

\subsection{Jacobians of curves of genus one\label{f5}}

For a curve $C$ of genus one over a field $k$, the Jacobian variety $J$ of $C$
is an elliptic curve over $k$. The problem is to compute a Weierstrass
equation for $J$ from an equation for $C$.

Weil (1954)\footnote{Weil, Andr\'{e}. Remarques sur un m\'{e}moire d'Hermite.
Arch. Math. (Basel) 5, (1954). 197--202.} showed that, when $C$ is defined by
an equation $Y^{2}=f(X)$, $\deg f=4$, then the Weierstrass equation of $J$ can
be computed using the invariant theory of the quartic of $f$, which goes back
to Hermite.

An et al. (2001)\footnote{An, Sang Yook; Kim, Seog Young; Marshall, David C.;
Marshall, Susan H.; McCallum, William G.; Perlis, Alexander R. Jacobians of
genus one curves. J. Number Theory 90 (2001), no. 2, 304--315.} showed how
formulas from classical invariant theory give Weierstrass equations for $J$
and for the map $C\rightarrow J$ when $\mathrm{char}(k)\neq2,3$ and $C$ is a
double cover of $\mathbb{P}{}^{1}$, a plane cubic, or a space quartic.

Tate and Rodrigues-Villegas found the Weierstrass equations over fields of
characteristic $2$ and $3$, where the classical invariant theory no longer
applies. Together with Artin they extended their result to an arbitrary base
scheme $S$ (Artin, Rodriguez-Villegas, Tate 2005). Specifically, let $C$ be
the family of curves over a scheme $S$ defined, as a subscheme of
$\mathbb{P}{}_{S}^{2}$, by a cubic $f\in\Gamma(S,\mathcal{O}{}_{S})[X,Y,Z]$,
and assume that the ten coefficients of $f$ have no common zero. Then there is
a Weierstrass equation
\[
g\colon\quad Y^{2}Z+a_{1}XYZ+a_{3}YZ^{2}=X^{3}+a_{2}X^{2}Z+a_{4}XZ^{2}%
+a_{6}Z^{3},\quad a_{i}\in\Gamma(S,\mathcal{O}{}_{S}),
\]
whose coefficients are given explicitly in terms of the coefficients of $f$,
such that the functor $\Pic_{C/S}^{0}$ is represented by the smooth locus of
the subscheme of $\mathbb{P}{}_{S}^{2}$ defined by $g$. A key ingredient of
the proof is a characterization, over sufficiently good base schemes, of the
group algebraic spaces that can be described by a Weierstrass equation.

\subsection{Expositions\label{f6}}

In 1961, Tate gave a series of lectures at Haverford College titled
\textquotedblleft Rational Points on Cubic Curves\textquotedblright\ intended
for bright undergraduates in mathematics. Notes were taken of the lectures,
and these were distributed in mimeographed form. The book, Silverman and Tate
1992, is a revision, and expansion, of the notes.

In the spring of 1960, the fall of 1967, 1975,\ldots, Tate gave courses on the
arithmetic of elliptic curves, whose informal notes have influenced later expositions.

\section{The $K$-theory of number fields\label{g}}

\subsection{$K$-groups and symbols\label{g1}}

Grothendieck defined $K_{0}(X)$ for $X$ a scheme in order to be able to state
his generalization of the Riemann-Roch theorem. The topologists soon adapted
Grothendieck's definition to topological spaces, and extended it to obtain
groups $K_{n}$ for all $n\in\mathbb{\mathbb{N}{}}$.

For a commutative ring $R$, $K_{0}(R)$ is just the Grothendieck group of the
category of finitely generated projective $R$-modules. In 1962, Bass and
Schanuel\footnote{Bass, H.; Schanuel, S. The homotopy theory of projective
modules. Bull. Amer. Math. Soc. 68 1962 425--428.} defined $K_{1}(R)$, and in
1967, Milnor\footnote{During a course at Princeton University; published as:
Milnor, John. Introduction to algebraic $K$-theory. Annals of Mathematics
Studies, No. 72. Princeton University Press, Princeton, N.J., 1971.} defined
$K_{2}(R)$. In the early 1970s, several authors suggested definitions for the
higher $K$-groups, which largely coincided when this could be checked.
Quillen's definition\footnote{Quillen, Daniel. Higher algebraic $K$-theory. I.
Algebraic $K$-theory, I: Higher $K$-theories (Proc. Conf., Battelle Memorial
Inst., Seattle, Wash., 1972), pp. 85--147. Lecture Notes in Math., Vol. 341,
Springer, Berlin 1973.} was the most flexible, and it is his that has been adopted.

The Steinberg group $\mathrm{ST}(R)$ of a ring $R$ is the group with
generators%
\[
x_{ij}(r),\quad i,j=1,2,3,\ldots,\quad i\neq j\text{,}\quad r\in R
\]
and relations%
\begin{align*}
x_{ij}(r)x_{ij}(s)  &  =x_{ij}(r+s)\\
\lbrack x_{ij}(r),x_{kl}(s)]  &  =\left\{
\begin{array}
[c]{ll}%
1 & \text{if }i\neq l\text{ and }j\neq k\\
x_{il}(rs) & \text{if }i\neq l\text{ and }j=k.
\end{array}
\right.
\end{align*}
The elementary matrices $E_{ij}(r)=I+re_{ij}$ in $\GL(R)$ satisfy these
relations, and so there is a homomorphism $x_{ij}(r)\mapsto E_{ij}%
(r)\colon\mathrm{ST}(R)\rightarrow\GL(R)$. The groups $K_{1}(R)$ and
$K_{2}(R)$ can be defined by the exact sequence%
\[
1\rightarrow K_{2}(R)\rightarrow\mathrm{ST}(R)\rightarrow\GL(R)\rightarrow
K_{1}(R)\rightarrow1.
\]

Let $F$ be a field. A \emph{symbol} on $F$ with values in a commutative group
$C$ is defined to be a bimultiplicative map%
\[
(\,\,\,,\,\,)\colon F^{\times}\times F^{\times}\rightarrow C
\]
such that $(a,1-a)=1$ whenever $a\neq0,1.$ Matsumoto
(1969)\footnote{Matsumoto, Hideya, Sur les sous-groupes arithm\'{e}tiques des
groupes semi-simples d\'{e}ploy\'{e}s. Ann. Sci. \'{E}cole Norm. Sup. (4) 2
1969 1--62. (1968 thesis ENS).} showed that the natural map%
\[
\{\,\,\,,\,\,\}\colon F^{\times}\times F^{\times}\rightarrow K_{2}F
\]
is a universal symbol, i.e., that $K_{2}(F)$ is the free abelian group on the pairs
$\{a,b\}$, $a,b\in F^{\times}$, subject only to the relations%
\begin{align*}
\{aa^{\prime},b\}  &  =\{a,b\}\{a^{\prime},b\} &  &  \text{all }a,a^{\prime},b\in
F^{\times}\\
\{a,bb^{\prime}\}  &  =\{a,b\}\{a,b^{\prime}\} &  &  \text{all }a,b,b^{\prime}\in
F^{\times}\\
\{a,1-a\}  &  =1 &  &  \text{all }a\neq0,1\text{ in }F^{\times}.
\end{align*}

\subsubsection{Examples of symbols}

\begin{enumerate}
\item The tame (Hilbert) symbol. Let $v$ be a discrete valuation of $F$, with
residue field $\kappa(v)$. Then%
\[
(a,b)_{v}=(-1)^{v(a)v(b)}\frac{a^{v(b)}}{b^{v(a)}}\quad\text{mod }%
\mathfrak{m}{}_{v}%
\]
is a symbol on $F$ with values in $\kappa(v)^{\times}.$

\item The Galois symbol (Tate 1970b, \S 1). For $m$ not divisible by
$\mathrm{char}(F)$, $H^{1}(F,\mu_{m})\simeq F^{\times}/F^{\times m}$, and the
cup-product pairing%
\[
H^{1}(F,\mu_{m})\times H^{1}(F,\mu_{m})\rightarrow H^{2}(F,\mu_{m}\otimes
\mu_{m})
\]
gives a symbol on $F$ with values in $H^{2}(F,\mu_{m}\otimes\mu_{m})$. When
$F$ contains the $m$th roots of $1$,%
\[
H^{2}(F,\mu_{m}\otimes\mu_{m})\simeq H^{2}(F,\mu_{m})\otimes\mu_{m}%
\simeq\Br(F)_{m}\otimes\mu_{m}%
\]
and the symbol was known classically.

\item The differential symbol (Tate ibid.). For $p=\mathrm{char}(F)$,
\[
f,g\mapsto\frac{df}{f}\wedge\frac{dg}{g}\colon F^{\times}\times F^{\times
}\rightarrow\Omega_{F/\mathbb{F}{}_{p}}^{2}%
\]
is a symbol.

\item On $\mathbb{C}{}$ there are no continuous symbols, but on $\mathbb{R}{}$
there is the symbol%
\[
(a,b)_{\infty}=\left\{
\begin{array}
[c]{rl}%
1 & \text{if }a>0\text{ or }b>0\\
-1 & \text{otherwise.}%
\end{array}
\right.
\]

\end{enumerate}

\subsection{The group $K_{2}F$ for
$F$ a global field\label{g2}}

Tate recognized that the study of the $K_{2}$ of a global field is related to
classical objects in number fields, and sheds new light on them. He largely
initiated the study of the $K$-groups of global fields and their rings of integers.

For a field $F$, $K_{0}F\simeq\mathbb{Z}{}$ is without particular interest. On
the other hand, $K_{1}F\simeq F^{\times}$. For a number field, there is an
exact sequence%
\[
0\rightarrow U_{F}\rightarrow F^{\times}\overset{(\ord_{v})}{\longrightarrow
}\bigoplus\nolimits_{v}\mathbb{Z}{}\rightarrow C_{F}\rightarrow0
\]
where $v$ runs over the finite primes of $F$. Dirichlet proved that
$U_{F}\approx\mu(F)\times\mathbb{Z}{}^{r_{1}+r_{2}-1}$, where $r_{1}$ and
$r_{2}$ are the numbers of real and complex primes, and Dedekind proved that
the class group $C_{F}$ is finite. Thus understanding $K_{1}F$ involves
understanding the two basic objects in classical algebraic number theory.

Let $F$ be a global field. For a noncomplex prime $v$ of $F$, let $\mu_{v}%
=\mu(F_{v})$ and let $m_{v}=\left\vert \mu_{v}\right\vert $. For a finite
prime $v$ of $F$, $\Br(F_{v})\simeq\mathbb{Q}{}/\mathbb{Z}{}$, and so the
Galois symbol with $m=m_{v}$ gives a homomorphism $\lambda_{v}\colon
K_{2}F_{v}\rightarrow\mu_{v}$. Similarly, $(\,\,\,,\,\,)_{\infty}$ gives a
homomorphism $\lambda_{v}\colon K_{2}F_{v}\rightarrow\mu_{v}$ when $v$ is
real. The $\lambda_{v}$ can be combined with the obvious maps $K_{2}%
F\rightarrow K_{2}F_{v}$ to give the homomorphism $\lambda_{F}$ in the
sequence%
\[
0\rightarrow\Ker(\lambda_{F})\rightarrow K_{2}F\overset{\lambda_{F}%
}{\longrightarrow}\bigoplus\nolimits_{v}\mu_{v}\rightarrow\mu_{F}%
\rightarrow0\text{;}%
\]
the direct sum is over the noncomplex primes of $F$, and the map from it sends
$(x_{v})_{v}$ to $\prod\nolimits_{v}\frac{m_{v}}{m_{F}}x_{v}$ where
$m_{F}=\left\vert \mu(F)\right\vert $. A product formula,%
\[
\prod(a,b)_{v}^{\frac{m_{v}}{m_{F}}}=1
\]
shows that the sequence is a complex, and Moore (1969)\footnote{Moore, Calvin
C. Group extensions of p-adic and adelic linear groups. Inst. Hautes
\'{E}tudes Sci. Publ. Math. No. 35 1968 157--222.} showed that the cokernel of
$\lambda_{F}$ is $\mu_{F}$. Thus, to compute $K_{2}F$, it remains to identify
$\Ker(\lambda_{F})$.

For $\mathbb{Q}{}$, Tate proved that $\Ker(\lambda_{F})$ is trivial, and then
observed that most of his argument was already contained in Gauss's first
proof of the quadratic reciprocity law.

For a global field $F$, Bass and Tate proved that $\Ker(\lambda_{F})$ is
finitely generated, and that it is finite of order prime to the characteristic
in the function field case. Garland (1971) proved that it is also finite in
the number field case.

In the function field case, Tate proved that%
\begin{equation}
\left\vert \Ker(\lambda_{F})\right\vert =(q-1)\cdot(q^{2}-1)\cdot\zeta
_{F}(-1)\text{.} \label{e12}%
\end{equation}
For a number field, the Birch-Tate conjecture%
\index{Birch-Tate conjecture}
says that%
\begin{equation}
\left\vert \Ker(\lambda_{F})\right\vert =\pm w_{2}(F)\cdot\zeta_{F}(-1)
\label{e11}%
\end{equation}
where $w_{2}(F)$ is the largest integer $m$ such $\Gal(F^{\mathrm{al}}/F)$
acts trivially on $\mu_{m}(F^{\mathrm{al}})\otimes\mu_{m}(F^{\mathrm{al}})$
(Birch 1971;\footnote{Birch, B. J. $K_{2}$ of global fields. 1969 Number
Theory Institute (Proc. Sympos. Pure Math., Vol. XX, State Univ. New York,
Stony Brook, N.Y., 1969), pp. 87--95. Amer. Math. Soc., Providence, R.I.,
1971.} Tate 1970b). The odd part of this conjecture was proved by
Wiles.\footnote{Wiles, A. The Iwasawa conjecture for totally real fields. Ann.
of Math. (2) 131 (1990), no. 3, 493--540.}

When Quillen defined the higher $K$-groups, he proved that%
\[
K_{2}(\mathcal{O}{}_{F})=\Ker\left(  K_{2}(F)\rightarrow\bigoplus
\nolimits_{v\text{ finite}}\mu_{v}\right)
\]
and so there is an exact sequence%
\[
0\rightarrow\Ker(\lambda_{F})\rightarrow K_{2}(\mathcal{O}{}_{F}%
)\rightarrow\bigoplus\nolimits_{v\text{ real}}\mu_{v}.
\]
Thus the computation of $\Ker(\lambda_{F})$ is closely related to that of
$K_{2}(\mathcal{O}{}_{F})$. Lichtenbaum\footnote{Lichtenbaum, Stephen. On the
values of zeta and $L$-functions. I. Ann. of Math. (2) 96 (1972), 338--360.}
generalized the Birch-Tate conjecture to the following statement: for all
totally real number fields $F$\footnote{The first test of the conjecture was
for $F=\mathbb{Q}{}$ and $i=1$. Since $\zeta_{\mathbb{Q}{}}(-1)=-1/12$ and
$K_{2}(\mathbb{Z}{})=\mathbb{Z}{}/2\mathbb{Z}{}$, the conjecture predicts that
$|K_{3}(\mathbb{Z}{})|$ has $24$ elements, but Lee and Szczarba showed that it
has $48$ elements. When a seminar speaker at Harvard mentioned this, and
scornfully concluded that the conjecture was false, Tate responded from the
audience \textquotedblleft Only for $2$\textquotedblright. In fact,
Lichtenbaum's conjecture is believed to be correct up to a power of $2$.}%
\[
\frac{\left\vert K_{4i-2}(\mathcal{O}{}_{F})\right\vert }{\left\vert
K_{4i-1}(\mathcal{O}{}_{F})\right\vert }=\left\vert \zeta_{F}(1-2i)\right\vert
,\quad\text{all }i\geq1.
\]

\subsubsection{The Galois symbol}

Tate proved (\ref{e12}) by using Galois symbols. For a global field $F$, he
proved that the map
\begin{equation}
K_{2}F\rightarrow H^{2}(F,\mu_{m}^{\otimes2}) \label{e13}%
\end{equation}
defined by the Galois symbol induces an isomorphism%
\begin{equation}
K_{2}F/(K_{2}F)^{m}\simeq H^{2}(F,\mu_{m}^{\otimes2}) \label{e23}%
\end{equation}
when $m$ is not divisible by $\mathrm{char}(F)$, and wrote \textquotedblleft I
don't know whether \ldots\ this holds for all fields\textquotedblright\ Tate
(1970b, p.208). Merkurjev and Suslin (1982)\footnote{Merkurjev, A. S.; Suslin,
A. A. $K$-cohomology of Severi-Brauer varieties and the norm residue
homomorphism. (Russian) Izv. Akad. Nauk SSSR Ser. Mat. 46 (1982), no. 5,
1011--1046, 1135--1136.} proved that it does hold for all fields.

Tate noted that the isomorphism (\ref{e23}) gives little information on
$\Ker(\lambda_{F})$ because $\bigcap\nolimits_{m}\left(  K_{2}F\right)  ^{m}$
is a subgroup of $\Ker(\lambda_{F})$ of index at most $2$, and $\Ker(\lambda
_{F})\subset(K_{2}F)^{m}$ for all $m$ not divisible by $8$. He then defined
more refined Galois symbols, which are faithful.

Fix a prime $\ell\neq\mathrm{char}(F)$, and let $\mathbb{Z}{}_{\ell
}(1)=\varprojlim_{n}\mu_{\ell^{n}}(F^{\mathrm{al}}))$. This is a free
$\mathbb{Z}{}_{\ell}$-module of rank $1$ with an action of $\Gal(F{}%
^{\mathrm{al}}/F)$, and we let $H^{r}(F,\mathbb{Z}{}_{\ell}(1)^{\otimes2})$
denote the Galois cohomology group defined using continuous cocycles (natural
topology on both $\Gal(F^{\mathrm{al}}/F)$ and $\mathbb{Z}{}_{\ell
}(1)^{\otimes2}$). Tate proves that the maps (\ref{e13}) with $m=\ell^{n}$
lift to a map
\[
K_{2}F\rightarrow H^{2}(F,\mathbb{Z}{}_{\ell}(1)^{\otimes2}),
\]
and that this map induces an isomorphism%
\[
K_{2}F(\ell)\rightarrow H^{2}(F,\mathbb{Z}{}_{\ell}(1)^{\otimes2}%
)_{\mathrm{tors}}\text{.}%
\]
As $K_{2}F$ is torsion, with no $\mathrm{char}(F)$-torsion, this gives a
purely cohomological description of $K_{2}F$.

\begin{nt}
The results of Tate in this subsection were announced in Tate 1970b, and
proved in Tate 1973b, 1976b, or in Tate's appendix to Bass and Tate 1973a.
\end{nt}

\subsection{The Milnor $K$-groups\label{g3}}

Milnor (1970)\footnote{Milnor, John. Algebraic $K$-theory and quadratic forms.
Invent. Math. 9 1969/1970 318--344.} defines the (Milnor) $K$-groups of a
field $F$ as follows: regard $F^{\times}$ as a $\mathbb{Z}{}$-module; then
$K_{\ast}^{M}F$ is the quotient of the tensor algebra of $F^{\times}$ by the
ideal generated by the elements%
\[
a\otimes(1-a),\quad a\neq0,1.
\]
This means that, for $n\geq2$, $K_{n}^{M}F$ is the quotient of $(K_{1}%
F)^{\otimes n}$ by the subgroup generated by the elements%
\[
a_{1}\otimes\cdots\otimes a_{n},\quad a_{i}+a_{i+1}=1\text{ for some }i.
\]
There is a canonical homomorphism $K_{\ast}F\rightarrow K_{\ast}^{M}F$ which
induces isomorphisms $K_{i}F\rightarrow K_{i}^{M}F$ for $i\leq2$. In the same
article, Milnor defined for each discrete valuation $v$ on $F$, a homomorphism%
\[
\partial_{v}\colon K_{\ast}F\rightarrow K_{\ast}\kappa(v)
\]
of degree $-1$, where $\kappa(v)$ is the residue field.

Milnor (ibid.) quotes a theorem of Tate: for a global field $F$,
\[
K_{n}^{M}F/2K_{n}^{M}F\simeq\bigoplus\nolimits_{v}K_{n}^{M}F_{v}/2K_{n}%
^{M}F_{v},\quad n\geq3,
\]
which implies that $K_{n}^{M}F/2K_{n}^{M}F\simeq(\mathbb{Z}{}/2\mathbb{Z}%
{})^{r_{1}}$ ($n\geq3)$ where $r_{1}$ is the number of real primes of $F$.
Bass and Tate (1973a) improve this statement by showing that%
\[
K_{n}^{M}F\simeq(\mathbb{Z}{}/2\mathbb{Z}{})^{r_{1}}\text{ for }n\geq3.
\]
The proof makes essential use of the \textquotedblleft transfer
maps\textquotedblright%
\[
\Tr\colon K_{\ast}E\rightarrow K_{\ast}F,
\]
defined whenever $E$ is a finite field extension of $F$. Since these had only
been defined for $n\leq2$, a major part of the article is taken up with
proving results on $K_{\ast}F$ for a general field, including the existence of
the transfer maps.

The theorem of Bass and Tate completes the determination of the Milnor
$K$-groups of a global field, except for $K_{1}$ and $K_{2}$.

\subsection{Other results on $K_{2}F$\label{g4}}

Let $F$ be a field containing a primitive $m$th root $\zeta$ of $1$ for some
$m>1$. Tate (1976) showed that, when $F$ is a global field, every element of
$K_{2}F$ killed by $m$ can be represented as $\{\zeta,a\}$ for some $a\in
F^{\times}$. Tate (1977b) examines whether this holds for other fields and
obtains a number of positive results.

For a finite extension of fields $E\supset F$, there is a transfer (or trace)
map $\Tr_{E/F}\colon K_{2}E\rightarrow K_{2}F$. As $K_{2}E$ is generated by
symbols $\{a,b\}$, in order to describe $\Tr_{E/F}$ it suffices to describe
its action on each symbol. This Rosset and Tate (1983c) do.

\section{The Stark conjectures\label{h}}

In a series of four papers, Stark examined the behaviour of Artin $L$-series
near $s=0$ (equivalently, $s=1$), and stated his now famous
conjectures.\footnote{Stark, H. M., Values of $L$-functions at $s=1$. I.
$L$-functions for quadratic forms. Advances in Math. 7 1971 301--343 (1971);
II. Artin $L$-functions with rational characters. ibid. 17 (1975), no. 1,
60--92; III. Totally real fields and Hilbert's twelfth problem. ibid. 22
(1976), no. 1, 64--84; IV. First derivatives at $s=0$. ibid. 35 (1980), no. 3,
197--235.} Tate gave a seminar at Harvard on Stark's conjectures in the spring
of 1978, after Stark had given some talks on the subject in the number theory
seminar in the fall of 1977. In 1980/81 Tate gave a course at the
Universit\'{e} de Paris-Sud (Orsay) in which he clarified and extended Stark's
work in important ways. The notes of Tate's course, when published in 1984,
included most of the results known at that date, and became the basic
reference for the Stark conjectures.

Let $\zeta_{k}(s)$ be the zeta function of a number field $k$. A celebrated
theorem of Dedekind shows that%
\begin{equation}
\zeta_{k}(s)\,\,\,\sim\,\,\,-\frac{R}{(e/h)}\,\,s^{r_{1}+r_{2}-1}\quad\text{as
}s\rightarrow0, \label{e18}%
\end{equation}
where $h$ is the class number of $k$, $R$ is its regulator, $e=\left\vert
\mu(k)\right\vert $, and $r_{1}+r_{2}-1$ is the rank of the group of units in
$k$.

Let $K$ be a finite Galois extension of $k$, with Galois group $G=G(K/k)$.
Stark's insight was that the decomposition of $\zeta_{K}(s)$ into a product of
Artin $L$-series indexed by the irreducible characters of $G$ should induce an
interesting decomposition of (\ref{e18}).

\subsubsection{Stark's main conjecture}

Let $\chi\colon G\rightarrow\mathbb{C}{}$ be the character of a
finite-dimensional complex representation $\rho\colon G\rightarrow\GL(V)$ of
$G$. For a finite set $S$ of primes of $k$ containing the infinite primes, let%
\[
L(s,\chi)=\prod\nolimits_{\mathfrak{p}{}\notin S}\frac{1}{\det(1-\rho
(\sigma_{\mathfrak{p}{}})N\mathfrak{\mathfrak{p}{}}^{-s}|V^{I_{\mathfrak{P}{}%
}})}%
\]
(Artin $L$-function relative to $S$; cf. \ref{a1}). Let $S_{K}$ be the set of
primes of $K$ lying over a prime in $S$, let $Y$ be the free $\mathbb{Z}{}%
$-module on $S_{K}$, and let $X$ be the submodule of $Y$ of elements $\sum
n_{w}w$ such that $\sum n_{w}=0$. Then $L(s,\chi)$ has a zero of multiplicity
$r(\chi)$ at $s=0$, where
\[
r(\chi)=\dim_{\mathbb{C}{}}\Hom_{G}(V^{\vee},X_{\mathbb{C}{}}).
\]
Let $U$ be the group of $S_{K}$-units in $K$. The unit theorem provides us
with an isomorphism%
\[
\lambda\colon U_{\mathbb{R}{}}\rightarrow X_{\mathbb{R}{}},\quad u\mapsto
\sum\nolimits_{w\in S_{K}}\log\left\vert u\right\vert _{w}w\text{.}%
\]
For each choice of an isomorphism of $G$-modules $f\colon X_{\mathbb{Q}{}%
}\rightarrow U_{\mathbb{Q}{}}$, Tate (1984, p.26) defines the Stark regulator,
$R(\chi,f)$, to be the determinant of the endomorphism of $\Hom_{G}(V^{\vee
},X_{\mathbb{C}{}})$ induced by $\lambda_{\mathbb{C}{}}\circ f_{\mathbb{C}{}}%
$. Then%
\[
L(s,\chi)\quad\sim\quad\frac{R(\chi,f)}{A(\chi,f)}\,\,s^{r(\chi)}\quad\text{
as }s\rightarrow0
\]
for a complex number $A(\chi,f)$. The main conjecture of Stark, as formulated
by Tate (1984, p.27), says that%
\[
A(\chi,f)^{\alpha}=A(\chi^{\alpha},f)\text{ for all automorphisms }%
\alpha\text{ of }\mathbb{C}{},
\]
where $\chi^{\alpha}=\alpha\circ\chi$. In particular, $A(\chi,f)$ is an
algebraic number, lying in the cyclotomic field $\mathbb{Q}{}(\chi)$. Tate
proves that the validity of the conjecture is independent of both $f$ and $S$,
and that it suffices to prove it for irreducible characters $\chi$ of
dimension $1$ (application of Brauer's theorem p.\pageref{1955b}).

\subsubsection{Characters with values in $\mathbb{Q}{}$}

When the character $\chi$ takes its values in $\mathbb{Q}{}$, Stark's
conjecture predicts that $A(\chi,f)\in\mathbb{Q}{}$ for all $f$. If, in
addition, $\chi$ is a $\mathbb{Z}{}$-linear combination of characters induced
from trivial characters, then the proof of the conjecture comes down to the
case of a trivial character, where it follows from (\ref{e18}). Some multiple
of $\chi$ has this form, and so this shows that some power of $A(\chi
,f)$ is in $\mathbb{Q}{}$ (Stark 1975). Tate proves (1984, Chapter II) that
$A(\chi,f)$ itself lies in $\mathbb{Q}{}$. His proof makes heavy use of the
cohomology of number fields, including the theorems in \ref{a3}.

\subsubsection{The case that $L(s,\chi)$ is nonzero at $s=0$}

When $r(\chi)=0$, the Stark regulator $R(\chi,f)=1$, and Stark's conjecture
becomes the statement:%
\[
L(0,\chi)^{\alpha}=L(0,\chi^{\alpha})\text{ for all automorphisms }%
\alpha\text{ of }\mathbb{C}{}.
\]
This special case of Stark's conjecture is also a special case of Deligne's
conjecture on the critical values of motives (Deligne 1979,
\S 6).\footnote{Deligne, P. Valeurs de fonctions L et p\'{e}riodes
d'int\'{e}grales. Proc. Sympos. Pure Math., XXXIII, Automorphic forms,
representations and L-functions (Proc. Sympos. Pure Math., Oregon State Univ.,
Corvallis, Ore., 1977), Part 2, pp. 313--346, Amer. Math. Soc., Providence,
R.I., 1979.} Using a refinement of Brauer's theorem (cf. p.\pageref{1955b}),
Tate writes $L(s,\chi)$ as a sum of partial zeta functions:%
\[
L(s,\chi)=\sum_{\sigma\in G(K/k)}\chi(\sigma)\cdot\zeta(s,\sigma),\quad
\zeta(s,\sigma)=\sum_{(\mathfrak{a}{},K/k)=\sigma}N\mathfrak{a}{}^{-s}%
\]
(Tate 1984, III 1). According to an important theorem of
Siegel,\footnote{Siegel, Carl Ludwig. \"{U}ber die Fourierschen Koeffizienten
von Modulformen. Nachr. Akad. Wiss. G\"{o}ttingen Math.-Phys. Kl. II 1970
15--56.} $\zeta(0,\sigma)\in\mathbb{Q}{}$, which proves Stark's conjecture in
this case.

\subsubsection{The case that $L(s,\chi)$ has a first order zero at $s=0$}

By contrast, when $r(\chi)=1$, the conjecture is still unknown, but it has
remarkable consequences. Let $\mathbb{C}{}[G]$ be the group algebra of $G$,
and let%
\[
e_{\chi}=\frac{\chi(1)}{\left\vert G\right\vert }\sum_{\sigma\in G}\chi
(\sigma^{-1})\cdot\sigma
\]
be the idempotent in $\mathbb{C}{}[G]$ that projects every representation of
$G$ onto its $\chi$-component. For an $a\in\mathbb{Q}{}(\chi)$ and a character
$\chi$ with $r(\chi)=1$, let%
\[
\pi(a,\chi)=\sum_{\alpha\in G(\mathbb{Q}{}(\chi)/\mathbb{Q}{})}a^{\alpha}\cdot
L^{\prime}(0,\chi^{\alpha})\cdot e_{\bar{\chi}^{\alpha}}\in\mathbb{C}{}[G].
\]
The character $\chi$ is realized on a $\mathbb{Q}{}[G]$-submodule $U_{W}$ of
$U_{\mathbb{Q}{}}$, and Stark's conjecture is true for $\chi$ if and only if
\[
\pi(a,\chi)X_{\mathbb{Q}{}}=\lambda(U_{W})\quad\quad\text{(inside
}X_{\mathbb{C}}\text{).}%
\]
(Tate 1984, III 2.1). More explicitly, let $\Psi$ be a set of irreducible
characters $\chi\neq1$ of $G$ such that $r(\chi)=1$, and assume that $\Psi$ is
stable under $\Aut(\mathbb{C}{})$. Let $(a_{\chi})_{\chi\in\Psi}$ be a family
of complex numbers such that $a_{\chi^{\alpha}}=(a_{\chi})^{\alpha}$ for all
$\alpha\in\Aut(\mathbb{C}{})$. If Stark's conjecture holds for the $\chi
\in\Psi$, then for each prime $v$ of $S$ and extension of $v$ to a prime $w$
of $K$, there exists an integer $m>0$ and an $S$-unit $\varepsilon$ of $K$
such that%
\begin{equation}
\lambda(\varepsilon)=m\sum_{\chi\in\Psi}a_{\chi}\cdot L^{\prime}(0,\chi)\cdot
e_{\bar{\chi}}\cdot w; \label{e19}%
\end{equation}
once $m$ has been fixed, $\varepsilon$ is unique up to a root of $1$ in $K$
(ibid. III \S 3). The units $\varepsilon$ arising (conjecturally) in this way
are called Stark units. They are analogous to the cyclotomic units in
cyclotomic fields.

\subsubsection{Finer conjectures when $K/k$ is abelian}

When $K/k$ is abelian, (\ref{e19}) can be made into a more precise form of
Stark's conjecture, which Tate denotes $\mathrm{St}(K/k,S)$ (Stark 1980; Tate
1984, IV 2). For a real prime $w$ of $K$ and certain hypotheses on $S$,
$\mathrm{St}(K/k,S)$ predicts the existence of a unit $\varepsilon(K,S,w)\in
U$ such that%
\[
\varepsilon(K,S,w)^{\sigma}=\exp(-2\zeta^{\prime}(0,\sigma)),\quad\text{all
}\sigma\in G.
\]
When we use $w$ to embed $K$ in $\mathbb{R}{}$, the $\varepsilon(K,S,w)$ lie
in the abelian closure of $k$ in $\mathbb{R}{}$. In the case that $k$ is
totally real, Tate (1984, 3.8) determines the subfield they generate; for
example, when $[K\colon\mathbb{Q}{}]=2$, they generate the abelian closure of
$k$ in $\mathbb{R}{}$. This has implications for Hilbert's 12th problem. To
paraphrase Tate (ibid. p.95):

\begin{quote}
If the conjecture $\mathrm{St}(K/k,S)$ is true in this situation, then the
formula%
\[
\varepsilon=\exp(-2\zeta^{\prime}(0,1))
\]
gives generators of abelian extensions of $k$ that are special values of
transcendental functions. Finding generators of class fields of this shape is
the vague form of Hilbert's 12th problem, and the Stark conjecture represents
an important contribution to this problem. However, it is a totally unexpected
contribution: Hilbert asked that we discover the functions that play, for an
arbitrary number field, the same role as the exponential function for
$\mathbb{Q}{}$ and the elliptic modular functions for a quadratic imaginary
field. In contrast, Stark's conjecture, by using $L$-functions directly,
bypasses the transcendental functions that Hilbert asked for. Perhaps a
knowledge of these last functions will be necessary for the proof of Stark's conjecture.
\end{quote}

\noindent Remarkably, $\mathrm{St}(K/k,S)$ is useful for the explicit
computation of class fields, and has even been incorporated into the computer
algebra system PARI/GP.

For an abelian extension $K/k$, Tate introduced another conjecture, combining
ideas of Brumer and Stark, and which he calls the Brumer-Stark conjecture. Let
$S$ be a set of primes of $k$ including a finite prime $\mathfrak{p}{}$ that
splits completely in $K$, and let $T=S\smallsetminus\{\mathfrak{p}{}\}$.
Assume that $T$ contains the infinite primes and the primes that ramify in
$K$. Let%
\[
\theta_{T}(0)=\sum_{\chi\text{ irreducible}}L(0,\chi)e_{\bar{\chi}}%
\in\mathbb{C}{}[G].
\]
Brumer conjectured that, for every ideal $\mathfrak{A}$ of $K$, $\mathfrak{A}%
{}^{e\theta_{T}(0)}$ is principal; the Brumer-Stark conjecture BS$(K/k,T)$
says that $\mathfrak{A}{}^{e\theta_{T}(0)}=(\alpha)$ for an $\alpha$
satisfying certain conditions on the absolute values $\left\vert
\alpha\right\vert _{w}$ ($w\in T$) and that $K(\alpha^{\frac{1}{e}})$ is an
abelian extension of $k$ (Tate 1984, 6.2).\footnote{Recall that $e=|\mu(k)|$.}
Tate proved this conjecture for $k=\mathbb{Q}{}$ (ibid. 6.7) and for quadratic
extensions $K/k$ (Tate 1981b).

\subsubsection{Function fields}

All of the conjectures make sense for a global field $k$ of characteristic
$p\neq0$. In this case, the Artin $L$-series are rational functions in
$q^{-s}$, where $q$ is the order of the field of constants, and Stark's main
conjecture follows easily from the known properties of these functions.
However, as Mazur pointed out, the Brumer-Stark conjecture is far from trivial
for function fields. Tate gave a seminar in Paris in early fall 1980 in which
he discussed the conjecture and some partial results he had obtained. Deligne
attended the seminar, and later gave a proof of the conjecture using his
one-motives. This proof is included in Chapter V of Tate 1984.

\subsubsection{$p$-adic analogues}

Tate's reformulation of Stark's conjecture helped inspire two $p$-adic
analogues of his main conjecture, one for $s=0$ (Gross) and one for $s=1$
(Serre) --- the absence of a functional equation for the $p$-adic $L$-series
makes these distinct conjectures. In a 1997 letter, Tate proposed a refinement
of Gross's conjecture. This letter was published, with additional comments, as
Tate 2004.\medskip

There is much numerical evidence for the Stark conjectures, found by Stark and
others. As Tate (1981a, p.977) notes: \textquotedblleft Taken all together,
the evidence for the conjectures seems to me to be
overwhelming\textquotedblright.

\section{Noncommutative ring theory\label{i}}

The Tate conjecture for divisors on a variety is related to the finiteness of
the Brauer group of the variety, which is defined to be the set of the
similarity classes of sheaves of (noncommutative) Azumaya algebras on the
variety. This connection led M.\thinspace Artin to an interest in
noncommutative rings, which soon broadened beyond Azumaya algebras. Tate wrote a number of articles on noncommutative rings
in collaboration with Artin and others.

\subsection{Regular algebras\label{i1}}

A ring $A$ is said to have \emph{finite global dimension} if there exists an
integer $d$ such that every $A$-module has a projective resolution of length
at most $d$. The smallest such $d$ is then called the \emph{global dimension}
of $A$. Serre\footnote{Serre, Jean-Pierre, Alg\`{e}bre locale.
Multiplicit\'{e}s. Seconde \'{e}dition, 1965. Lecture Notes in Mathematics, 11
Springer-Verlag, Berlin-New York 1965} showed that a \textit{commutative} ring
is noetherian of finite global dimension if and only if it is regular.

Let $k$ be a field. Artin and Schelter (1987)\footnote{Artin, Michael;
Schelter, William F.; Graded algebras of global dimension 3. Adv. in Math. 66
(1987), no. 2, 171--216.} defined a finitely generated $k$-algebra to be
\emph{regular} if it is of the form%
\begin{equation}
A=k\oplus A_{1}\oplus A_{2}\oplus\cdots, \label{e10}%
\end{equation}
and

\begin{enumerate}
\item $A$ has finite global dimension (defined in terms of \emph{graded }$A$-modules),

\item $A$ has polynomial growth (i.e., $\dim A_{n}$ is bounded by a polynomial
function in $n$), and

\item $A$ is Gorenstein (i.e., the $k$-vector space $\Ext_{A}^{i}(k,A)$ has
dimension $1$ when $i$ is the global dimension of $A$, and is zero otherwise).
\end{enumerate}

\noindent The only commutative graded $k$-algebras satisfying these conditions
and generated in degree $1$ are the polynomial rings. It is expected that the
regular algebras have many of the good properties of polynomial rings. For
example, Artin and Schelter conjecture that they are noetherian domains. The
dimension of a regular algebra is its global dimension.

In collaboration with Artin and others, Tate studied regular algebras,
especially the classification of those of low dimension.

From now on, we require regular $k$-algebras to be generated in degree $1$.
Such an algebra is a quotient of a tensor algebra by a homogeneous ideal.

A regular $k$-algebra of dimension one is a polynomial ring, and one of
dimension two is the quotient of the free associative algebra $k\langle
X,Y\rangle$ by a single quadratic relation, which can be taken to be $XY-cYX$
$(c\neq0)$ or $XY-YX-Y^{2}$. Thus, the first interesting dimension is three.
Artin and Schelter (ibid.) showed that a regular $k$-algebra of dimension
three either has three generators and three relations of degree two, or two
generators and two relations of degree three. Moreover, they showed that the
algebras fall into thirteen families. While the generic members of each family
are regular, they were unable to show that all the algebras in the families
are regular. Artin, Tate, Van den Bergh (1990a) overcame this problem, and
consequently gave a complete classification of these algebras. Having found an
explicit description of all the algebras, they were able to show that they are
all noetherian.

These two articles introduced new geometric techniques into noncommutative
ring theory. They showed that the regular algebras of dimension 3 correspond
to certain triples $(E,\mathcal{L}{},\sigma)$ where $E$ is a one-dimensional
scheme of arithmetic genus $1$ which is embedded either as a cubic divisor in
$\mathbb{P}{}^{2}$ or as a divisor of bidegree $(2,2)$ in $\mathbb{P}{}%
^{1}\times\mathbb{P}{}^{1}$, $\mathcal{L}{}=\mathcal{O}{}_{E}(1)$ is an
invertible sheaf on $E$, and $\sigma$ is an automorphism of $E$. The scheme
$E$ parametrizes the point modules for $A$, i.e., the graded cyclic right
$A$-modules, generated in degree zero, such that $\dim_{k}(M_{n})=1$ for all
$n\geq0$. The geometry of $(E,\sigma)$ is reflected in the structure of the
point modules, and Artin, Tate, van den Bergh (1991a) exploit this relation to
prove that the $3$-dimensional regular algebra corresponding to a triple
$(E,\mathcal{L},\sigma{})$ is finite over its centre if and only if the
automorphism $\sigma$ has finite order. They also show that noetherian regular
$k$-algebras of dimension $\leq4$ are domains.

\subsection{Quantum groups\label{i2}}

A bialgebra $A$ over a field $k$ is a $k$-module equipped with compatible
structures of an associative algebra with identity and of a coassociative
coalgebra with coidentity. A bialgebra is called a Hopf algebra if it admits
an antipodal map (linear map $S\colon A\rightarrow A$ such that certain
diagrams commute).

A bialgebra is said to be commutative if it is commutative as a $k$-algebra.
The commutative bialgebras (resp. Hopf algebras) over $k$ are exactly the
coordinate rings of affine monoid schemes (resp. affine group schemes) over
$k$.

Certain Hopf algebras (not necessarily commutative) are called quantum groups.
For example, there is a standard one-parameter family $\mathcal{O}{}%
(\GL_{n}(q))$, $q\in k^{\times}$, of Hopf algebras that takes the value
$\mathcal{O}{}(\GL_{n})$ for $q=1$. This can be regarded as a one-parameter
deformation of $\mathcal{O}{}(\GL_{n})$ by Hopf algebras, or of $\GL_{n}$ by
quantum groups.

Artin, Schelter, and Tate (1991b) construct a family of deformations of
$\mathcal{O}{}(\GL_{n})$, depending on $1+\binom{n}{2}$ parameters, which
includes the family $\mathcal{O}{}(\GL_{n}(q))$. The algebras in the family
are all twists of $\mathcal{O}{}(\GL_{n}(q))$ by $2$-cocycles. They first
construct a family of deformations of $\mathcal{O(}{}M_{n})$ by bialgebras
that are graded algebras generated in degree $1$, have the same Hilbert series
as the polynomial ring in $n^{2}$ variables, and are noetherian domains. The
family of deformations of $\mathcal{O}{}(\GL_{n})$ is then obtained by
inverting the quantum determinant. The algebras in the family of deformations
of $\mathcal{O}{}(M_{n})$ are regular in the sense of (\ref{i1}), and so this
gives a large class of regular algebras with the expected good properties.

\subsection{Sklyanin algebras\label{i3}}

As noted in (\ref{i1}), regular algebras of degree $3$ over a field $k$
correspond to certain triples $(E,\mathcal{L}{},\sigma)$ with $E$ a curve,
$\mathcal{L}{}$ an invertible sheaf of degree $3$ on $E$, and $\sigma$ an
automorphism of $E$. When $E$ is a nonsingular elliptic curve and $\sigma$ is
translation by a point $P$ in $E(k)$, the algebra $A(E,\mathcal{L}{},\sigma)$
is called a \emph{Sklyanin algebra}. Let $U=$ $\Gamma(E,\mathcal{L}{})$. This
is a $3$-dimensional $k$-vector space, and we can identify $U\otimes U$ with
$\Gamma(E,\mathcal{L}{}\boxtimes\mathcal{L}{})$. The algebra $A(E,\mathcal{L}%
{},\sigma)$ is the quotient of the tensor algebra $T(U)$ of $U$ by $\{f\in
U\otimes U\mid f(x,x+P)=0\}$. It is essentially independent of $\mathcal{L}{}%
$, because any two invertible sheaves of degree $3$ differ by a translation.
More generally, there is Sklyanin algebra $A(E,\mathcal{L}{},\sigma)$ for
every triple consisting of a nonsingular elliptic curve, an invertible sheaf
$\mathcal{L}$ of degree $d$ on $E$, and a translation by a point in $E(k)$.
The algebra $A(E,\mathcal{L}{},\sigma)$ has dimension $d$.

Artin, Schelter, and Tate (1994c) give a precise description of the centres of
Sklyanin algebras of dimension three, and Smith and Tate (1994d) extend the
description to those of dimension four.

Tate and van den Bergh (1996) prove that every $d$-dimensional Sklyanin
algebra $A(E,\mathcal{L}{},\sigma)$ is a noetherian domain, is Koszul, has the
same Hilbert series as a polynomial ring in $d$ variables, and is regular in
the sense of (\ref{i1}); moreover, if $\sigma$ has finite order, then
$A(E,\mathcal{L}{},\sigma)$ is finite over its centre.

\section{Miscellaneous articles\label{misc}}

\subsubsection{1951\textup{a} Tate, John. On the relation between extremal
points of convex sets and homomorphisms of algebras. Comm. Pure Appl. Math. 4,
(1951). 31--32.}

Tate considers a convex set $K$ of linear functionals on a commutative algebra
$A$ over $\mathbb{R}{}$. Under certain hypotheses on $A$ and $K$, he proves
that the extremal points of $K$ are exactly the homomorphisms from $A$ into
$\mathbb{R}{}$.

\subsubsection{1951\textup{b} Artin, Emil; Tate, John T. A note on finite ring
extensions. J. Math. Soc. Japan 3, (1951). 74--77.}

\begin{minipage}{4in}
Artin and Tate prove that if $S$ is a commutative finitely generated algebra
over a noetherian ring $R$, and $T$ is a subalgebra of $S$ such that $S$ is
finitely generated as a $T$-module, then $T$ is also finitely generated over
$R$. This statement generalizes a lemma of Zariski, and is now known as the
Artin-Tate lemma.%
\index{Artin-Tate lemma}
There are various generalizations of it to noncommutative rings.
\end{minipage}\hfill\begin{minipage}{1in}
\begin{tikzpicture}[baseline=(current bounding box.center)]
\matrix(m)[matrix of math nodes, row sep=1.3em, column sep=2.5em,
text height=1.5ex, text depth=0.25ex]
{S\\
T\\
R\rlap{ \scriptsize{noetherian}}\\};
\path[-,font=\scriptsize]
(m-1-1) edge node[right] {finite} (m-2-1)
(m-2-1) edge node[right] {$\Rightarrow$ fg} (m-3-1)
(m-1-1) edge [bend right] node[left] {fg} (m-3-1);
\end{tikzpicture}
\end{minipage}

\subsubsection{1952\textup{a} Tate, John. Genus change in inseparable
extensions of function fields. Proc. Amer. Math. Soc. 3, (1952). 400--406.}

Let $C$ be a complete normal geometrically integral curve over a field $k$ of
characteristic $p$, and let $C^{\prime}$ be the curve obtained from $C$ by an
extension of the base field $k\rightarrow k^{\prime}$. If $k^{\prime}$ is
inseparable over $k$, then $C^{\prime}$ need not be normal, and its
normalization $\tilde{C}^{\prime}$ may have genus $g(\tilde{C}^{\prime})$ less
than the genus $g(C)$ of $C$. However, Tate proves that
\begin{equation}
\frac{p-1}{2}\text{ divides }g(C)-g(\tilde{C}^{\prime}). \label{e7}%
\end{equation}
In particular, the genus of $C$ can't change if $g(C)<$ $(p-1)/2$ (which
implies that $C$ is smooth in this case).

Statement (\ref{e7}) is widely used. Tate derives it from a \textquotedblleft
Riemann-Hurwitz formula\textquotedblright\ for purely inseparable coverings,
which he proves using the methods of the day (function fields and
repartitions). A modern proof has been given of (\ref{e7}%
)\footnote{Schr\"{o}er, Stefan, On genus change in algebraic curves over
imperfect fields. Proc. Amer. Math. Soc. 137 (2009), no. 4, 1239--1243.}, but
not, as far as I know, of the more general formula.

\subsubsection{1952\textup{b} Lang, Serge; Tate, John. On Chevalley's proof of
Luroth's theorem. Proc. Amer. Math. Soc. 3, (1952). 621--624.}

Chevalley (1951, p.~106)\footnote{Chevalley, Claude. Introduction to the
Theory of Algebraic Functions of One Variable. Mathematical Surveys, No. VI.
American Mathematical Society, New York, N. Y., 1951.} proved L\"{u}roth's
theorem\ in the following form: let $k_{0}$ be a field, and let $k=k_{0}(X)$
be the field of rational functions in the symbol $X$ (i.e., $k$ is the field
of fractions of the polynomial ring $k_{0}[X]$); then every intermediate field
$k^{\prime}$, $k_{0}\subsetneqq k^{\prime}\subset k$, is of the form
$k_{0}(f)$ for some $f\in k$.

A classical proof of L\"{u}roth's theorem uses the Riemann-Hurwitz formula.
Let $k$ be a function field in one variable over a field $k_{0}$, and let
$k^{\prime}$ be an intermediate field; the Riemann-Hurwitz formula shows that,
if $k/k^{\prime}$ is separable, then%
\[
g(k^{\prime})\leq g(k).
\]
Therefore, if $k$ has genus zero, so also does $k^{\prime}$; if, in addition,
$k$ has a prime of degree $1$, so also does $k^{\prime}$, and so $k^{\prime}$
is a rational field (by a well-known criterion).

However, if $k/k^{\prime}$ is not separable, it may happen that $g(k^{\prime
})>g(k)$. Chevalley proved L\"{u}roth's theorem in nonzero characteristic by
showing directly that, when $k=k_{0}(X)$, every intermediate field $k$ has
genus zero. Lang and Tate generalized Chevalley's argument to prove:

\begin{quote}
Let $k$ be a function field in one variable over a field $k_{0}$, and let
$k^{\prime}$ be an intermediate field; if $k$ is separably generated over
$k_{0}$, then $g(k^{\prime})\leq g(k).$
\end{quote}

\noindent In other words, they showed that Chevalley's argument doesn't
require that $k=k_{0}(X)$ but only that it be separably generated over $k_{0}%
$. They also prove a converse statement:

\begin{quote}
A field of genus zero that is not separably generated over its field of
constants contains subfields of arbitrarily high genus.
\end{quote}

\noindent Finally, to complete their results, they exhibit a field of genus
zero, not separably generated over its field of constants.

\subsubsection{1955\textup{b} Brauer, Richard; Tate, John. On the characters
of finite groups. Ann. of Math. (2) 62, (1955). 1--7.}

\label{1955b}Recall (p.\pageref{Brauer1947}) that Brauer's theorem says that
every character $\chi$ of a finite group $G$ can be expressed in the form%
\[
\chi=\sum\nolimits_{i}n_{i}\Ind\chi_{i},\quad n_{i}\in\mathbb{Z}{},
\]
with the $\chi_{i}$ \textit{one-dimensional} characters on subgroups of $G$
(as conjectured by Artin). Brauer and Tate found what is probably the simplest
known proof of Brauer's theorem. Recall that a group is said to be elementary
if it can be expressed as the product of a cyclic group with a $p$-group for
some prime $p$. An elementary group is nilpotent, and so every irreducible
character of it is induced from a one-dimensional character on a subgroup. Let
$G$ be a finite group, and let $\mathcal{H}{}$ be a set of subgroups of $G$.
Consider the following three $\mathbb{Z}{}$-submodules of the space of
complex-valued functions on $G$:%
\begin{align*}
X(G)  &  =\mathrm{span}\{\text{irreducible characters of }G\}\text{ }%
\quad\quad\text{(module of virtual characters)}\\
Y  &  =\mathrm{span}\{\text{characters of }G\text{ induced from an irreducible
character of an }H\text{ in }\mathcal{H}\text{{}}\}\\
U  &  =\{\text{class functions }\chi\text{ on }G\text{ such that }\chi|H\in
X(H)\text{ for all }H\text{ in }\mathcal{H}{}\}.
\end{align*}
Brauer and Tate show that%
\[
U\supset X(G)\supset Y\text{,}%
\]
that $U$ is a ring, and that $Y$ is an ideal of $U$. Using this, they show
that if $\mathcal{H}{}$ consists of the elementary subgroups of $G$, then
$U=Y$, thereby elegantly proving not only Artin's conjecture (the equality
$X(G)=Y$), but also the main theorem of Brauer 1953\footnote{Brauer, Richard,
A characterization of the characters of groups of finite order. Ann. of Math.
(2) 57, (1953). 357--377.} (the equality $U=X(G)$).

\subsubsection{1957 Tate, John. Homology of Noetherian rings and local rings.
Illinois J. Math. 1 (1957), 14--27.}

Tate makes systematic use of the skew-commutative graded differential algebras
over a noetherian commutative ring $R$ to obtain results concerning $R$ and
its quotient rings${}$. The differential of such an $R$-algebra allows it to
be regarded as a complex, and Tate proves that every quotient $R/\mathfrak{a}%
{}$ of $R$ has a free resolution that is an $R$-algebra (in the above sense).
Such resolutions are now called Tate resolutions.%
\index{Tate resolutions}%

Let $R$ be a local noetherian ring with maximal ideal $\mathfrak{m}{}$. The
Betti series of $R$ is defined to be the formal power series $\mathcal{R}%
{}=\sum_{r\geq0}b_{r}Z^{r}$ with $b_{r}$ equal to the length of $\mathrm{Tor}%
_{r}^{R}(R/\mathfrak{m}{},R/\mathfrak{m}{})$. Serre (1956)\footnote{Serre,
Jean-Pierre, Sur la dimension homologique des anneaux et des modules
noeth\'{e}riens. Proceedings of the international symposium on algebraic
number theory, Tokyo \& Nikko, 1955, pp. 175--189. Science Council of Japan,
Tokyo, 1956.} showed that $R$ is regular if and only if $\mathcal{R}{}$ is a
polynomial, in which case $\mathcal{R}{}=(1+Z)^{d}$ with $d=\dim(R)$. Tate
showed that $\mathcal{R}{}=(1+Z)^{d}/(1-Z^{2})^{b_{1}-d}$ if $R$ is a complete
intersection. In general, he showed that the natural homomorphism%
\[
\bigwedge\nolimits^{\ast}\mathrm{Tor}_{1}^{R}(R/\mathfrak{m}{},R/\mathfrak{m}%
{})\rightarrow\mathrm{Tor}^{R}(R/\mathfrak{m}{},R/\mathfrak{m}{})
\]
is injective and realizes $\mathrm{Tor}^{R}(R/\mathfrak{m}{},R/\mathfrak{m}%
{})$ as a free module over $\bigwedge\nolimits^{\ast}\mathrm{Tor}_{1}%
^{R}(R/\mathfrak{m}{},R/\mathfrak{m}{})$ with a homogeneous basis. If $R$ is
regular, then the homomorphism is an isomorphism; conversely, if the
homomorphism is an isomorphism on one homogeneous component of degree $\geq2$,
then $R$ is regular.

\subsubsection{1962\textup{a} Fr\"{o}hlich, A.; Serre, J.-P.; Tate, J. A
different with an odd class. J. Reine Angew. Math. 209 1962 6--7.}

Let $A$ be a Dedekind domain with field of fractions $K$, and let $B$ be the
integral closure of $A$ in a finite separable extension of $K$. The different
$\mathfrak{D}{}$ of $B/A$ is an ideal in $B$, and its norm $\mathfrak{d}{}$ is
the discriminant ideal of $B/A$. The ideal class of $\mathfrak{d}{}$ is always
a square, and Hecke (1954, \S 63)\footnote{Hecke, Erich, Vorlesungen \"{u}ber
die Theorie der algebraischen Zahlen. 2te Aufl. (German) Akademische
Verlagsgesellschaft, Geest \& Portig K.-G., Leipzig, 1954.} proved that the
ideal class of $\mathfrak{D}{}$ is a square when $K$ is a number field, but
the authors show that it need not be a square otherwise. Specifically, they
construct examples of affine curves over perfect fields whose coordinate rings
$A$ have extensions $B$ for which the ideal class of the different is not a
square.\footnote{Hecke's theorem can be proved for global fields of
characteristic $p\neq0$ by methods similar to those of Hecke (Armitage, J. V.,
On a theorem of Hecke in number fields and function fields. Invent. Math. 2
1967 238--246).}

This is not major result. However, Martin Taylor\footnote{Taylor, M. J.
Obituary: Albrecht Fr\"{o}hlich, 1916--2001. Bull. London Math. Soc. 38
(2006), no. 2, 329--350.} writes:

\begin{quote}
[This article and Fr\"{o}hlich's earlier work on discriminants] seems to have
marked the start of [his] interest in parity questions. He would go on to be
interested in whether conductors of real-valued characters were squares; this
in turn led to questions about the signs of Artin root numbers --- an issue
that lay right at the heart of his work on Galois modules.
\end{quote}

\noindent Fr\"{o}hlich's work on Artin root numbers and Galois module
structures was his most important.

\subsubsection{1963 Sen, Shankar; Tate, John. Ramification groups of local
fields. J. Indian Math. Soc. (N.S.) 27 1963 197--202 (1964).}

Let $F$ be a field, complete with respect to a discrete valuation, and let $K$
be a finite Galois extension of $F$. Assume initially that the residue field
is finite, and let $W$ be the Weil group of $K/F$ (extension of $G(K/F)$ by
$K^{\times}$ determined by the fundamental class of $K/F$). Shafarevich showed
that there is a homomorphism $s$ making the following diagram commute%
\[
\begin{CD}
1@>>>K^{\times}@>>>W@>>>G(K/F)@>>>1\\
@.@VVrV@VVsV@|@.\\
1@>>>G(K^{\mathrm{ab}}/K)@>>>G(K^{\mathrm{ab}}/F)@>>>G(K/F)@>>>1,
\end{CD}
\]
where $r$ is the reciprocity map. For a real $t>0$, let $G(K^{\mathrm{ab}%
}/K)^{t}$ denote the $t$th ramification subgroup of $G(K^{\mathrm{ab}}/K)$.
Then%
\begin{equation}
r^{-1}(G(K^{\mathrm{ab}}/K)^{t})=U_{K}^{t}\overset{\textup{{\tiny def}}}%
{=}\left\{
\begin{array}
[c]{l}%
\{u\in K^{\times}\mid\ord_{K}(u)=0\}\text{ if }t=0\\
\{u\in K^{\times}\mid\ord_{K}(u-1)\geq t\}\text{ if }t>0.
\end{array}
\right.  \label{e24}%
\end{equation}

Artin and Tate (1961) proved the existence of the Shafarevich map $s$ for a
general class formation. When the residue field of $F$ is algebraically
closed, the groups $\pi_{1}(U_{K})$ (fundamental group of $U_{K}$ regarded as
a pro-algebraic group) form a class formation, and so the above diagram exists
with $K^{\times}$ replaced by $\pi_{1}(U_{K})$. In this case,
\begin{equation}
r^{-1}(G(K^{\mathrm{ab}}/K)^{t})=\pi_{1}(U_{K}^{t}). \label{e25}%
\end{equation}

In both cases, Sen and Tate give a description of the subgroups $s^{-1}%
(G(K^{\mathrm{ab}}/F)^{t})$ of $W$ generalizing those in (\ref{e24}) and
(\ref{e25}), which can be considered the case $K=F$. Specifically, let
$G(K/F)_{x}$, $x\geq0$, denote with ramification groups of $K/F$ with the
lower numbering, and let%
\[
\varphi(x)=\int_{0}^{x}\frac{du}{(G(K/F)_{0}\colon G(K/F)_{u})}\text{ for
}x\geq0.
\]
For $u\in W$, let $m(u)>0$ be the smallest integer such that $u^{m(u)}\in
K^{\times}$ (resp. $\pi_{1}(U_{K})$). Then%
\[
W^{\varphi(x)}=\left\{  u\in W\,\middle|\, u^{m(u)}\in U_{K}^{m(u)\cdot x}\text{ (resp. 
}\pi_{1}(U_{K}^{m(u)\cdot x})\text{)}\right\}  \text{.}%
\]

\subsubsection{1964\textup{c} Tate, John. Nilpotent quotient groups. Topology
3 1964 suppl. 1 109--111.}

For a finite group $G$, subgroup $S$, and positive integer $p,$ there are
restriction maps $r$ and transfer maps $t$,
\[
H^{i}(G,\mathbb{Z}{}/p\mathbb{Z}{})\overset{r^{i}}{\longrightarrow}%
H^{i}(S,\mathbb{Z}{}/p\mathbb{Z}{})\overset{t^{i}}{\longrightarrow}%
H^{i}(G,\mathbb{Z}{}/p\mathbb{Z}{}),\quad i\geq0,
\]
whose composite is multiplication by $(G\colon S)$.

Let $S$ be Sylow $p$-subgroup of $G$ (so $p$ is prime). If $S$ has a normal
$p$-complement in $G$, then the restriction maps are isomorphisms, and Atiyah
asked whether the converse is true. Thompson pointed out that the answer is
yes, and that results of his and Huppert show that one need only require that
$r^{1}$ is an isomorphism. Tate gives a very short cohomological proof of a
somewhat stronger result.

Specifically, for a finite group $G$, define a descending sequence of normal
subgroups of $G$ as follows:%
\[
G_{0}=G,\quad G_{n+1}=\left(  G_{n}\right)  ^{p}[G,G_{n}]\text{ for }%
n\geq0,\quad G_{\infty}=\bigcap\nolimits_{n=0}^{\infty}G_{n}%
\]
($p$ not necessarily prime). Thus, $G/G_{1}$ (resp. $G/G_{\infty}$) is the
largest quotient group of $G$ that is abelian of exponent $p$ (resp. nilpotent
and $p$-primary). Let $S$ be a subgroup of $G$ of index prime to $p$. The
following three conditions are (obviously) equivalent,

\begin{itemize}
\item the restriction map $r^{1}\colon H^{1}(G,\mathbb{Z}{}/p\mathbb{Z}%
{})\rightarrow H^{1}(S,\mathbb{Z}{}/p\mathbb{Z}{})$ is an isomorphism,

\item the map $S/S_{1}\rightarrow G/G_{1}$ is an isomorphism,

\item $S\cap G_{1}=S_{1}$,
\end{itemize}

\noindent and Tate proves that they imply

\begin{itemize}
\item $S\cap G_{n}=S_{n}$ for all $1\leq n\leq\infty$.
\end{itemize}

\noindent When $S$ is a Sylow $p$-subgroup of $G$, $S\cap G_{\infty}=1$, and
so the conditions imply that $G_{\infty}$ is a normal $p$-complement of $S$ in
$G$ (thereby recovering the Huppert-Thompson theorem).

\subsubsection{1968\textup{b} Tate, John. Residues of differentials on curves.
Ann. Sci. \'{E}cole Norm. Sup. (4) 1 1968 149--159.}

Tate defines the residues of differentials on curves as the traces of certain
\textquotedblleft finite potent\textquotedblright\ linear maps. From his
definition, all the standard theorems on residues follow naturally and easily.
In particular, the residue formula%
\[
\sum\nolimits_{P\in C}\res_{p}(\omega)=0\quad\quad\text{(}C\text{ a complete
curve)}%
\]
follows directly, without computation, from the finite dimensionality of the
cohomology groups $H^{0}(C,\mathcal{O}{}_{C})$ and $H^{1}(C,\mathcal{O}{}%
_{C})$ \textquotedblleft almost as though one had an abstract Stokes's Theorem
available\textquotedblright.

A linear map $\theta\colon V\rightarrow V$ is finite potent if $\theta^{n}V$
is finite dimensional for some $n$. The trace $\mathrm{Tr}_{V}(\theta)$ of
such a map can be defined to be its trace on any finite dimensional subspace
$W$ of $V$ such that $\theta W\subset W$ and $\theta^{n}V\subset W$ for some
$n$. Many of the properties of the usual trace continue to hold, but not all.
For example, there exist finite potent maps such that\footnote{Pablos Romo,
Fernando, On the linearity property of Tate's trace. Linear Multilinear
Algebra 55 (2007), no. 4, 323--326. Argerami, Martin; Szechtman, Fernando;
Tifenbach, Ryan, On Tate's trace. Linear Multilinear Algebra 55 (2007), no. 6,
515--520.}%
\[
\mathrm{Tr}_{V}(\theta_{1}+\theta_{2}+\theta_{3})\neq\mathrm{Tr}_{V}%
(\theta_{1})+\mathrm{Tr}_{V}(\theta_{2})+\mathrm{Tr}_{V}(\theta_{3}).
\]
Tate defines the residue of a differential $f\,dg$ at a closed point $p$ of a
curve $C$ to be the trace of the commutator $[f_{p},g_{p}]$, where $f_{p}$,
$g_{p}$ are representatives of $f$, $g$ in a certain subspace of
$\End(k(C)_{p})$.

Tate's approach to residues has found its way into the text books (e.g.,
Iwasawa 1993\footnote{Iwasawa, Kenkichi Algebraic functions. Translated from
the 1973 Japanese edition by Goro Kato. Translations of Mathematical
Monographs, 118. American Mathematical Society, Providence, RI, 1993.}).
Elzein\footnote{Elzein, Fouad, R\'{e}sidus en g\'{e}om\'{e}trie
alg\'{e}brique. C. R. Acad. Sci. Paris S\'{e}r. A-B 272 (1971), A878--A881.}
used Tate's ideas to give a definition of the residue that recaptures both
Leray's in the case of a complex algebraic variety and Grothendieck's in the
case of a smooth integral morphisms of relative dimension $n$.

Others have adapted his proof of the residue formula to other situations; for
example, Arbarello et al (1989)\footnote{Arbarello, E.; De Concini, C.; Kac,
V. G. The infinite wedge representation and the reciprocity law for algebraic
curves. Theta functions---Bowdoin 1987, Part 1 (Brunswick, ME, 1987),
171--190, Proc. Sympos. Pure Math., 49, Part 1, Amer. Math. Soc., Providence,
RI, 1989.} use it to prove an \textquotedblleft abstract reciprocity
law\textquotedblright\ for tame symbols on a curve over an algebraically
closed field, and Beilinson et al (2002)\footnote{Beilinson, Alexander; Bloch,
Spencer; Esnault, H\'{e}l\`{e}ne, $\varepsilon$-factors for Gauss-Manin
determinants. Dedicated to Yuri I. Manin on the occasion of his 65th birthday.
Mosc. Math. J. 2 (2002), no. 3, 477--532.} use it to prove a product formula
for $\varepsilon$-factors in the de Rham setting.

In reading Tate's article, Beilinson recognized that a certain linear algebra
construction there can be reformulated as the construction of a canonical
central extension of Lie algebras. This led to the notion of a Tate extension%
\index{Tate extension}
in various settings; see Beilinson and Schechtman
1988,\footnote{Be\u{\i}linson, A. A.; Schechtman, V. V. Determinant bundles
and Virasoro algebras. Comm. Math. Phys. 118 (1988), no. 4, 651--701.} and
Beilinson and Drinfeld 2004, 2.7.\footnote{Beilinson, Alexander; Drinfeld,
Vladimir. Chiral algebras. American Mathematical Society Colloquium
Publications, 51. American Mathematical Society, Providence, RI, 2004.}

\subsubsection{1978\textup{a} Cartier, P.; Tate, J. A simple proof of the main
theorem of elimination theory in algebraic geometry. Enseign. Math. (2) 24
(1978), no. 3-4, 311--317.}

The authors give an elementary one-page proof of the homogeneous form of
Hilbert's theorem of zeros:

\begin{quote}
let $\mathfrak{a}{}$ be a graded ideal in a polynomial ring $k[X_{0}%
,\ldots,X_{n}]$ over a field $k$; either the radical of $\mathfrak{a}$
contains the ideal $(X_{0},\ldots,X_{n})$, or $\mathfrak{a}{}$ has a
nontrivial zero in an algebraic closure of $k$.
\end{quote}

\noindent From this, they quickly deduce the main theorem of elimination
theory, both in its classical form and in its modern form:

\begin{quote}
let $A=\bigoplus\nolimits_{d\geq0}A_{d}$ be a graded commutative algebra such
that $A$ is generated as an $A_{0}$-algebra by $A_{1}$ and each $A_{0}$-module
$A_{d}$ is finitely generated; then the map of topological spaces
$\mathrm{proj}(A)\rightarrow\mathrm{spec(}A_{0})$ is closed.
\end{quote}

\subsubsection{1989 Gross, B.; Tate, J. Commentary on algebra. A century of
mathematics in America, Part II, 335--336, Hist. Math., 2, Amer. Math. Soc.,
Providence, RI, 1989.}

For the bicentenary of Princeton University in 1946, there was a three-day
conference in which various distinguished mathematicians discussed Problems in
Mathematics. Artin, Brauer, and others contributed to the discussion on
algebra, and in 1989 Gross and Tate wrote a commentary on their remarks. For example:

\begin{quote}
Artin's belief that \textquotedblleft whatever can be said about non-Abelian
class field theory follows from what we know now,\textquotedblright\ and that
\textquotedblleft our difficulty is not in the proofs, but in learning what to
prove,\textquotedblright\ seems overly optimistic.
\end{quote}

\subsubsection{1994\textup{b} Tate, John. The non-existence of certain Galois
extensions of $\mathbb{Q}{}$ unramified outside $2$. Arithmetic geometry
(Tempe, AZ, 1993), 153--156, Contemp. Math., 174, Amer. Math. Soc.,
Providence, RI, 1994.}

In a 1973 letter to Tate, Serre suggested that certain two-dimensional mod $p$
representations of $\Gal(\mathbb{Q}{}^{\mathrm{al}}/\mathbb{Q}{})$ should be
modular. In response, Tate verified this for $p=2$ by showing that every
two-dimensional mod $2$ representation unramified outside $2$ has zero trace.
The article is based on his letter.

Serre's suggestion became Serre's conjecture on the modularity of
two-dimensional mod $p$ representations, which attracted much attention
because of its relation to the modularity conjecture for elliptic curves over
$\mathbb{Q}{}$ and Fermat's last theorem. Serre's conjecture was recently
proved by an inductive argument that uses Tate's result as one of the base
cases.\footnote{Khare, Chandrashekhar; Wintenberger, Jean-Pierre. Serre's
modularity conjecture. I. Invent. Math. 178 (2009), no. 3, 485--504. II. Ibid.
505--586.}

\subsubsection{1996\textup{a} Tate, John; Voloch, Jos\'{e} Felipe . Linear
forms in $p$-adic roots of unity. Internat. Math. Res. Notices 1996, no. 12,
589--601.}

The authors make the following conjecture: for a semi-abelian variety $A$ over
$\mathbb{C}{}_{p}$ and a closed subvariety $X$, there exists a lower bound
$c>0$ for the $p{}$-adic distance of torsion points of $A$, not in $X$, to
$X$. Here, as usual, $\mathbb{C}{}_{p}$ is the completion of an algebraic
closure of $\mathbb{Q}{}_{p}$. They prove the conjecture for the torus
\[
A=\Spec\mathbb{C}{}_{p}[T_{1},T_{1}^{-1},\ldots,T_{n},T_{n}^{-1}].
\]
This comes down to proving the following explicit statement: for every
hyperplane%
\[
a_{1}T_{1}+\cdots+a_{n}T_{n}=0
\]
in $\mathbb{C}_{p}^{n}$, there exists a constant $c>0$, depending on
$(a_{1},\ldots,a_{n})$, such that, for any $n$-tuple $\zeta_{1},\ldots
,\zeta_{n}$ of roots of $1$ in $\mathbb{C}{}_{p}$, either $a_{1}\zeta
_{1}+\cdots+a_{1}\zeta_{n}=0$ or $|a_{1}\zeta_{1}+\cdots+a_{1}\zeta_{n}|\geq
c$.

\subsubsection{2002\textup{a} Tate, John. On a conjecture of Finotti. Bull.
Braz. Math. Soc. (N.S.) 33 (2002), no. 2, 225--229.}

In his study of the Teichm\"{u}ller points in canonical lifts of elliptic
curves, Finotti was led to a conjecture on remainders of division by
polynomials.\footnote{L. R. A. Finotti, Canonical and minimal degree liftings
of curves, J. Math. Sci. Univ. Tokyo 11 (2004), no. 1, 1--47 (Ph.D. thesis,
Univ. Texas, 2001).} He checked it by computer for all primes $p\leq877$, and
Tate proved it in general. The statement is:

\begin{quote}
Let $k$ be a field of characteristic $p=2m+1\geq5$. Let $F\in k[X]$ be a monic
cubic polynomial, and let $A$ be the coefficient of $X^{p-1}$ in $F^{m}$. Let
$G\in k[X]$ be a polynomial of degree $3m+1$ such that $G^{\prime}%
=F^{m}-AX^{p-1}$. Then the remainder in the division of $G^{2}$ by
$X^{p}F^{m+1}$ has degree $\leq5m+2=\frac{5p-1}{2}$.
\end{quote}

\subsection*{Acknowledgements}

I thank B. Gross for help with dates, J-P. Serre for correcting a
misstatement, and J. Tate for answering my queries and pointing out some mistakes.

\subsection*{Added September 2012}
I should have mentioned the work of Tate on liftings of Galois representations, as included in Part II of: Serre, J.-P. Modular forms of weight one and Galois representations. Algebraic number fields: $L$-functions and Galois properties (Proc. Sympos., Univ. Durham, Durham, 1975), pp. 193--268. Academic Press, London, 1977.
See also: Variations on a theorem of Tate. Stefan Patrikis. arXiv:1207.6724.

Also, ``An oft cited (1979) letter from Tate to Serre on computing local heights on elliptic curves.''
was posted on the arXiv by Silverman (arXiv:1207.5765).

The collected works of Tate, which will include other unpublished letters, is in preparation. 

\clearpage \markboth{}{}

\section*{Bibliography of Tate's articles}
\addcontentsline{toc}{section}{Bibliography}

\setlength{\parindent}{-2em}

\section*{1950s}

\hspace{-2em}1950 Tate, John, Fourier Analysis in Number Fields and Hecke's
Zeta Functions, Ph.D. thesis, Princeton University. Published as 1967b.

1951a Tate, John. On the relation between extremal points of convex sets and
homomorphisms of algebras. Comm. Pure Appl. Math. 4, (1951). 31--32.

1951b Artin, Emil; Tate, John T. A note on finite ring extensions. J. Math.
Soc. Japan 3, (1951). 74--77.

1952a Tate, John. Genus change in inseparable extensions of function fields.
Proc. Amer. Math. Soc. 3, (1952). 400--406.

1952b Lang, Serge; Tate, John. On Chevalley's proof of Luroth's theorem. Proc.
Amer. Math. Soc. 3, (1952). 621--624.

1952c Tate, John. The higher dimensional cohomology groups of class field
theory. Ann. of Math. (2) 56, (1952). 294--297.

1954 Tate, John, The Cohomology Groups of Algebraic Number Fields, pp. 66-67
in Proceedings of the International Congress of Mathematicians, Amsterdam,
1954. Vol. 2. Erven P. Noordhoff N. V., Groningen; North-Holland Publishing
Co., Amsterdam, 1954. iv+440 pp.19

1955a Kawada, Y.; Tate, J. On the Galois cohomology of unramified extensions
of function fields in one variable. Amer. J. Math. 77, (1955). 197--217.

1955b Brauer, Richard; Tate, John. On the characters of finite groups. Ann. of
Math. (2) 62, (1955). 1--7.

1957 Tate, John. Homology of Noetherian rings and local rings. Illinois J.
Math. 1 (1957), 14--27.

1958a Mattuck, Arthur; Tate, John. On the inequality of Castelnuovo-Severi.
Abh. Math. Sem. Univ. Hamburg 22 1958 295--299.

1958b Tate, J. WC-groups over $\mathfrak{p}{}$-adic fields. S\'{e}minaire
Bourbaki; 10e ann\'{e}e: 1957/1958. Textes des conf\'{e}rences; Expos\'{e}s
152 \`{a} 168; 2e \'{e}d. corrig\'{e}e, Expos\'{e} 156, 13 pp. Secr\'{e}tariat
math\'{e}matique, Paris 1958 189 pp (mimeographed).

1958c Lang, Serge; Tate, John. Principal homogeneous spaces over abelian
varieties. Amer. J. Math. 80 1958 659--684.

1958d Tate, John, Groups of Galois Type (published as Chapter VII of Lang
1967; reprinted as Lang 1996).\footnote{Lang calls his Chapter VII an
\textquotedblleft unpublished article of Tate\textquotedblright, but gives no
date. In his MR review, Shatz writes that \textquotedblleft It appears
here in almost the same form the reviewer remembers from the original seminar
of Tate in 1958. \textquotedblright}

1959a Tate, John. Rational points on elliptic curves over complete fields,
manuscript 1959. Published as part of 1995.

1959b. Tate, John. Applications of Galois cohomology in algebraic
geometry. (Written by S. Lang based on letters of Tate 1958--1959).
Chapter X of: Lang, Serge. Topics in cohomology of groups. Translated
from the 1967 French original by the author. Lecture Notes in
Mathematics, 1625. Springer-Verlag, Berlin, 1996. vi+226 pp.

\section*{1960s}

\hspace{-2em}1961 Artin, E., and Tate., J. Class Field Theory, Harvard
University, Department of Mathematics, 1961.\footnote{The original notes don't
give a date or a publisher. I copied this information from the footnote p. 162
of 1967a. The volume was prepared by the staff of
the Institute of Advanced Study, but it was distributed by the Harvard
University Mathematics Department.} Notes from the Artin-Tate seminar on class
field theory given a Princeton University 1951--1952. Reprinted as 1968c,
1990b; second edition 2009.

1962a Fr\"{o}hlich, A.; Serre, J.-P.; Tate, J. A different with an odd class.
J. Reine Angew. Math. 209 1962 6--7.

1962b Tate, John. Principal homogeneous spaces for Abelian varieties. J. Reine
Angew. Math. 209 1962 98--99.

1962c Tate, John, Rigid analytic spaces. Private notes, reproduced with(out)
his permission by I.H.E.S (1962). Published as 1971b; Russian translation 1969a.

1962d Tate, John. Duality theorems in Galois cohomology over number fields.
1963 Proc. Internat. Congr. Mathematicians (Stockholm, 1962) pp. 288--295
Inst. Mittag-Leffler, Djursholm

1963 Sen, Shankar; Tate, John. Ramification groups of local fields. J. Indian
Math. Soc. (N.S.) 27 1963 197--202 (1964).

1964a Tate, John. Algebraic cohomology classes, Woods Hole 1964, 25 pages. In:
Lecture notes prepared in connection with seminars held at the Summer
Institute on Algebraic Geometry, Whitney Estate, Woods Hole, MA, July 6--July
31, 1964. Published as 1965b; Russian translation 1965c.

1964b Tate, John (with Lubin and Serre). Elliptic curves and formal groups, 8
pages. In: Lecture notes prepared in connection with seminars held at the
Summer Institute on Algebraic Geometry, Whitney Estate, Woods Hole, MA, July
6--July 31, 1964.

1964c Tate, John. Nilpotent quotient groups. Topology 3 1964 suppl. 1 109--111.

1965a Lubin, Jonathan; Tate, John. Formal complex multiplication in local
fields. Ann. of Math. (2) 81 1965 380--387.

1965b Tate, John T. Algebraic cycles and poles of zeta functions. 1965
Arithmetical Algebraic Geometry (Proc. Conf. Purdue Univ., 1963) pp. 93--110
Harper \& Row, New York

1965c Tate, John. Algebraic cohomology classes. (Russian) Uspehi Mat. Nauk 20
1965 no. 6 (126) 27--40.

1965d Tate, John. Letter to Cassels on elliptic curve formulas. Published as 1975b. 

1966a Tate, John. Multiplication complexe formelle dans les corps locaux. 1966
Les Tendances G\'{e}om. en Alg\`{e}bre et Th\'{e}orie des Nombres pp. 257--258
\'{E}ditions du Centre National de la Recherche Scientifique, Paris.

1966b Tate, John. Endomorphisms of abelian varieties over finite fields.
Invent. Math. 2 1966 134--144.

1966c Tate, J. The cohomology groups of tori in finite Galois extensions of
number fields. Nagoya Math. J. 27 1966 709--719.

1966d Lubin, Jonathan; Tate, John. Formal moduli for one-parameter formal Lie
groups. Bull. Soc. Math. France 94 1966 49--59.

1966e Tate, John T. On the conjectures of Birch and Swinnerton-Dyer and a
geometric analog. 1966. S\'{e}minaire Bourbaki: Vol. 1965/66, Expose 306.

1966f Tate, John. Letter to Springer, January 13, 1966. (Contains proofs of
some of the theorems announced in 1963a.)

1967a Tate, J. T. Global class field theory. 1967 Algebraic Number Theory
(Proc. Instructional Conf., Brighton, 1965) pp. 162--203 Thompson, Washington, D.C.

1967b Tate, J. T. Fourier analysis in number fields and Hecke's
zeta-functions. 1967 Algebraic Number Theory (Proc. Instructional Conf.,
Brighton, 1965) pp. 305--347 Thompson, Washington, D.C.

1967c Tate, J. T. $p$-divisible groups. 1967 Proc. Conf. Local Fields
(Driebergen, 1966) pp. 158--183 Springer, Berlin.

1967d Tate, John. Shafarevich, I. R. The rank of elliptic curves. (Russian)
Dokl. Akad. Nauk SSSR 175 1967 770--773.

1968a Serre, Jean-Pierre; Tate, John. Good reduction of abelian varieties.
Ann. of Math. (2) 88 1968 492--517.

1968b Tate, John. Residues of differentials on curves. Ann. Sci. \'{E}cole
Norm. Sup. (4) 1 1968 149--159.

1968c Artin, E.; Tate, J. Class field theory. W. A. Benjamin, Inc., New
York-Amsterdam 1968 xxvi+259 pp.

1969a Tate, John; Rigid analytic spaces. (Russian) Mathematics: periodical
collection of translations of foreign articles, Vol. 13, No. 3 (Russian), pp.
3--37. Izdat. \textquotedblleft Mir\textquotedblright, Moscow, 1969.

1969b Tate, John, Classes d'isog\'enie des vari\'et\'es ab\'eliennes sur un
corps fini (d'apr\`es T. Honda) S\'eminaire Bourbaki 352, (1968/69).

1969c Tate, John, $K_2$ of global fields, AMS Taped Lecture (Cambridge, Masss., Oct. 1969).

\section*{1970s}

\hspace{-2em} 1970a Tate, John; Oort, Frans. Group schemes of prime order.
Ann. Sci. \'{E}cole Norm. Sup. (4) 3 1970 1--21.

1970b Tate, John. Symbols in arithmetic. Actes du Congr\`es International des
Math\`ematiciens (Nice, 1970), Tome 1, pp. 201--211. Gauthier-Villars, Paris, 1971.

1971 Tate, John. Rigid analytic spaces. Invent. Math. 12 (1971), 257--289.

1973a Bass, H.; Tate, J. The Milnor ring of a global field. Algebraic
K-theory, II: \textquotedblleft Classical\textquotedblright\ algebraic
K-theory and connections with arithmetic (Proc. Conf., Seattle, Wash.,
Battelle Memorial Inst., 1972), pp. 349--446. Lecture Notes in Math., Vol.
342, Springer, Berlin, 1973.

1973b Tate, J. Letter from Tate to Iwasawa on a relation between $K_2$ and Galois
cohomology. Algebraic K-theory, II: \textquotedblleft
Classical\textquotedblright\ algebraic K-theory and connections with
arithmetic (Proc. Conf., Seattle Res. Center, Battelle Memorial Inst., 1972),
pp. 524--527. Lecture Notes in Math., Vol. 342, Springer, Berlin, 1973.

1973c Mazur, B.; Tate, J. Points of order 13 on elliptic curves. Invent. Math.
22 (1973/74), 41--49.

1974a Tate, John T. The arithmetic of elliptic curves. Invent. Math. 23
(1974), 179--206.

1974b Tate, J. The 1974 Fields medals. I. An algebraic geometer. Science 186
(1974), no. 4158, 39--40.

1975a Tate, J. The work of David Mumford. Proceedings of the International
Congress of Mathematicians (Vancouver, B. C., 1974), Vol. 1, pp. 11--15.
Canad. Math. Congress, Montreal, Que., 1975.

1975b Tate, J. Algorithm for determining the type of a singular fiber in an
elliptic pencil. Modular functions of one variable, IV (Proc. Internat. Summer
School, Univ. Antwerp, Antwerp, 1972), pp. 33--52. Lecture Notes in Math.,
Vol. 476, Springer, Berlin, 1975.

1976a Tate, J. Problem 9: The general reciprocity law. Mathematical
developments arising from Hilbert problems (Proc. Sympos. Pure Math., Northern
Illinois Univ., De Kalb, Ill., 1974), pp. 311--322. Proc. Sympos. Pure Math.,
Vol. XXVIII, Amer. Math. Soc., Providence, R. I., 1976.

1976b Tate, John. Relations between $K_{2}$ and Galois cohomology. Invent.
Math. 36 (1976), 257--274.

1977a Tate, J. On the torsion in $K_{2}$ of fields. Algebraic number theory
(Kyoto Internat. Sympos., Res. Inst. Math. Sci., Univ. Kyoto, Kyoto, 1976),
pp. 243--261. Japan Soc. Promotion Sci., Tokyo, 1977.

1977b Tate, J. T. Local constants. Prepared in collaboration with C. J.
Bushnell and M. J. Taylor. Algebraic number fields: L-functions and Galois
properties (Proc. Sympos., Univ. Durham, Durham, 1975), pp. 89--131. Academic
Press, London, 1977.

1978a Cartier, P.; Tate, J. A simple proof of the main theorem of elimination
theory in algebraic geometry. Enseign. Math. (2) 24 (1978), no. 3-4, 311--317.

1978b Tate, John. Fields medals. IV. Mumford, David; An instinct for the key
idea. Science 202 (1978), no. 4369, 737--739.

1979 Tate, J. Number theoretic background. Automorphic forms, representations
and L-functions (Proc. Sympos. Pure Math., Oregon State Univ., Corvallis,
Ore., 1977), Part 2, pp. 3--26, Proc. Sympos. Pure Math., XXXIII, Amer. Math.
Soc., Providence, R.I., 1979.

\section*{1980s}

\hspace{-2em} 1981a Tate, John. On Stark's conjectures on the behavior of
$L(s,\chi)$ at $s=0$. J. Fac. Sci. Univ. Tokyo Sect. IA Math. 28 (1981), no.
3, 963--978 (1982).

1981b Tate, John. Brumer-Stark-Stickelberger. Seminar on Number Theory,
1980--1981 (Talence, 1980--1981), Exp. No. 24, 16 pp., Univ. Bordeaux I,
Talence, 1981.

1981c Tate, John. On conjugation of abelian varieties of CM type. Handwritten
notes. 1981.

1983a Tate, J. Variation of the canonical height of a point depending on a
parameter. Amer. J. Math. 105 (1983), no. 1, 287--294.

1983b Mazur, B.; Tate, J. Canonical height pairings via biextensions.
Arithmetic and geometry, Vol. I, 195--237, Progr. Math., 35, Birkh\"{a}user
Boston, Boston, MA, 1983.

1983c Rosset, Shmuel; Tate, John. A reciprocity law for $K_{2}$-traces.
Comment. Math. Helv. 58 (1983), no. 1, 38--47.

1984 Tate, John. Les conjectures de Stark sur les fonctions $L$ d'Artin en
$s=0$. Notes of a course at Orsay written by Dominique Bernardi and Norbert Schappacher.
Progress in Mathematics, 47. Birkh\"{a}user Boston, Inc., Boston, MA, 1984.

1986 Mazur, B.; Tate, J.; Teitelbaum, J. On $p$-adic analogues of the
conjectures of Birch and Swinnerton-Dyer. Invent. Math. 84 (1986), no. 1, 1--48.

1987 Mazur, B.; Tate, J. Refined conjectures of the \textquotedblleft Birch
and Swinnerton-Dyer type\textquotedblright. Duke Math. J. 54 (1987), no. 2, 711--750.

1989 Gross, B.; Tate, J. Commentary on algebra. A century of mathematics in
America, Part II, 335--336, Hist. Math., 2, Amer. Math. Soc., Providence, RI, 1989.

\section*{1990s}

\hspace{-2em} 1990a Artin, M.; Tate, J.; Van den Bergh, M. Some algebras
associated to automorphisms of elliptic curves. The Grothendieck Festschrift,
Vol. I, 33--85, Progr. Math., 86, Birkh\"{a}user Boston, Boston, MA, 1990.

1990b Artin, Emil; Tate, John. Class field theory. Second edition. Advanced
Book Classics. Addison-Wesley Publishing Company, Advanced Book Program,
Redwood City, CA, 1990. xxxviii+259 pp. ISBN: 0-201-51011-1

1991a Artin, M.; Tate, J.; Van den Bergh, M. Modules over regular algebras of
dimension 3. Invent. Math. 106 (1991), no. 2, 335--388.

1991b Artin, Michael; Schelter, William; Tate, John. Quantum deformations of
$\mathrm{GL}_n$. Comm. Pure Appl. Math. 44 (1991), no. 8-9, 879--895.

1991c Mazur, B.; Tate, J. The $p$-adic sigma function. Duke Math. J. 62
(1991), no. 3, 663--688.

1992 Silverman, Joseph H.; Tate, John. Rational points on elliptic curves.
Undergraduate Texts in Mathematics. Springer-Verlag, New York, 1992. x+281 pp.

1994a Tate, John. Conjectures on algebraic cycles in $l$-adic cohomology.
Motives (Seattle, WA, 1991), 71--83, Proc. Sympos. Pure Math., 55, Part 1,
Amer. Math. Soc., Providence, RI, 1994.

1994b Tate, John. The non-existence of certain Galois extensions of
$\mathbb{Q}{}$ unramified outside $2$. Arithmetic geometry (Tempe, AZ, 1993),
153--156, Contemp. Math., 174, Amer. Math. Soc., Providence, RI, 1994.

1994c Artin, Michael; Schelter, William; Tate, John. The centers of
$3$-dimensional Sklyanin algebras. Barsotti Symposium in Algebraic Geometry
(Abano Terme, 1991), 1--10, Perspect. Math., 15, Academic Press, San Diego,
CA, 1994.

1994d Smith, S. P.; Tate, J. The center of the 3-dimensional and 4-dimensional
Sklyanin algebras. Proceedings of Conference on Algebraic Geometry and Ring
Theory in honor of Michael Artin, Part I (Antwerp, 1992). K-Theory 8 (1994),
no. 1, 19--63.

1995 Tate, John. A review of non-Archimedean elliptic functions. Elliptic
curves, modular forms, \& Fermat's last theorem (Hong Kong, 1993), 162--184,
Ser. Number Theory, I, Int. Press, Cambridge, MA, 1995.

1996a Tate, John; Voloch, Jos\'{e} Felipe. Linear forms in p-adic roots of
unity. Internat. Math. Res. Notices 1996, no. 12, 589--601.

1996b Tate, John; van den Bergh, Michel. Homological properties of Sklyanin
algebras. Invent. Math. 124 (1996), no. 1-3, 619--647.

1997a Tate, John. Finite flat group schemes. Modular forms and Fermat's last
theorem (Boston, MA, 1995), 121--154, Springer, New York, 1997.

1997b Tate, J. The work of David Mumford. Fields Medallists' lectures,
219--223, World Sci. Ser. 20th Century Math., 5, World Sci. Publ., River Edge,
NJ, 1997.

1999 Katz, Nicholas M.; Tate, John. Bernard Dwork (1923--1998). Notices Amer.
Math. Soc. 46 (1999), no. 3, 338--343.

\section*{2000s}

\hspace{-2em} 2000 Tate, John. The millennium prize problems I, Lecture by
John Tate at the Millenium Meeting of the Clay Mathematical Institute, May
2000, Paris. Video available from the CMI website.

2001 Tate, John. Galois cohomology. Arithmetic algebraic geometry (Park City,
UT, 1999), 465--479, IAS/Park City Math. Ser., 9, Amer. Math. Soc.,
Providence, RI, 2001.

2002 Tate, John. On a conjecture of Finotti. Bull. Braz. Math. Soc. (N.S.) 33
(2002), no. 2, 225--229.

2004 Tate, John. Refining Gross's conjecture on the values of abelian
L-functions. Stark's conjectures: recent work and new directions, 189--192,
Contemp. Math., 358, Amer. Math. Soc., Providence, RI, 2004.

2005 Artin, Michael; Rodriguez-Villegas, Fernando; Tate, John. On the
Jacobians of plane cubics. Adv. Math. 198 (2005), no. 1, 366--382.

2006 Mazur, Barry; Stein, William; Tate, John. Computation of p-adic heights
and log convergence. Doc. Math. 2006, Extra Vol., 577--614 (electronic).

2008 Tate, John. Foreword to $p$-adic geometry, 9--63, Univ. Lecture Ser., 45, Amer. Math. Soc., Providence, RI, 2008.

2009 Artin, Emil; Tate, John. Class field theory. New edition. TeXed and
slightly revised from the original 1961 version. AMS Chelsea Publishing,
Providence, RI, 2009. viii+194 pp.

\section*{2010s}

\hspace{-2em}2011 Raussen, Martin; Skau, Christian. Interview with Abel Laureate John Tate. Notices Amer. Math. Soc. 58 (2011), no. 3, 444--452.

2011 Tate, John. Stark's basic conjecture. Arithmetic of $L$-functions, 7--31, IAS/Park City Math. Ser., 18, Amer. Math. Soc., Providence, RI.

\clearpage \addcontentsline{toc}{section}{Index} \printindex
\end{document}